\documentclass{amsart}

\usepackage{amsmath}
\usepackage{amsfonts}
\usepackage{amssymb}
\usepackage{amsthm}
\usepackage[dpi=600,PostScript=dvips]{diagrams}

\makeatletter
\def\@secnumfont{\relax}
\makeatother

\newcommand{\al}{\alpha}
\newcommand{\be}{\beta}
\newcommand{\ga}{\gamma}
\newcommand{\ali}{{\alpha^{-1}}}
\newcommand{\bei}{{\beta^{-1}}}
\newcommand{\gai}{{\gamma^{-1}}}

\DeclareMathOperator{\id}{id}

\newcommand{\cp}[2]{\Delta_{#1,#2}}
\newcommand{\cpcop}[2]{\Delta^{\rm cop}_{#1,#2}}
\newcommand{\la}{\lambda_{\alpha}}
\newcommand{\lb}{\lambda_{\beta}}
\newcommand{\lai}{\lambda_{\alpha^{-1}}}

\newcommand{\lab}{\lambda_{\alpha \beta}}
\newcommand{\lbabi}{\lambda_{\beta \alpha \beta^{-1}}}
\newcommand{\lu}{\lambda_1}

\newcommand{\pa}{\varphi_{\alpha}}
\newcommand{\pai}{\varphi_{\alpha^{-1}}}
\newcommand{\pb}{\varphi_{\beta}}
\newcommand{\pbi}{\varphi_{\beta^{-1}}}
\newcommand{\pab}{\varphi_{\alpha \beta}}
\newcommand{\pabi}{\varphi_{(\al \be)^{-1}}}
\newcommand{\ta}{\theta_{\alpha}}
\newcommand{\tai}{\theta_{\alpha^{-1}}}
\newcommand{\tb}{\theta_{\beta}}

\newcommand{\rt}{\widetilde{R}}

\newcommand{\rtab}{\widetilde{R}_{\al,\be}}

\newcommand{\rab}{R_{\al,\be}}
\newcommand{\ruu}{R_{1,1}}

\newcommand{\raai}{R_{\al,\ali}}

\newcommand{\rba}{R_{\be,\al}}
\newcommand{\rbai}{R_{\be,\ali}}

\newcommand{\tba}{\sigma_{\be,\al}}

\newcommand{\tuu}{\sigma_{1,1}}
\newcommand{\rh}[2]{\rho_{#1,#2}}

\newcommand{\aua}{a_{(1,\alpha)}}
\newcommand{\ada}{a_{(2,\alpha)}}
\newcommand{\ata}{a_{(3,\alpha)}}

\newcommand{\adb}{a_{(2,\beta)}}
\newcommand{\atb}{a_{(3,\beta)}}

\newcommand{\adai}{a_{(2,\alpha^{-1})}}

\newcommand{\adbi}{a_{(2,\beta^{-1})}}

\newcommand{\auu}{a_{(1,1)}}
\newcommand{\adu}{a_{(2,1)}}

\newcommand{\Lua}{\Lambda_{(1,\alpha)}}

\newcommand{\Ldai}{\Lambda_{(2,\alpha^{-1})}}

\newcommand{\hua}{h_{(1,\alpha)}}
\newcommand{\hda}{h_{(2,\alpha)}}
\newcommand{\hta}{h_{(3,\alpha)}}

\newcommand{\hdb}{h_{(2,\beta)}}

\newcommand{\hdai}{h_{(2,\alpha^{-1})}}

\newcommand{\huu}{h_{(1,1)}}
\newcommand{\hdu}{h_{(2,1)}}

\newcommand{\cua}{c_{(1,\alpha)}}

\newcommand{\cdb}{c_{(2,\beta)}}

\newcommand{\ctg}{c_{(3,\gamma)}}

\newcommand{\mua}{m_{(1,\alpha)}}

\newcommand{\mub}{m_{(1,\beta)}}
\newcommand{\mdb}{m_{(2,\beta)}}

\newcommand{\mdg}{m_{(2,\gamma)}}

\newcommand{\muai}{m_{(1,\alpha^{-1})}}
\newcommand{\mdai}{m_{(2,\alpha^{-1})}}
\newcommand{\mtai}{m_{(3,\alpha^{-1})}}

\newcommand{\mdbi}{m_{(2,\beta^{-1})}}

\newcommand{\muu}{m_{(1,1)}}

\newcommand{\mza}{m_{(0,\alpha)}}

\newcommand{\mzai}{m_{(0,\alpha^{-1})}}

\newcommand{\mzab}{m_{(0,\alpha \beta)}}

\newcommand{\muabi}{m_{(1,(\alpha \beta)^{-1})}}

\newcommand{\xua}{x_{(1,\alpha)}}
\newcommand{\xda}{x_{(2,\alpha)}}
\newcommand{\xta}{x_{(3,\alpha)}}

\newcommand{\xub}{x_{(1,\beta)}}
\newcommand{\xdb}{x_{(2,\beta)}}

\newcommand{\xuai}{x_{(1,\alpha^{-1})}}
\newcommand{\xdai}{x_{(2,\alpha^{-1})}}
\newcommand{\xtai}{x_{(3,\alpha^{-1})}}

\newcommand{\xdbi}{x_{(2,\beta^{-1})}}

\newcommand{\xuu}{x_{(1,1)}}
\newcommand{\xdu}{x_{(2,1)}}

\newcommand{\xza}{x_{(0,\alpha)}}

\newcommand{\xzu}{x_{(0,1)}}

\newcommand{\yua}{y_{(1,\alpha)}}

\newcommand{\ydb}{y_{(2,\beta)}}

\newcommand{\nda}{n_{(2,\alpha)}}

\newcommand{\nuai}{n_{(1,\alpha^{-1})}}

\newcommand{\nuu}{n_{(1,1)}}

\newcommand{\nza}{n_{(0,\alpha)}}

\newcommand{\oM}{{\overline{M}}}
\newcommand{\oF}{{\overline{F}}}
\newcommand{\oN}{{\overline{N}}}
\newcommand{\oL}{{\overline{L}}}
\newcommand{\of}{{\overline{f}}}

\newcommand{\vrh}[2]{\varrho_{#1,#2}}
\newcommand{\oHs}{{\kern.2ex \overline{\kern-0.2ex H^*}}}
\newcommand{\oHsa}{{\overline{H^*_\al}}}
\newcommand{\pc}{{\scriptscriptstyle \square}}
\newcommand{\HH}{H^\pc}

\newcommand{\kk}{\Bbbk}
\newcommand{\cop}{{\rm cop}}
\newcommand{\co}{{\rm co}}
\newcommand{\opp}{{\rm op}}
\newcommand{\opcop}{{\rm op,cop}}
\newcommand{\rat}{{\rm rat}}
\DeclareMathOperator{\homo}{Hom}  \DeclareMathOperator{\conv}{Conv}
\DeclareMathOperator{\tr}{tr}

\newcommand{\p}{$\pi$\nobreakdash}

\theoremstyle{plain}
\newtheorem{theorem}{Theorem}
\newtheorem{lemma}{Lemma}

\newtheorem{corollary}{Corollary}

\theoremstyle{definition}

\newtheorem{example}{Example}
\newtheorem*{notation}{Notation}
\newtheorem*{SWnotation}{Sweedler's notation}

\theoremstyle{remark}
\newtheorem*{remark}{Remark}

\numberwithin{equation}{subsection}
\renewcommand{\theequation}{\thesubsection.\alph{equation}}

\newcommand{\libelle}[1]{(#1)}
\newcounter{pastempo}

\newenvironment{defi}[1][8.8.d]
 {\setcounter{pastempo}{\value{equation}}
  \begin{list}{}{\settowidth{\labelwidth}{\libelle{#1}}
                 \settowidth{\labelsep}{\ \,}
                 \setlength{\leftmargin}{\labelwidth}
                 \addtolength{\leftmargin}{\labelsep}
                 \usecounter{equation}
                 \setcounter{equation}{\value{pastempo}}
                 }}
 {\end{list}}

\newarrow{Inclu}{C}{-}{-}{-}{>}
\newarrow{DashInclu}{C}{dash}{}{dash}{>}
\newarrow{DotLine}{.}{.}{.}{.}{.}
\newarrow{Line}{-}{-}{-}{-}{-}
\newarrow{ToDash}{-}{-}{}{.}{.}
\newarrow{DashTo}{dash}{dash}{}{-}{>}
\newarrow{DashIncluLine}{C}{dash}{}{dash}{}
\newarrow{DashLine}{}{dash}{}{dash}{}

\begin{document}
\title{Hopf group-coalgebras}
\author[A. Virelizier]{Alexis Virelizier}
\address{Institut de Recherche Math\'ematique Avanc\'ee - Universit\'{e} Louis Pasteur - C.N.R.S.\\
7 rue Ren\'{e} Descartes \\ 67084 Strasbourg Cedex, France.} \email{virelizi@math.u-strasbg.fr}
\date{December 11, 2000}

\begin{abstract}
We study algebraic properties of Hopf group-coalgebras, recently introduced by Turaev. We show the existence of
integrals and traces for such coalgebras, and we generalize the main properties of quasitriangular and ribbon Hopf
algebras to the setting of Hopf group-coalgebras.
\end{abstract}
\maketitle

\setcounter{tocdepth}{1} \tableofcontents

\section*{Introduction}

Recently, Turaev \cite{Tur1} introduced, for a group $\pi$, the notion of a modular crossed \p-category and showed
that such a category gives rise to a $3$-dimensional homotopy quantum field theory with target space $K(\pi,1)$.
Examples of \p-categories can be constructed from so-called Hopf \p-coalgebras defined in \cite{Tur1}.

The notion of a Hopf \p-coalgebra generalizes that of a Hopf algebra. Hopf \p-coalgebras are used by the author in
\cite{Vir} to construct Hennings-like (see \cite{He2,KR1}) and Kuperberg-like (see \cite{Ku1}) invariants of
principal \p-bundles over link complements and over $3$-manifolds. The aim of the present paper is to lay the
algebraic foundations for \cite{Vir}, specifically the existence of integrals and traces for a Hopf \p-coalgebra.

Let us briefly recall some definitions of \cite{Tur1}. Given a (discrete) group $\pi$, a Hopf \p-coalgebra is a
family $H=\{H_\al\}_{\al \in \pi}$ of algebras (over a field $\Bbbk$) endowed with a comultiplication
$\Delta=\{\cp{\al}{\be}: H_{\al \be} \to H_\al \otimes H_\be\}_{\al,\be \in \pi}$, a counit $\epsilon: H_1 \to
\kk$, and an antipode $S=\{S_\al:H_\al \to H_\ali\}_{\al \in \pi}$ which verify some compatibility conditions. A
crossing for $H$ is a family of algebra isomorphisms $\varphi=\{\varphi_\be : H_\al \to H_{\be \al \bei}
\}_{\al,\be \in \pi}$ which preserves the comultiplication and the counit, and which yields an action of $\pi$ in
the sense that $\pb \varphi_{\be'} = \varphi_{\be \be'}$. A crossed Hopf \p-coalgebra $H$ is quasitriangular
(resp.\@ ribbon) when it is endowed with an $R$-matrix $R=\{R_{\al,\be} \in H_\al \otimes H_\be \}_{\al,\be \in
\pi}$ (resp.\@ an $R$-matrix and a twist $\theta=\{\theta_\al \in H_\al\}_{\al \in \pi}$) verifying some axioms
which generalize the classical ones given in \cite{Drin2} (resp.\@ \cite{RT}). The case $\pi=1$ is the standard
setting of Hopf algebras. When $\pi$ is commutative and $\varphi$ is trivial, one recovers the definition of a
quasitriangular or ribbon \p-colored Hopf algebra given by Ohtsuki~\cite{Oh}.

Basic notions of the theory of Hopf algebras can be extended to the setting of Hopf \p-coalgebras. In particular,
a (right) \p-integral for a Hopf \p-coalgebra $H$ is a family of $\kk$-forms $\lambda=\{\la:H_\al\to \kk\}_{\al
\in \pi}$ such that $(\la \otimes \id_{H_\be}) \cp{\al}{\be} = \lambda_{\al\be} \, 1_\be$ for all $\al, \be \in
\pi$. When $H$ is crossed, a \p-trace for $H$ is a family of $\kk$-forms $\tr=\{\tr_\al:H_\al \to \kk\}_{\al \in
\pi}$ which verifies $\tr_\al (xy)= \tr_\al (yx)$, $\tr_\ali (S_\al (x)) = \tr_\al (x)$, and $\tr_{\be \al
\bei}(\pb(x))=\tr_\al(x)$ for all $\al,\be \in \pi$ and $x,y \in H_\al$. These notions were introduced in
\cite{Vir} for topological purposes.

In the first part of the paper (Sect.\@ \ref{basicdefi}-\ref{semicosemi}), we mainly focus on finite dimensional
Hopf \p-coalgebras, that is Hopf \p-coalgebras $H=\{H_\al\}_{\al \in \pi}$ with each $H_\al$ finite dimensional.
The first main result is the existence and uniqueness (up to a scalar multiple) of a \p-integral for such a Hopf
\p-coalgebra. To prove this result, we study rational \p-graded modules, introduce the notion of a Hopf
\p-comodule, and generalize the fundamental theorem of Hopf modules (see \cite{LS}) to Hopf \p-comodules.

As for Hopf algebras, any finite dimensional Hopf \p-coalgebra contains a distinguished \p-grouplike element.
Generalizing \cite{Rad1}, we study the relationships between this element, the antipode, and the \p-integrals. As
a corollary, we give an upper bound for the order of $S_\ali S_\al$ whenever $\al \in \pi$ has a finite order.

The notions of semisimplicity and cosemisimplicity can be extended to the setting of Hopf \p-coalgebras. We show
that a finite dimensional Hopf \p-coalgebra $H=\{H_\al\}_{\al \in \pi}$ is semisimple (that is each $H_\al$ is
semisimple) if and only if $H_1$ is semisimple. We define the cosemisimplicity for \p-comodules and \p-coalgebras,
and we use \p-integrals to give necessary and sufficient criteria for a Hopf \p-coalgebra to be cosemisimple.

In the second part of the paper (Sect.\@ \ref{quasitriangularity}-\ref{s:pitrace}), we study quasitriangular Hopf
\p-coalgebras. The main result is the existence of \p-traces for a semisimple (resp. cosemisimple) finite
dimensional unimodular ribbon Hopf \p-coalgebra. To prove this result, we generalize the main properties of
quasitriangular Hopf algebras (see~\cite{Drin,Rad2}). In particular, we introduce and study the (generalized)
Drinfeld elements of a quasitriangular Hopf \p-coalgebra $H$, we compute the distinguished \p-grouplike element of
$H$ by using the $R$-matrix, and we show that the twist of a ribbon Hopf \p-coalgebra leads to a \p-grouplike
element which implements the square of the antipode by conjugation.

The paper is organized as follows. In Section~\ref{basicdefi}, we review the basic definitions and properties of
Hopf \p-coalgebras. In Section~\ref{modcomod}, we discuss the notions of a rational \p-graded module and of a Hopf
\p-comodule. In Section~\ref{piint}, we use these notions to establish the existence and uniqueness of
\p-integrals. Section~\ref{s:disting} is devoted to the study of the distinguished \p-grouplike element. In
Section~\ref{semicosemi}, we discuss the notion of a semisimple (resp.\@ cosemisimple) Hopf \p-coalgebra. In
Section~\ref{quasitriangularity}, we study crossed, quasitriangular, and ribbon Hopf \p-coalgebras. Finally, we
construct \p-traces in Section~\ref{s:pitrace}.\\

{\it Acknowledgements.} The author would like to thank his advisor Vladimir Turaev for his useful suggestions and
constructive criticism.

\section{Basic definitions}\label{basicdefi}
Throughout the paper, we let $\pi$ be a discrete group (with neutral element 1) and $\Bbbk$ be a field (although
much of what we do is valid over any commutative ring). We set $\kk^*=\kk \setminus\{0\}$. All algebras are
supposed to be over $\Bbbk$, associative, and unitary. The tensor product $\otimes=\otimes_\Bbbk$ is always
assumed to be over $\Bbbk$. If $U$ and $V$ are $\Bbbk$-spaces, $\sigma_{U,V}:U \otimes V \to V \otimes U$ will
denote the flip defined by $\sigma_{U,V}(u \otimes v)=v \otimes u$.

\subsection{$\pi$-coalgebras}\label{coalgebra}
We recall the definition of a \p-coalgebra, following \cite[\S 11.2]{Tur1}. A \emph{\p-coalgebra} (over $\Bbbk$)
is a family $C=\{C_\al\}_{\al \in \pi}$ of $\Bbbk$-spaces endowed with a family $\Delta=\{\cp{\al}{\be}:C_{\al
\be} \to C_\al \otimes C_\be \}_{\al,\be \in \pi}$ of $\Bbbk$-linear maps (the {\it comultiplication}) and a
$\Bbbk$-linear map $\epsilon:C_1 \to \Bbbk$ (the \emph{counit}) such that
\begin{defi}
  \item \label{coass} $\Delta$ is coassociative in the sense that, for any $\al,\be,\ga \in \pi$,
        $$ (\cp{\al}{\be} \otimes \id_{C_\ga})
        \cp{\al \be}{\ga}=(\id_{C_\al} \otimes
        \cp{\be}{\ga}) \cp{\al}{\be \ga};$$
  \item \label{counit} for all $\al \in \pi$,
        $$ (\id_{C_\al} \otimes \epsilon)
        \cp{\al}{1}=\id_{C_\al}=(\epsilon \otimes \id_{C_\al}) \cp{1}{\al}.$$
\end{defi}

Note that $(C_1,\cp{1}{1},\epsilon)$ is a coalgebra in the usual sense of the word.

\begin{SWnotation}
We extend the Sweedler notation for a comultiplication in the following way: for any $\al,\be \in \pi$ and $c \in
C_{\al \be}$, we write
  $$
  \cp{\al}{\be}(c)=\sum_{(c)} \cua \otimes \cdb \in C_\al \otimes C_\be,
  $$
or shortly, if we leave the summation implicit,
  $
  \cp{\al}{\be}(c)= \cua \otimes \cdb
  $.\\
The coassociativity axiom~\eqref{coass} gives that, for any $\al,\be, \ga \in \pi$ and $c \in C_{\al \be \ga}$,
  $$
  c_{(1,\al \be)(1,\al)} \otimes c_{(1,\al \be)(2,\be)} \otimes
  c_{(2,\ga)} = c_{(1,\al)} \otimes c_{(2,\be \ga)(1,\be)} \otimes c_{(2,\be \ga)(2,\ga)}.
  $$
This element of $C_\al \otimes C_\be \otimes C_\ga$ is written as $\cua \otimes \cdb \otimes \ctg$. For any $c \in
C_{\al_1 \dotsm \al_n}$, by iterating the
procedure, we define inductively $c_{(1,\al_1)} \otimes \dotsb \otimes c_{(n,\al_n)}$.
\end{SWnotation}

\subsection{Convolution algebras}\label{convo}
Let $C=(\{C_\al\},\Delta,\epsilon)$ be a \p-coalgebra and $A$ be an algebra with multiplication $m$ and unit
element $1_A$. For any $f \in \homo_\kk(C_\al,A)$ and $g \in \homo_\kk(C_\be,A)$, we define their
\emph{convolution product} by
 $$
 f*g=m(f \otimes g) \cp{\al}{\be} \in \homo_\kk(C_{\al \be},A).
 $$
Using~\eqref{coass} and~\eqref{counit}, one verifies that the $\Bbbk$-space
 $$
 \conv(C,A)=\oplus_{\al \in \pi}\homo_\kk(C_\al,A),
 $$
endowed with the convolution product $*$ and the unit element $\epsilon 1_A$, is a \p-graded algebra, called
\emph{convolution algebra}.

In particular, for $A=\kk$, the \p-graded algebra $\conv(C,\kk)=\oplus_{\al \in \pi}C^*_\al$ is called \emph{dual}
to $C$ and is denoted by $C^*$.

\subsection{Hopf $\pi$-coalgebras}\label{hopfpicoal}
Following \cite[\S 11.2]{Tur1}, a \emph{Hopf \p-coalgebra} is a \p-coalgebra $H=(\{H_\al\},\Delta,\epsilon)$
endowed with a family $S=\{S_\al : H_\al \to H_\ali \}_{\al \in \pi}$ of $\Bbbk$-linear maps (the \emph{antipode})
such that
\begin{defi}
  \item \label{hopfalg}  each $H_\al$ is an algebra with multiplication
                         $m_\al$
                         and unit element $1_\al \in H_\al$;
  \item \label{hopfmor}  $\epsilon:H_1 \to \kk$ and
                         $\cp{\al}{\be}: H_{\al \be} \to H_\al \otimes H_\be$
                         (for all $\al,\be \in \pi$) are algebra homomorphisms;
  \item \label{antipode}  for any $\al \in \pi$,
                          $$
                          m_\al (S_\ali \otimes \id_{H_\al}) \cp{\ali}{\al} =
                          \epsilon 1_\al =  m_\al (\id_{H_\al} \otimes S_\ali)
                          \cp{\al}{\ali} .
                          $$
\end{defi}

Remark that it is not a self-dual notion and that $(H_1,m_1,1_1,\cp{1}{1},\epsilon,S_1)$ is a (classical) Hopf
algebra.

The Hopf \p-coalgebra $H$ is said to be \emph{finite dimensional} if, for all $\al \in \pi$, $H_\al$ is finite
dimensional (over $\Bbbk$). Note that it does not mean that $\oplus_{\al \in \pi} H_\al$ is finite dimensional
(unless $H_\al=0$ for all but a finite number of $\al \in \pi$).

The antipode $S=\{ S_\al \}_{\al \in \pi}$ of $H$ is said to be \emph{bijective} if each $S_\al$ is bijective.
Unlike \cite[\S 11.2]{Tur1}, we do not suppose that the antipode of a Hopf \p-coalgebra $H$ is bijective. However,
we will show that it is bijective whenever  $H$ is finite dimensional (see Corollary~\ref{corinteg}(a)) or
quasitriangular (see Lemma~\ref{B2}(c)).

A useful remark is that if $H=\{H_\al\}_{\al \in \pi}$ is a Hopf \p-coalgebra with antipode $S=\{S_\al\}_{\al \in
\pi}$, then Axiom~\eqref{antipode} says that $S_\al$ is the inverse of $\id_{H_\ali}$ in the convolution algebra
$\conv(H,H_\ali)$ for all $\al \in \pi$.

In the next lemma, generalizing \cite[Proposition~4.0.1]{sweed}, we show that the antipode of a Hopf \p-coalgebra
is anti-multiplicative and anti-comultiplicative.
\begin{lemma}\label{antipodepptes}
Let $H=(\{H_\al,m_\al,1_\al\},\Delta,\epsilon,S)$ be a Hopf \p-coalgebra. Then
\begin{enumerate}
\renewcommand{\labelenumi}{{\rm (\alph{enumi})}}
 \item $S_\al(ab)=S_\al(b) S_\al(a)$ for any $\al \in \pi$ and $a,b \in H_\al$;
 \item $S_\al(1_\al)=1_\ali$ for any $\al \in \pi$;
 \item $\cp{\bei}{\ali} S_{\al \be}=\sigma_{H_\ali,H_\bei}(S_\al
       \otimes S_\be) \cp{\al}{\be}$ for any $\al,\be \in \pi$;
 \item $\epsilon S_1=\epsilon$.
\end{enumerate}
\end{lemma}

\begin{proof}
The proof is essentially the same as in the Hopf algebra setting. For example, to show Part (c), fix $\al,\be \in
\pi$ and consider the algebra $\conv(H, H_\bei \otimes H_\ali)$ with convolution product $*$ and unit element
$e=\epsilon 1_\bei \otimes 1_\ali$. Using Axioms \eqref{counit}, \eqref{hopfmor}, and \eqref{antipode}, one easily
checks that $\cp{\bei}{\ali} S_{\al \be} * \cp{\bei}{\ali} = e$ and $\cp{\bei}{\ali} *
\sigma_{H_\ali,H_\bei}(S_\al \otimes S_\be) \cp{\al}{\be}=e$. Hence we can conclude that $\cp{\bei}{\ali} S_{\al
\be} =\sigma_{H_\ali,H_\bei}(S_\al \otimes S_\be) \cp{\al}{\be}$.
\end{proof}

\begin{corollary}\label{subgroup}
Let $H=\{H_\al\}_{\al \in \pi}$ be a Hopf \p-coalgebra. Then $\{ \al \in \pi \, | \, H_\al \neq 0 \}$ is a
subgroup of $\pi$.
\end{corollary}
\begin{proof}
Set $G=\{ \al \in \pi \, | \, H_\al \neq 0 \}$. Firstly $1_1\neq 0$ (since $\epsilon(1_1)=1 \neq 0$) and so $1 \in
G$. Then let $\al,\be \in G$. Using \eqref{hopfmor}, $\cp{\al}{\be}(1_{\al \be})=1_\al \otimes 1_\be \neq 0$.
Therefore $1_{\al \be} \neq 0$ and so $\al \be \in G$. Finally, let $\al \in G$. By Lemma~\ref{antipodepptes}(b),
$S_\ali(1_\ali)=1_\al \neq 0$. Thus $1_\ali \neq 0 $ and hence $\ali \in G$.
\end{proof}

\subsubsection{Opposite Hopf $\pi$-coalgebra}\label{op}
Let $H=\{H_\al\}_{\al \in \pi}$ be a Hopf \p-coalgebra. Suppose that the antipode $S=\{ S_\al \}_{\al \in \pi}$ of
$H$ is bijective. For any $\al \in \pi$, let $H^\opp_\al$ be the opposite algebra to $H_\al$. Then
$H^\opp=\{H^\opp_\al\}_{\al \in \pi}$, endowed with the comultiplication and counit of $H$ and with the antipode
$S^\opp=\{ S^\opp_\al=S^{-1}_\ali \}_{\al \in \pi}$, is a Hopf \p-coalgebra, called \emph{opposite to} $H$.

\subsubsection{Coopposite Hopf $\pi$-coalgebra}\label{coop}
Let $C=(\{C_\al\},\Delta,\epsilon)$ be a \p-coalgebra. Set
 $$
 C^\cop_\al=C_\ali \quad \text{and} \quad
 \cpcop{\al}{\be}=\sigma_{C_\bei,C_\ali} \cp{\bei}{\ali}.
 $$
Then $C^\cop=(\{C^\cop_\al\},\Delta^\cop,\epsilon)$ is a \p-coalgebra, called \emph{coopposite to} $C$. If $H$ is
a Hopf \p-coalgebra whose antipode $S=\{ S_\al \}_{\al \in \pi}$ is bijective, then the coopposite \p-coalgebra
$H^\cop$, where $H^\cop_\al=H_\ali$ as an algebra, is a Hopf \p-coalgebra with antipode $S^\cop=\{
S^\cop_\al=S^{-1}_\al\}_{\al \in \pi}$.

\subsubsection{Opposite and coopposite Hopf $\pi$-coalgebra}\label{opcoop} Let $H=(\{H_\al\},\Delta,\epsilon,S)$ be a
Hopf \p-coalgebra. Even if the antipode of $H$ is not bijective, one can always define a Hopf \p-coalgebra
\emph{opposite and coopposite to} $H$ by setting $H^\opcop_\al=H_\ali^\opp$,
$\cp{\al}{\be}^\opcop=\cpcop{\al}{\be}$, $\epsilon^\opcop=\epsilon$, and $S^\opcop_\al=S_\ali$.
\subsubsection{The dual Hopf algebra}\label{dualhopf}
Let $H=(\{H_\al,m_\al,1_\al\}, \Delta, \epsilon, S)$ be a finite dimensional Hopf \p-coalgebra. The \p-graded
algebra $H^*=\oplus_{\al \in \pi} H_\al^*$ dual to $H$ (see \S\ref{convo}) inherits a structure of a Hopf algebra
by setting, for all $\al \in \pi$ and $f \in H_\al^*$,
  $$
  \Delta(f)= m_\al^*(f) \in (H_\al
  \otimes H_\al)^* \cong H_\al^* \otimes H_\al^*,
  $$
$\epsilon(f)=f(1_\al)$, and $S(f)=f \circ S_\ali$. Note that if $H_\al \neq 0$ for infinitely many $\al \in \pi$,
then $H^*$ is infinite dimensional.

\subsubsection{The case $\pi$ finite}\label{pifinite}
Let us first remark that, when $\pi$ is a finite group, there is a one-to-one correspondence between (isomorphic
classes of) \p-coalgebras and (isomorphic classes of) \p-graded coalgebras. Recall that a coalgebra
$(C,\Delta,\epsilon)$ is \emph{\p-graded} if $C$ admits a decomposition as a direct sum of $\kk$-spaces
$C=\oplus_{\al \in \pi} C_\al$ such that, for any $\al \in \pi$,
  $$
  \Delta(C_\al)\subseteq  \sum_{\be \ga=\al} C_\be \otimes C_\ga
  \quad \text{and} \quad \epsilon(C_\al)=0 \;\, \text{if $\al \neq 1$}.
  $$
Let us denote by $p_\al:C \to C_\al$ the canonical projection. Then $\{C_\al\}_{\al \in \pi}$ is a \p-coalgebra
with comultiplication $\{(p_\al \otimes p_\be) \Delta_{|C_{\al\be}} \}_{\al, \be \in \pi}$ and counit
$\epsilon_{|C_1}$. Conversely, if $C=(\{C_\al\},\Delta,\epsilon)$ is a \p-coalgebra, then $\tilde{C}=\oplus_{\al
\in \pi} C_\al$ is a \p-graded coalgebra with comultiplication $\tilde{\Delta}$ and counit $\tilde{\epsilon}$
given on the summands by
  $$
  \tilde{\Delta}_{|C_\al}=\sum_{\be \ga = \al} \cp{\be}{\ga} \quad \text{and} \quad
  \tilde{\epsilon}_{|C_\al}=\begin{cases}
                               \epsilon & \text{if $\al =1$}\\
                               0 & \text{if $\al \neq 1$}
                            \end{cases}.
  $$
Let now $H=(\{H_\al,m_\al,1_\al\}, \Delta, \epsilon, S)$ be a Hopf \p-coalgebra, where $\pi$ is a finite group.
Then the coalgebra $(\tilde{H},\tilde{\Delta},\tilde{\epsilon})$, defined as above, is a Hopf algebra with
multiplication $\tilde{m}$, unit element $\tilde{1}$, and antipode $\tilde{S}$ given by
  $$
  \tilde{m}_{|H_\al \otimes H_\be}=\begin{cases}
                                      m_\al & \text{if $\al = \be$} \\
                                      0 & \text{if $\al \neq \be$}
                                   \end{cases},
  \quad \tilde{1}= \sum_{\al \in \pi} 1_\al, \quad \text{and} \quad
  \tilde{S}=\sum_{\al \in \pi} S_\al.
  $$

When $H$ is finite dimensional and $\pi$ is finite, the Hopf algebra $H^*$ (see \S\ref{dualhopf}) is simply the
dual Hopf algebra $\tilde{H}^*$.

\begin{remark}
When $\pi$ is finite, the structure of \p-comodules over a \p-coalgebra $C$ (Theorem~\ref{theorational}), the
existence of \p-integrals for a finite dimensional Hopf \p-coalgebra $H$ (Theorem~\ref{piintegexists}) and their
relations with the distinguished \p-grouplike element (Theorem~\ref{intfortrace}) can be deduced from the
classical theory of coalgebras and Hopf algebras by using $\tilde{C}$ or $\tilde{H}$ (defined as in
\S\ref{pifinite}). Nevertheless, for the general case, self-contained proofs must be given.

In general, the results relating to a quasitriangular Hopf \p-coalgebra (see Sect.\@
\ref{quasitriangularity}-\ref{s:pitrace}) cannot be deduced from the classical theory of quasitriangular Hopf
algebras. Indeed, even if $\pi$ is finite, an $R$-matrix for a Hopf $\pi$-coalgebra $H$ (whose definition involves
an action of $\pi$, see \S\ref{quasitrig}) does not necessarily lead to a usual $R$-matrix for the Hopf algebra
$\tilde{H}$.
\end{remark}

\section{Modules and comodules}\label{modcomod}
In this section, we introduce and discuss the notions of $\pi$-comodules, rational $\pi$-graded modules, and Hopf
$\pi$-comodules. They are used in Section~\ref{piint} to show the existence of integrals.

\subsection{$\pi$-comodules}\label{picomod}
Let  $C=(\{C_\al\},\Delta,\epsilon)$ be a \p-coalgebra. A \emph{right \p-comodule over} $C$ is a family
$M=\{M_\al\}_{\al \in \pi}$ of $\Bbbk$-spaces endowed with a family $\rho=\{\rh{\al}{\be}:M_{\al \be} \to M_\al
\otimes C_\be \}_{\al, \be \in \pi}$  of $\Bbbk$-linear maps (the \emph{structure maps}) such that
\begin{defi}
  \item \label{coasscomod} for any $\al,\be,\ga \in \pi$,
                           $$
                           (\rh{\al}{\be} \otimes \id_{C_\ga}) \rh{\al
                           \be}{\ga}=(\id_{M_\al} \otimes
                           \cp{\be}{\ga}) \rh{\al}{\be \ga};
                           $$
  \item \label{counitcomod} for any $\al \in \pi$,
                            $$
                            (\id_{M_\al} \otimes \epsilon)
                            \rh{\al}{1}=\id_{M_\al}.
                            $$
\end{defi}

Note that $M_1$ endowed with the structure map $\rh{1}{1}$ is a (usual) right $C_1$-comodule.

If $\pi$ is finite and $\tilde{C}=\oplus_{\al\in \pi} C_\al$ is the \p-graded coalgebra defined as in
\S\ref{pifinite}, then $M$ leads to a $\pi$-graded comodule $\tilde{M}=\oplus_{\al \in \pi} M_\al$ over
$\tilde{C}$ with comodule map $\tilde{\rho}=\sum_{\al, \be \in \pi} \rh{\al}{\be}$ (see \cite{NT1}).

A \emph{\p-subcomodule of $M$} is a family $N=\{N_\al\}_{\al \in \pi}$, where $N_\al$ is a $\Bbbk$-subspace of
$M_\al$, such that $\rh{\al}{\be}(N_{\al \be}) \subset N_\al \otimes C_\be$ for all $\al,\be \in \pi$. Then $N$ is
a right \p-comodule over $C$ with induced structure maps.

A \emph{\p-comodule morphism} between two right \p-comodules $M$ and $M'$ over $C$ (with structure maps $\rho$ and
$\rho'$) is a family $f=\{f_\al : M_\al \to M'_\al\}_{\al \in \pi}$ of $\Bbbk$-linear maps such that
$\rh{\al}{\be}' f_{\al \be}=(f_\al \otimes \id_{C_\be}) \rh{\al}{\be}$ for all $\al,\be \in \pi$.

\begin{SWnotation}
We extend the notation of Section~\ref{coalgebra} by setting, for any $\al,\be \in \pi$ and $m \in M_{\al \be}$,
  $$
  \rh{\al}{\be}(m)=\mza \otimes \mub \in M_\al \otimes C_\be.
  $$
Axiom \eqref{coasscomod} gives that, for any $\al,\be, \ga \in \pi$ and $m \in M_{\al\be\ga}$,
  $$
  m_{(0,\al \be)(0,\al)} \otimes m_{(0,\al \be)(1,\be)} \otimes m_{(1,\ga)}
  = m_{(0,\al)} \otimes m_{(1,\be \ga)(1,\be)} \otimes m_{(1,\be \ga)(2,\ga)}.
  $$
This element of $M_\al \otimes C_\be \otimes C_\ga$ is written as $ \mza \otimes \mub \otimes \mdg$. For any $m
\in M_{\al_0 \al_1 \cdots \al_n}$,   we define inductively $m_{(0,\al_0)} \otimes m_{(1,\al_1)} \otimes \cdots
\otimes m_{(n,\al_n)}$ by iterating the procedure.
\end{SWnotation}

Let $N=\{N_\al\}_{\al \in \pi}$ be a \p-subcomodule of a right \p-comodule $M=\{M_\al\}_{\al \in \pi}$ over a
\p-coalgebra $C$. One easily checks that $M/N=\{M_\al/N_\al\}_{\al \in \pi}$ is a right \p-comodule over $C$, with
structure maps naturally induced from the structure maps of $M$. Moreover this is the unique structure of a right
\p-comodule over $C$ on $M/N$ which makes the canonical projection $p=\{ p_\al : M_\al \to M_\al/N_\al \}_{\al \in
\pi}$ a \p-comodule morphism.

If $f=\{f_\al:M_\al \to M_\al' \}_{\al \in \pi}$ is a \p-comodule morphism between two right \p-comodules $M$ and
$M'$, then $\ker(f)=\{ \ker(f_\al) \}_{\al \in \pi}$ is a \p-subcomodule of $M$, $f(M)=\{ f_\al(M_\al) \}_{\al \in
\pi}$ is a \p-subcomodule of $M'$, and the canonical isomorphism $\Bar{f}=\{\Bar{f}_\al:M_\al/\ker(f_\al) \to
f_\al(M_\al) \}_{\al \in \pi}$ is a \p-comodule isomorphism.

\begin{example}\label{coinvariants}
Let $H$ be a Hopf \p-coalgebra and $M=\{ M_\al \}_{\al \in \pi}$ be a right \p-comodule over $H$ with structure
maps $\rho=\{\rh{\al}{\be}\}_{\al,\be \in \pi}$. The \emph{coinvariants of $H$ on $M$} are the elements of the
$\Bbbk$-space
 $$
 \{ m=(m_\al)_{\al \in \pi} \in \underset{\al \in \pi}{\Pi} M_\al \; | \;
 \rh{\al}{\be}(m_{\al \be})=m_\al \otimes 1_\be \text{ for all }
 \al,\be \in \pi \}.
 $$
For any $\al \in \pi$, let $M^{\co H}_\al$ be the image of the (canonical) projection of this set onto $M_\al$. It
is easy to verify that $M^{\co H}=\{ M^{\co H}_\al \}_{\al \in \pi}$ is a right \p-subcomodule of $M$, called the
{\it \p-subcomodule of coinvariants}.
\end{example}

\subsection{Rational $\pi$-graded modules}\label{ratgradmod}
Throughout this subsection, $C=(\{C_\al\},\Delta,\epsilon)$ will denote a \p-coalgebra and $C^*=\oplus_{\al \in
\pi} C^*_\al$ its dual \p-graded algebra (see \S\ref{convo}). In this subsection we explore the relationships
between right \p-comodules over $C$ and \p-graded left $C^*$-modules.

Let $M=\oplus_{\al \in \pi} M_\al$ be a \p-graded left $C^*$-module with action $\psi : C^* \otimes M \to M$. Set
$\oM_\al=M_\ali$. For any $\al, \be \in \pi$, define
\begin{equation}\label{ratmap1}
 \rh{\al}{\be} : \oM_{\al \be} \to
 \homo_\kk(C^*_\be,\oM_\al) \quad \text{by} \quad \rh{\al}{\be}(m)(f)=\psi(f \otimes m).
\end{equation}
There is a  natural embedding
  $$
  \oM_\al \otimes C_\be \hookrightarrow \homo_\kk(C^*_\be,\oM_\al) \quad
  m \otimes c \mapsto (f \mapsto f(c)m) .
  $$
Regard this embedding as inclusion, so that $\oM_\al \otimes C_\be \subset \homo_\kk(C^*_\be,\oM_\al)$. The
\p-graded left $C^*$-module $M$ is said to be \emph{rational} provided $\rh{\al}{\be}(\oM_{\al\be}) \subset
\oM_\al \otimes C_\be$ for all $\al,\be \in \pi$. In this case, the restriction of $\rh{\al}{\be}$ onto $\oM_\al
\otimes C_\be$ will also be denoted by
\begin{equation}\label{ratmap2}
\rh{\al}{\be}: \oM_{\al \be} \to \oM_\al \otimes C_\be.
\end{equation}

The definition given here generalizes that of a rational $\pi$-graded left module given in \cite{NT1} and agrees
with it when $\pi$ is finite.

The next theorem generalizes \cite[Theorem 6.3]{NT1} and \cite[Theorem 2.1.3]{sweed}.

\begin{theorem}\label{theorational}
Let $C$ be a \p-coalgebra. Then
\begin{enumerate}
\renewcommand{\labelenumi}{{\rm (\alph{enumi})}}
 \item There is a one-to-one correspondence between (isomorphic classes of)
        right \p-comodules over $C$ and (isomorphic classes of) rational \p-graded left
        $C^*$-modules.
  \item Every graded submodule of a rational \p-graded left
        $C^*$-module is rational.
  \item Any \p-graded left $C^*$-module $L=\oplus_{\al \in \pi} L_\al$ has a
        unique maximal rational graded submodule, noted $L^\rat$, which is equal to the
        sum of all rational graded submodules of $L$. Moreover, if
        $\rho=\{\rh{\al}{\be} \}_{\al,\be \in \pi}$ is defined as in \eqref{ratmap1},
        then $(L^\rat)_\ga = \underset{\substack{\al,\be
\in \pi \\ \al \be =\gai}}{\cap} \rh{\al}{\be}^{-1}(\oL_\al \otimes C_\be)$ for any $\ga \in \pi$.
\end{enumerate}
\end{theorem}
Before proving the theorem, we needs two lemmas. Recall that a left module $M$ over a \p-graded algebra
$A=\oplus_{\al \in \pi} A_\al$ is \emph{graded} if $M$ admits a decomposition as a direct sum of $\Bbbk$-spaces
$M=\oplus_{\al \in \pi} M_\al$ such that $A_\al M_\be \subset M_{\al \be}$ for all $\al,\be \in \pi$. A submodule
$N$ of $M$ is graded if $N=\oplus_{\al \in \pi} (N \cap M_\al)$. The quotient $M/N$ is then a left \p-graded
$A$-module by setting $(M/N)_\al=(M_\al+N)/N$ for all $\al \in \pi$. This is the unique structure of a \p-graded
$A$-module on $M/N$ which makes the canonical projection $M \to M/N$ a graded $A$-morphism.

Let $M=\{M_\al\}_{\al \in \pi}$ be a family of $\Bbbk$-spaces and $\rho=\{\rh{\al}{\be} :M_{\al \be} \to M_\al
\otimes C_\be \}_{\al, \be \in \pi}$ be a family of $\Bbbk$-linear maps. Set $\oM=\oplus_{\al \in \pi} \oM_\al$,
where $\oM_\al=M_\ali$. Let $\psi_\rho: C^* \otimes \oM \to \oM$ be the $\Bbbk$-linear map defined on the summands
by
\begin{multline*}
  C^*_\al \otimes \oM_\be \rTo^{\id_{C^*_\al} \otimes \, \rh{( \al \be)^{-1}}{\al}} C^*_\al \otimes \oM_{\al \be}
  \otimes C_\al \rLine^{\sigma_{C^*_\al,\oM_{\al \be}} \otimes \, \id_{C_\al}} \\
  \rTo \oM_{\al \be} \otimes C^*_\al \otimes
  C_\al \rTo^{\id_{\oM_{\al \be}} \otimes \, \langle, \rangle} \oM_{\al \be} \otimes \kk \cong \oM_{\al \be},
\end{multline*}
where $\langle\, ,\rangle$ denotes the natural pairing between $C^*_\al$ and $C_\al$.
\begin{lemma}\label{rational1}
$(M,\rho)$ is a right \p-comodule over $C$ if and only if $(\oM,\psi_\rho)$ is a \p-graded left $C^*$-module.
\end{lemma}
\begin{proof}
Suppose that $(M,\rho)$ is a right \p-comodule over $C$. Firstly, for any $m \in \oM_\al$, $\psi_\rho(\epsilon
\otimes m)=\mzai \epsilon(\muu)=m$, by \eqref{counitcomod}. Secondly, for any $f \in C^*_\al$, $g \in C^*_\be$,
and $m \in \oM_\ga$,
\begin{eqnarray*}
\psi_\rho(fg \otimes m)
  & = & m_{(0,(\al \be \ga)^{-1})} fg(m_{(1,\al \be)}) \\
  & = & m_{(0,(\al \be \ga)^{-1})} f(\mua) g(\mdb) \\
  & = & \psi_\rho(f \otimes m_{(0,(\be \ga)^{-1})}g(\mub)) \\
  & = & \psi_\rho(f \otimes  \psi_\rho(g \otimes m)).
\end{eqnarray*}
Moreover, by construction, $\psi_\rho(C^*_\al \otimes \oM_\be) \subset \oM_{\al \be}$ for any $\al,\be \in \pi$.
Hence $(\oM,\psi_\rho)$ is a \p-graded left $C^*$-module.

Conversely, suppose that $(\oM,\psi_\rho)$ is a left \p-graded $C^*$-module. Since, for all $\al \in \pi$ and  $m
\in M_\al=\oM_\ali$, $(\id_{M_\al} \otimes \epsilon) \rh{\al}{1}(m)=\psi_\rho(\epsilon \otimes m)=m$,
\eqref{counitcomod} is verified. To show that \eqref{coasscomod} is satisfied, let $\al,\be,\ga \in \pi$ and $m
\in M_{\al \be \ga}$. Set
 $$
 \delta=(\rh{\al}{\be} \otimes \id_{C_\ga}) \rh{\al
 \be}{\ga} (m) - (\id_{M_\al} \otimes \cp{\be}{\ga}) \rh{\al}{\be \ga}(m) \in M_\al \otimes C_\be \otimes C_\ga.
 $$
Suppose that $\delta \neq 0$. Then there exists $F \in (M_\al \otimes C_\be \otimes C_\ga)^*$ such that $F(\delta)
\neq 0$. Now $M_\al^* \otimes C_\be^* \otimes C_\ga^*$ is dense in the linear topological space $(M_\al \otimes
C_\be \otimes C_\ga)^*$ endowed with the $(M_\al \otimes C_\be \otimes C_\ga)$-topology (see \cite[page 70]{abe}).
Thus $(M_\al^* \otimes C_\be^* \otimes C_\ga^*) \cap (F+\delta^\bot) \neq \emptyset$, where $\delta^\bot=\{ f \in
(M_\al \otimes C_\be \otimes C_\ga)^* \, | \; f(\delta)=0\}$.  Then there exists $G \in M_\al^* \otimes C_\be^*
\otimes C_\ga^*$ such that $G(\delta) \neq 0$. Now, for all $f \in M^*_\al$, $g \in C^*_\be$, and $h \in C^*_\ga$,
\begin{eqnarray*}
(f \otimes g \otimes h)(\rh{\al}{\be} \otimes \id_{C_\ga}) \rh{\al \be}{\ga} (m)
  & = & f\circ \psi_\rho(g \otimes  \psi_\rho(h \otimes m)) \\
  & = & f \circ \psi_\rho(gh \otimes m) \\
  & = & (f \otimes g \otimes h) (\id_{M_\al} \otimes
        \cp{\be}{\ga}) \rh{\al}{\be \ga}(m),
\end{eqnarray*}
i.e., $(f \otimes g \otimes h)(\delta)=0$. Therefore $G(\delta)=0$, which is a contradiction. We conclude that
$\delta=0$ and then $(\rh{\al}{\be} \otimes \id_{C_\ga}) \rh{\al \be}{\ga}=(\id_{M_\al} \otimes \cp{\be}{\ga})
\rh{\al}{\be \ga}$. Hence $(M,\rho)$ is a right \p-comodule over $C$.
\end{proof}

\begin{lemma}\label{rational2}
Let $M=\oplus_{\al \in \pi} M_\al$ be a rational \p-graded left $C^*$-module. Then $\oM=\{ \oM_\al \}_{\al \in
\pi}$, endowed with the structure maps $\rho = \{ \rh{\al}{\be} \}_{\al \be \in \pi}$ defined by \eqref{ratmap2},
is a right \p-comodule over $C$.
\end{lemma}
\begin{proof}
Let  $\psi_\rho: C^* \otimes \overline{\oM} \to \overline{\oM}$ be the map defined as in Lemma~\ref{rational1}. It
is easy to verify that $(\overline{\oM},\psi_\rho)=(M,\psi)$. Thus $(\overline{\oM},\psi_\rho)$ is a \p-graded
left $C^*$-module and hence, by Lemma~\ref{rational1}, $(\oM,\rho)$ is a right \p-comodule over $C$.
\end{proof}

\begin{proof}[Proof of Theorem~\ref{theorational}]
Part (a) follows directly from Lemmas~\ref{rational1} and \ref{rational2}. To show Part (b), let $N$ be a graded
submodule of a rational \p-graded left $C^*$-module $(M,\psi)$. Let $\rho_{\al,\be}:\oN_{\al\be} \to
\homo_\kk(C_\be^*,\oN_\al)$ defined by $\rho_{\al,\be}(m)(f)=\psi(f \otimes m)$. Suppose that there exist $\al,\be
\in \pi$ and $n \in \oN_{\al \be}$ such that $\rh{\al}{\be}(n) \not \in \oN_\al \otimes C_\be$. Since $M$ is
rational, we can write $\rh{\al}{\be}(n)=\sum_{i=1}^{k}n_i \otimes c_i \in \oM_\al \otimes C_\be$. Without loss of
generality, we can assume that the $c_i$ are $\Bbbk$-linearly independent and $n_1 \not \in \oN_\al$. Let $f \in
C_\be^*$ such that $f(c_1)=1$ and $f(c_i)=0$ for $i \geq 2$. Now $\psi(f \otimes n) =\sum_{i=1}^k n_i f(c_i) =n_1
\not \in \oN_\al=N_\ali$, contradicting the fact that $N$ is a graded submodule of $M$. Thus
$\rh{\al}{\be}(\oN_{\al \be}) \subset \oN_\al \otimes C_\be$ for all $\al,\be \in \pi$. Hence $N$ is rational.

Let us show Part (c). Denote by $\cdot$ the left action of $C^*$ on $L$. Set $\oL_\al = L_\ali$ and $
\rh{\al}{\be} : \oL_{\al \be} \to \homo_\kk(C^*_\be,\oL_\al)$ given by $\rh{\al}{\be}(m)(f)=f \cdot m$. Recall
$\oL_\al \otimes C_\be$ can be viewed as a subspace of $\homo_\kk(C^*_\be,\oL_\al)$ via the embedding $\oL_\al
\otimes C_\be \hookrightarrow \homo_\kk(C^*_\be,\oL_\al)$ given by $m \otimes c \mapsto (f \mapsto f(c)m)$. Define
$M_\ga=\cap_{\al \be =\gai} \rh{\al}{\be}^{-1}(\oL_\al \otimes C_\be)$ for any $\ga \in \pi$, and set
$M=\oplus_{\ga \in \pi} M_\ga$. Fix $\al,\be \in \pi$, $f \in C^*_\al$, and $m \in M_\be$. Let $u,v \in \pi$ such
that $uv=(\al \be)^{-1}$. We can write $\rh{u}{v \al}(m)=\sum_{i=1}^k l_i \otimes c_i \in \oL_u \otimes C_{v\al}$.
Now, for any $g \in C^*_v$,
 $
 g \cdot (f \cdot m)  =  (gf) \cdot m
    =  \sum_{i=1}^k  gf(c_i)\,l_i
    =  \sum_{i=1}^k g(f(c_{i(2,\al)}) c_{i(1,v)})\, l_i
 $.
Then $\rh{u}{v}(f \cdot m)= \sum_{i=1}^k l_i \otimes f(c_{i(2,\al)}) c_{i(1,v)} \in \oL_u \otimes C_v$ and so $f
\cdot m \in \rh{u}{v}^{-1}(\oL_u \otimes C_v)$. Hence $f \cdot m \in \cap_{uv =(\al \be)^{-1}}
\rh{u}{v}^{-1}(\oL_u \otimes C_v)=M_{\al \be}$. Therefore $M$ is a graded submodule of $L$. Moreover one easily
checks at this point that $\rh{\al}{\be}(\oM_{\al\be}) \subset \oM_\al \otimes C_\be$ for any $\al,\be$ in $\pi$.
Thus $M$ is rational.

Suppose now that $N$ is another rational graded submodule of $L$ and denote by $\varrho=\{ \vrh{\al}{\be}
\}_{\al,\be \in \pi}$ its corresponding \p-comodule structure maps (see Lemma~\ref{rational2}). Let $\ga \in \pi$
and $\al, \be \in \pi$ such that $\al \be=\gai$. By the definition of $\rh{\al}{\be}$ and $\vrh{\al}{\be}$ and of
the embedding $\oN_\al \otimes C_\be \subset \oL_\al \otimes C_\be \subset \homo_\kk(C^*_\be,\oL_\al)$, it follows
that $\rho_{\al,\be|N}=\vrh{\al}{\be} : \oN_{\al \be} \to \oN_\al \otimes C_\be$. Thus
$\rh{\al}{\be}(N_\ga)=\vrh{\al}{\be}(\oN_{\al \be}) \subset \oN_\al \otimes C_\be \subset \oL_\al \otimes C_\be$,
and so $N_\ga \subset \rh{\al}{\be}^{-1}(\oL_\al \otimes C_\be)$. This holds for all $\al, \be \in \pi$ such that
$\al \be=\gai$. Thus $N_\ga \subset \cap_{\al \be =\gai} \rh{\al}{\be}^{-1}(\oL_\al \otimes C_\be)=M_\ga$ for any
$\ga \in \pi$. Hence $N \subset M$. Therefore $M$ is the unique maximal rational graded submodule of $L$ and is
the sum of all rational graded submodules of $L$.
\end{proof}
Remark that, using Lemma~\ref{rational2} and Theorem~\ref{theorational}(c), a unique ``maximal" right \p-comodule
$\overline{M^\rat}$ over $C$ can be associated to any \p-graded left $C^*$-module $M$.

\subsection{Hopf $\pi$-comodules}\label{defiHopfmod}
In this subsection, we introduce and discuss the notion of a Hopf $\pi$-comodule.

Let $H=(\{H_\al,m_\al,1_\al\},\Delta,\epsilon,S)$ be a Hopf \p-coalgebra. A \emph{right Hopf \p-comodule over} $H$
is a right \p-comodule $M=\{M_\al \}_{\al \in \pi}$ over $H$ such that
\begin{defi}
 \item \label{hopfmod1} $M_\al$ is a right $H_\al$-module for any $\al \in \pi$;
 \item \label{hopfmod2} Let us denote by $\psi_\al:M_\al \otimes H_\al \to H_\al$
                        the right action of $H_\al$
                        on $M_\al$ and by $\rho=\{ \rh{\al}{\be} \}_{\al,
                        \be \in \pi}$ the \p-comodule maps of $M$. These structures are
                        required to be compatible in the sense that,
                        for any $\al,\be \in \pi$, the following diagram is commutative:
      \begin{diagram}[labelstyle=\scriptstyle]
        M_{\al \be} \otimes H_{\al \be} & \rTo^{\psi_{\al \be}} & M_{\al \be}
        & \rTo^{\rh{\al}{\be}} & M_\al \otimes H_\be \\
        \dTo^{\rh{\al}{\be} \otimes \cp{\al}{\be}} &&&& \uTo_{ \psi_\al
        \otimes m_\be}\\
        M_\al \otimes H_\be \otimes H_\al \otimes H_\be &&
        \rTo_{\id_{M_\al} \otimes \sigma_{H_\be,H_\al} \otimes \id_{H_\be}}
        && M_\al \otimes H_\al \otimes H_\be \otimes H_\be
      \end{diagram}
\end{defi}

When $\pi=1$, one recovers the definition of a Hopf module (see \cite{LS}).

Note that Axiom \eqref{hopfmod2} means that $\rh{\al}{\be}:M_{\al\be} \to M_\al \otimes H_\be$ is $H_{\al
\be}$-linear, where $M_\al \otimes H_\be$ is endowed with the right $H_{\al \be}$-module structure given by $$(m
\otimes h) \cdot a= \psi_\al(m \otimes \aua) \otimes h \,\adb.$$

A \emph{Hopf \p-subcomodule of $M$} is a \p-subcomodule $N=\{N_\al\}_{\al \in \pi}$ of $M$ such that $N_\al$ is a
$H_\al$-submodule of $M_\al$ for any $\al \in \pi$. Then $N$ is a right Hopf \p-comodule over $H$.

A \emph{Hopf \p-comodule morphism} between two right Hopf \p-comodules $M$ and $M'$ is a \p-comodule morphism
$f=\{f_\al : M_\al \to M'_\al\}_{\al \in \pi}$ between $M$ and $M'$ such that $f_\al$ is $H_\al$-linear for any
$\al \in \pi$.

\begin{example} \label{trivial}
Let $H$ be a Hopf \p-coalgebra and $M=\{M_\al \}_{\al \in \pi}$ be a right \p-comodule over $H$, with structure
maps $\rho=\{ \rh{\al}{\be} \}_{\al,\be \in \pi}$. For any $\al \in \pi$, set $(M \otimes H)_\al=M_\al \otimes
H_\al$. The multiplication in $H_\al$ induces a structure of a right $H_\al$-module on $(M \otimes H)_\al$ by
setting $(m \otimes h)\vartriangleleft a = m \otimes ha$. Define the \p-comodule structure maps $\xi_{\al,\be} :(M
\otimes H)_{\al \be} \to (M \otimes H)_\al \otimes H_\be$ by
  $$
  \xi_{\al,\be} (m \otimes h)=\mza \otimes h_{(1,\al)}
  \otimes \mub h_{(2,\be)}.
  $$
Here we write as usual $\rh{\al}{\be}(m)= \mza \otimes \mub$ and $\cp{\al}{\be}(h)=h_{(1,\al)} \otimes
h_{(2,\be)}$. One easily verifies that $M \otimes H=\{ ( M \otimes H)_\al \}_{\al \in \pi}$ is a right Hopf
\p-comodule over $H$, called \emph{trivial}.
\end{example}
In the next theorem, we show that a Hopf \p-comodule can be canonically  decomposed. This generalizes the
fundamental theorem of Hopf modules (see \cite[Proposition 1]{LS}).
\begin{theorem}\label{decompHopfmod}
Let $H$ be a Hopf \p-coalgebra and $M=\{M_\al \}_{\al \in \pi}$ be a right Hopf \p-comodule over $H$. Consider the
\p-subcomodule of coinvariants $M^{\co H}$ of $M$ (see Example~\ref{coinvariants}) and the trivial right Hopf
\p-comodule $M^{\co H} \otimes H$ (see Example~\ref{trivial}). Then there exists a Hopf \p-comodule isomorphism $M
\cong M^{\co H} \otimes H$.
\end{theorem}
\begin{proof}
We will denote by $\cdot$ (resp.\@ $\vartriangleleft$) the right action of $H_\al$ on $M_\al$ (resp.\@ on $(M^{\co
H} \otimes H)_\al$) and by $\rho=\{ \rh{\al}{\be} \}_{\al,\be \in \pi}$ (resp.\@ $\xi=\{ \xi_{\al,\be} \}_{\al,\be
\in \pi}$) the \p-comodule structure maps of $M$ (resp.\@ of $M^{\co H} \otimes H$). For any $\al \in \pi$, define
$P_\al:M_1 \to M_\al$ by $P_\al(m)=\mza \cdot S_\ali(\muai)$. Remark first that, for any $m \in M_1$,
$(P_\al(m))_{\al \in \pi}$ is a coinvariant of $H$ on $M$. Indeed, if $m \in M_1$,
\begin{eqnarray*}
\lefteqn{ \rh{\al}{\be}( P_{\al\be}(m))}\\
    & = & \rh{\al}{\be}(\mzab \cdot S_{(\al \be)^{-1}} (\muabi)) \\
    & = & \rh{\al}{\be}(\mzab) \cdot \cp{\al}{\be} S_{(\al \be)^{-1}} (\muabi)
          \text{ \quad by \eqref{hopfmod2}}\\
    & = & \mza \cdot S_\ali(\mtai) \otimes \mub S_\bei(\mdbi) \quad \text{ by Lemma~\ref{antipodepptes}(c)}\\
    & = & \mza \cdot S_\ali(\epsilon(\muu)\mdai) \otimes  1_\be \text{ \quad by \eqref{antipode}}\\
    & = & \mza \cdot S_\ali(\muai) \otimes  1_\be \text{ \quad by \eqref{counit}}\\
    & = & P_\al(m) \otimes  1_\be.
\end{eqnarray*}
For any $\al \in \pi$, define $f_\al : (M^{\co H} \otimes H)_\al \to  M_\al$ by $f(m \otimes h)= m \cdot h$. Then
$f_\al$ is $H_\al$-linear since $f_\al(m \otimes h) \cdot a = (m \cdot h) \cdot a = m \cdot ha = f_\al((m \otimes
h) \vartriangleleft a)$ for all $m \in M_\al^{\co H}$ and $h,a \in H_\al$. Moreover $(f_\al \otimes
\id_{H_\be})\xi_{\al,\be}= \rh{\al}{\be}f_{\al \be}$ for all $\al,\be \in \pi$. Indeed let $m \in M_{\al \be}^{\co
H}$ and $h \in H_{\al \be}$. By the definition of $M_{\al\be}^{\co H}$, there exists a coinvariant $(m_\ga)_{\ga
\in \pi}$ of $H$ on $M$ such that $m=m_{\al \be}$. In particular $\rh{\al}{\be}(m)=m_\al \otimes 1_\be$. Thus
\begin{eqnarray*}
 (f_\al \otimes \id_{H_\be})\xi_{\al,\be}(m \otimes h)
   & = & m_\al \cdot \hua \otimes \hdb\\
   & = & \rh{\al}{\be}(m) \cdot \cp{\al}{\be}(h) \\
   & = & \rh{\al}{\be}(m \cdot h) \text{ \quad by \eqref{hopfmod2}}\\
   & = & \rh{\al}{\be}(f_{\al \be}(m \otimes h)).
\end{eqnarray*}
Then $f=\{ f_\al \}_{\al \in \pi}: M^{\co H} \otimes H \to M$ is a Hopf \p-comodule morphism. To show that $f$ is
an isomorphism, we construct its inverse. For any $\al \in \pi$, define $g_\al : M_\al \to M_\al^{\co H} \otimes
H_\al$ by $g_\al= (P_\al \otimes \id_{H_\al}) \rh{1}{\al}$. The map $g_\al$ is well-defined since $(P_\ga(m))_{\ga
\in \pi}$ is a coinvariant of $H$ on $M$ for all $m \in M_1$, and is $H_\al$-linear since, for any $x \in M_\al$
and $a \in H_\al$,
\begin{eqnarray*}
 g_\al(x \cdot a)
  & = & (P_\al \otimes \id_{H_\al}) \rh{1}{\al}( x \cdot a) \\
  & = & (P_\al \otimes \id_{H_\al}) (\rh{1}{\al}(x) \cdot
        \cp{1}{\al}(a))
        \text{ \quad by \eqref{hopfmod2}}\\
  & = & P_\al(\xzu \cdot \auu) \otimes \xua \ada \\
  & = & (\xza \cdot \aua) \cdot S_\ali(\xuai \adai) \otimes
        \xda \ata \text{ \quad by \eqref{hopfmod2}}\\
  & = & \xza \cdot (\aua S_\ali(\adai) S_\ali(\xuai) ) \otimes
        \xda \ata \\
  & = & \xza \cdot S_\ali(\xuai) \otimes \xda
        \epsilon(\auu) \ada \text{ \quad by \eqref{antipode}}\\
  & = & \xza \cdot S_\ali(\xuai) \otimes \xda a \text{ \quad by \eqref{counit}}\\
  & = & g_\al(x) \vartriangleleft a.
\end{eqnarray*}
Moreover $(g_\al \otimes \id_{H_\be}) \rh{\al}{\be}=\xi_{\al,\be} g_{\al \be} $ for all $\al,\be \in \pi$. Indeed,
for any $x \in M_{\al \be}$,
\begin{eqnarray*}
 \xi_{\al,\be}(g_{\al \be}(x))
  & = & \xi_{\al,\be}(P_{\al \be}(\xzu) \otimes x_{(1,\al \be)}) \\
  & = & P_{\al \be}(\xzu)_{(0,\al)} \otimes x_{(1,\al \be)(1,\al)}
        \otimes P_{\al \be}(\xuu)_{(1,\be)} x_{(1,\al \be)(2,\be)},
\end{eqnarray*}
and so, since $(P_\ga(\xzu))_{\ga \in \pi}$ is a \p-coinvariant of $H$ on $M$,
\begin{eqnarray*}
 \xi_{\al,\be}(g_{\al \be}(x))
   & = & P_\al(\xzu) \otimes \xua \otimes \xdb \\
   & = & g_\al(\xza) \otimes \xub \\
   & = & (g_\al \otimes \id_{H_\be}) \rh{\al}{\be}(x).
\end{eqnarray*}
Thus $g=\{ g_\al \}_{\al \in \pi} : M \to M^{\co H} \otimes H$ is a Hopf \p-comodule morphism. It remains now to
verify that $g_\al f_\al=\id_{M^{\co H}_\al \otimes H_\al}$ and $f_\al g_\al=\id_{M_\al}$ for any $\al \in \pi$.
Let $m \in M^{\co H}_\al$ and $h \in H_\al$. By the definition of $M_{\al}^{\co H}$, there exists a coinvariant
$(m_\ga)_{\ga \in \pi}$ of $H$ on $M$ such that $m=m_\al$. In particular, $\rh{1}{\al}(m)=m_1 \otimes 1_\al$ and
$P_\al(m_1)=m_\al \cdot S_\ali(1_\ali)=m \cdot 1_\al = m$. Then
\begin{eqnarray*}
 g_\al  f_\al(m \otimes h)
    & = &  g_\al(m \cdot h) \\
    & = & g_\al(m) \vartriangleleft h \text{ \quad since $g_\al$ is
          $H_\al$-linear}\\
    & = & (P_\al(m_1) \otimes 1_\al) \vartriangleleft h \\
    & = & m \otimes h.
\end{eqnarray*}
Finally, for all $x \in M_\al$,
\begin{eqnarray*}
 f_\al g_\al(x)
    & = & (\xza \cdot S_\ali(\xuai) ) \cdot \xda\\
    & = & \xza \cdot (S_\ali(\xuai) \,\xda)\\
    & = & \xza \epsilon(\xuu) \cdot 1_\al \text{ \quad by \eqref{antipode}}\\
    & = & x \text{ \quad by \eqref{counitcomod}}.
\end{eqnarray*}
Hence $g=f^{-1}$ and $f$ and $g$ are Hopf \p-comodule isomorphisms.
\end{proof}

\section{Existence and uniqueness of $\pi$-integrals}\label{piint}
In this section, we introduce and discuss the notion of a \p-integral for a Hopf \p-coalgebra. In particular, by
generalizing the arguments of \cite[\S 5]{sweed}, we show that, in the finite dimensional case, the space of left
(resp.\@ right) \p-integrals is one dimensional.

\subsection{$\pi$-integrals}\label{s:defiint}
We first recall that a left (resp.\@ right) integral for a Hopf algebra $(A,\Delta,\epsilon,S)$ is an element
$\Lambda \in A$ such that $x \Lambda = \epsilon(x) \Lambda$ (resp.\@ $\Lambda x = \epsilon(x) \Lambda$) for all $x
\in A$. A left (resp.\@ right) integral for the dual Hopf algebra $A^*$ is a $\kk$-linear form $\lambda\in A^*$
verifying $(f \otimes \lambda) \Delta=f(1_A) \lambda$ (resp.\@ $(\lambda \otimes f) \Delta=f(1_A) \lambda$) for
all $f \in A^*$. Let us extend this notion to the setting of a Hopf \p-coalgebra.

Let $H=(\{H_\al,m_\al,1_\al\},\Delta,\epsilon,S)$ be a Hopf \p-coalgebra. A \emph{left (resp.\@ right) \p-integral
for} $H$ is a family of $\kk$-linear forms $\lambda=(\la)_{\al \in \pi} \in \Pi_{\al \in \pi} H^*_\al$ such that,
for all $\al,\be \in \pi$,
\begin{equation}\label{defiint}
  (\id_{H_\al} \otimes \lb) \cp{\al}{\be} = \lab \, 1_\al
  \quad \text{(resp.} \quad (\la \otimes \id_{H_\be}) \cp{\al}{\be} = \lab \, 1_\be \,\text{).}
\end{equation}

Note that $\lu$ is a usual left (resp.\@ right) integral for the Hopf algebra $H_1^*$.

If we use the multiplication of the dual \p-graded algebra $H^*$ of $H$ (see \S\ref{convo}), we have that
$\lambda=(\la)_{\al \in \pi} \in \Pi_{\al \in \pi} H^*_\al$ is a left (resp.\@ right) \p-integral for $H$ if and
only if, for all $\al,\be \in \pi$ and $f \in H_\al^*$ (resp.\@ $g \in H_\be^*$),
  $$
  f \lb=f(1_\al)\, \lab \quad
  \text{(resp.} \quad \la  g=g(1_\al)\, \lab \,\text{).}
  $$

A \p-integral $\lambda=(\la )_{\al \in \pi}$ for $H$ is said to be \emph{non-zero} if $\lb \neq 0$ for some $\be
\in \pi$.

\begin{lemma}\label{intnonzero}
Let $\lambda=(\la)_{\al \in \pi}$ be a non-zero left (resp.\@ right) \p-integral for $H$. Then $\la \neq 0$ for
all $\al \in \pi$ such that $H_\al \neq 0$. In particular $\lu \neq 0$.
\end{lemma}
\begin{proof}
Let $\lambda=(\la )_{\al \in \pi}$ be a left \p-integral for $H$ such that $\lb \neq 0$ for some $\be \in \pi$ and
let $\al \in \pi$ such that $H_\al \neq 0$. Then $H_{\be \ali} \neq 0$ (by Corollary~\ref{subgroup}) and so
$1_{\be \ali}\neq 0$. Using \eqref{defiint}, we have that $(\id_{H_{\be\ali}} \otimes \la)\cp{\be\ali}{\al}=\lb
1_{\be\ali} \neq 0$. Hence $\la \neq 0$. The right case can be done similarly.
\end{proof}

\begin{remark}
Let $H$ be a finite dimensional Hopf \p-coalgebra. Consider the Hopf algebra $H^*$ dual to $H$ (see
\S\ref{dualhopf}). If $H_\al=0$ for all but a finite number of $\al \in \pi$, then $\lambda=(\la)_{\al \in \pi}\in
\Pi_{\al \in \pi} H^*_\al$ is a left (resp.\@ right) \p-integral for $H$ if and only if $\sum_{\al \in \pi} \la$
is a left (resp.\@ right) integral for $H^*$. If $H_\al \neq 0$ for infinitely many $\al \in \pi$, then $H^*$ is
infinite dimensional and thus does not have any non-zero left or right integral (see \cite{sweed2}). Nevertheless
we show in the next subsection that $H$ always has a non-zero \p-integral.
\end{remark}

\subsection{The space of $\pi$-integrals is one dimensional}\label{s:exists}
It is known (see \cite[Corollary 5.1.6]{sweed}) that the space of left (resp.\@ right) integrals for a finite
dimensional Hopf algebra is one dimensional. In this subsection, we generalize this result to finite dimensional
Hopf \p-coalgebras.

Let $H=\{H_\al\}_{\al \in \pi}$ be a Hopf \p-coalgebra (not necessarily finite dimensional). The dual \p-graded
algebra $H^*$ of $H$ (see \S\ref{convo}) is a \p-graded left $H^*$-module via left multiplication. Let
$(H^*)^\rat$ be its maximal rational \p-graded submodule (see Theorem~\ref{theorational}(c)). Denote by
$\HH=\overline{(H^*)^\rat}=\{\HH_\al \}_{\al \in \pi}$ the right \p-comodule over $H$ which corresponds to it by
Lemma~\ref{rational2}. Recall that $\HH_\al \subset H^*_\ali$ for any $\al \in \pi$. The \p-comodule structure
maps of $\HH$ will be denoted by $\rho=\{ \rh{\al}{\be} \}_{\al, \be \in \pi}$.
\begin{lemma}\label{lieqcoinv}
Let $\lambda=(\la)_{\al \in \pi} \in \Pi_{\al \in \pi} H^*_\al$. Then $\lambda$ is a left \p-integral for $H$ if
and only if $( \lai )_{\al \in \pi}$ is a coinvariant of $H$ on $\HH$ (see Example~\ref{coinvariants}).
\end{lemma}
\begin{proof}
Suppose that  $\lambda$ is a left \p-integral for $H$. Fix $\ga \in \pi$. Let $\al, \be \in \pi$ such that $\al
\be =\ga$.  We have that $\rh{\al}{\be}(\lambda_\gai)= \lambda_\ali \otimes 1_\be \in \overline{H^*_\al} \otimes
H_\be$ since $ f \lambda_\gai= f(1_\be)\,\lai $ for all $f \in H_\be^*$. Therefore $\lambda_\gai \in \cap_{\al \be
=\ga} \, \rh{\al}{\be}^{-1}(\overline{H^*_\al} \otimes H_\be) =H^{* \rat}_\gai=\HH_\ga$, see
Theorem~\ref{theorational}(c). Hence, since $\rh{\al}{\be}(\lambda_{(\al \be)^{-1}})= \lambda_\ali \otimes 1_\be$
for all $\al, \be \in \pi$, $( \lai )_{\al \in \pi}$  is a coinvariant of $H$ on $\HH$. Conversely, suppose that
$( \lai )_{\al \in \pi}$ is a coinvariant of $H$ on $\HH$. Let $\al , \be \in \pi$. Then $\rh{(\al
\be)^{-1}}{\al}(\lb) =\lambda_{\al \be} \otimes 1_\al$, i.e., $f \lb= f(1_\al)\,\lambda_{\al \be}$ for all $f \in
H^*_\al$. Hence $\lambda$ is a left \p-integral for $H$.
\end{proof}

For all $\al \in \pi$, we define a right $H_\al$-module structure on $\HH_\al$  by setting
 $$
 (f\leftharpoondown a)(x) =
 f( x S_\al(a))
 $$
for any $f \in \HH_\al$, $a \in H_\al$, and $x \in H_\ali$.
\begin{lemma}\label{hopfmodint}
$\HH$ is a right Hopf \p-comodule over $H$.
\end{lemma}
\begin{proof}
Let us first show that for any $\al,\be \in \pi$, $f \in \HH_{\al \be}$, $a \in H_{\al \be}$, and $g \in H^*_\be$,
\begin{equation}\label{eqdemohopfmod}
g(f \leftharpoondown a)= f_{(0,\al)} \leftharpoondown \aua \langle g, f_{(1,\be)} \adb\rangle ,
\end{equation}
where $\langle\, ,\rangle$ denotes the natural pairing between $H^*_\be$ and $H_\be$. Remark first that
\begin{eqnarray*}
 1_\be \otimes S_{\al \be}(a)
   & = & \epsilon(\adu)\, 1_\be \otimes S_{\al \be}( a_{(1,\al \be)}) \text{ \quad by \eqref{counit}}\\
   & = & S_\bei(\adbi) \atb \otimes S_{\al \be}( a_{(1,\al \be)}) \text{ \quad by \eqref{antipode}}\\
   & = & S_\al( a_{(1,\al)})_{(1,\be)} \adb \otimes
         S_\al( a_{(1,\al)})_{(2,(\al \be)^{-1})} \text{ \quad by Lemma~\ref{antipodepptes}(c).}
\end{eqnarray*}
Then, for all $x \in H_\ali$,
\begin{eqnarray*}
 \lefteqn{\xub \otimes  x_{(2,(\al \be)^{-1})} S_{\al \be}(a)} \\
   & \quad = & \xub S_\al( a_{(1,\al)})_{(1,\be)} \adb \otimes
          x_{(2,(\al \be)^{-1})} S_\al( a_{(1,\al)})_{(2,(\al
          \be)^{-1})}\\
   & \quad = & ( x S_\al( a_{(1,\al)}))_{(1,\be)} \adb
         \otimes ( x S_\al( a_{(1,\al)} ))_{(2,(\al
          \be)^{-1})}  \text{ \quad by \eqref{hopfmor},}
\end{eqnarray*}
and so
\begin{eqnarray*}
 g(f \leftharpoondown a) (x)
   & = & \langle g, \xub\rangle  \langle  f \leftharpoondown a,x_{(2,(\al
         \be)^{-1})}\rangle  \\
   & = & \langle g, \xub\rangle  \langle  f,x_{(2,(\al
         \be)^{-1})} S_{\al \be}(a)\rangle  \\
   & = & \langle g,( x S_\al( a_{(1,\al)}) )_{(1,\be)} \adb\rangle
         \langle f, ( x S_\al( a_{(1,\al)} ) )_{(2,(\al
          \be)^{-1})}\rangle \\
   & = &  ( (\adb \rightharpoonup g)f ) \leftharpoondown
         \aua  (x),
\end{eqnarray*}
where $\rightharpoonup$ is the left $H_\be$-action on $H_\be^*$ defined by $(b \rightharpoonup l)(y)=l(yb)$ for
any $l \in H_\be^*$ and $b,y \in H_\be$. Then
\begin{eqnarray*}
 g(f \leftharpoondown a)
   & = & ( (\adb \rightharpoonup g)f) \leftharpoondown
         \aua \\
   & = & (f_{(0,\al)}\langle  \adb \rightharpoonup g , f_{(1,\be)}\rangle  ) \leftharpoondown
         \aua \text{ \quad by definition of $\rh{\al}{\be}$}\\
   & = & f_{(0,\al)}\leftharpoondown \aua \langle g , f_{(1,\be)}\adb\rangle ,
\end{eqnarray*}
and hence \eqref{eqdemohopfmod} is proved.

Recall that the \p-comodule structure map $\rh{\al}{\be}$ of $\HH$ is, via the natural embedding $\HH_\al \otimes
H_\be\subset \oHsa \otimes H_\be \hookrightarrow \homo_\kk(H^*_\be, \oHsa)$, the restriction onto $\HH_\al \otimes
H_\be$ of the map $\xi_{\al,\be}: \HH_{\al \be} \to \homo_\kk(H^*_\be, \oHsa)$ defined by
$\xi_{\al,\be}(f)(g)=gf$. Let $\ga \in \pi$. By \eqref{eqdemohopfmod}, we have that, for any $\al,\be \in \pi$
such that $\al \be = \ga$, $f \in \HH_\ga$, and $a \in H_\ga$,
  $$
  \xi_{\al,\be}(f \leftharpoondown a) = f_{(0,\al)} \leftharpoondown \aua \otimes
  f_{(1,\be)} \adb \in (\HH_\al \leftharpoondown \aua) \otimes H_\be \subset \oHsa \otimes H_\be.
$$
Therefore, by Theorem~\ref{theorational}(c), $ f \leftharpoondown a \in \cap_{\al \be =\ga}
\xi_{\al,\be}^{-1}(\oHsa \otimes C_\be)=\HH_\ga$. Hence the action of $H_\ga$ on $\HH_\ga$ is well-defined. This
is a right action because $S_\ga$ is unitary and anti-multiplicative (see Lemma~\ref{antipodepptes}). Finally,
Axiom \eqref{hopfmod2} is satisfied since \eqref{eqdemohopfmod} says that $\rh{\al}{\be}(f \leftharpoondown
a)=f_{(0,\al)}\leftharpoondown \aua \otimes f_{(1,\be)} \adb$ for any $\al,\be \in \pi$, $f \in \HH_{\al \be}$,
and $a \in H_{\al \be}$. Thus $\HH$ is a right Hopf \p-comodule over $H$.
\end{proof}

By Theorem~\ref{decompHopfmod}, the Hopf \p-comodule $\HH$ is isomorphic to the Hopf \p-comodule $(H^\pc)^{\co H}
\otimes H$. Let $f=\{ f_\al : (H^\pc)^{\co H}_\al \otimes H_\al \to \HH_\al \}_{\al \in \pi}$ be the right Hopf
\p-comodule isomorphism between them as in the proof Theorem~\ref{decompHopfmod}. Recall that $f_\al (m \otimes
h)= m \leftharpoondown h$ for any $\al \in \pi$, $m \in (H^\pc)^{\co H}_\al$, and $h \in H_\al$.
\begin{lemma}\label{Sinj}
If there exists a non-zero left \p-integral for $H$, then $S_\al$ is injective for all $\al \in \pi$.
\end{lemma}
\begin{proof}
Suppose that  $\lambda=(\la )_{\al \in \pi}$ is a non-zero left \p-integral for $H$.  Let $\al \in \pi$. If
$H_\al=0$, then the result is obvious. Let us suppose that $H_\al \neq 0$. Then $H_\ali \neq 0$ by
Corollary~\ref{subgroup} and so $\lai \neq 0$ (by Lemma~\ref{intnonzero}). Let $h \in H_\al$ such that
$S_\al(h)=0$. By Lemma~\ref{lieqcoinv}, $\lai \in H^{\pc \co H}_\al$. Now $f_\al(\lai \otimes h)=\lai
\leftharpoondown h=0$ (since $S_\al(h)=0$). Thus $\lai \otimes h=0$ (since $f_\al$ is an isomorphism) and so $h=0$
(since $\lai \neq 0$).
\end{proof}

\begin{theorem}\label{piintegexists}
Let $H$ be a finite dimensional Hopf \p-coalgebra. Then the space of left (resp.\@ right) \p-integrals for $H$ is
one-dimensional.
\end{theorem}
\begin{proof}
For any $\al,\be \in \pi$, since $H$ is finite dimensional and $\oHsa=H^*_\ali$, we have that  $\dim\oHsa \otimes
H_\be=\dim\homo_\kk(H^*_\be,\oHsa) < +\infty$. Therefore the natural embedding $\oHsa \otimes H_\be
\hookrightarrow \homo_\kk(H^*_\be,\oHsa)$ is an isomorphism. Thus $H^*$ is a rational \p-graded $H^*$-module (see
\S\ref{ratgradmod}) and so $\HH_\al=H_\ali^*$ for all $\al \in \pi$. Now $\dim (H^\pc)^{\co H}_1=1$ since
$(H^\pc)^{\co H}_1 \otimes H_1 \cong \HH_1$, $\dim H_1=\dim H_1^*=\dim \HH_1 < +\infty$, and $\dim H_1 \neq 0$ (by
Corollary~\ref{subgroup}). Hence there exists a \p-coinvariant $(\psi_\al)_{\al \in \pi}$ of $H$ on $\HH$ such
that $\psi_1 \neq 0$. Set $\la = \psi_\ali$ for any $\al \in \pi$. By Lemma~\ref{lieqcoinv}, $\lambda=( \la )_{\al
\in \pi}$ is then a left \p-integral for $H$. Moreover $\lu = \psi_1 \neq 0$ and so $\lambda$ is non-zero.

Suppose now that $\delta=( \delta_\al )_{\al \in \pi}$ is another left \p-integral for $H$. Let $\al \in \pi$ such
that $H_\al \neq 0$. By Lemma~\ref{Sinj}, $S_\al$ and $S_\ali$ are injective (since there exists a non-zero left
integral for $H$) and so $\dim H_\al=\dim H_\ali$. Therefore $\dim (H^\pc)^{\co H}_\al =1$ since $(H^\pc)^{\co
H}_\al \otimes H_\al \cong \HH_\al$ and  $0 \neq \dim H_\al=\dim \HH_\al < +\infty$. Now $\lai, \delta_\ali \in
(H^\pc)^{\co H}_\al$ by Lemma~\ref{lieqcoinv} and $\lai \neq 0$ (by Lemma~\ref{intnonzero}). Hence there exists
$k_\al \in \kk$ such that $\delta_\ali = k_\al \, \lai$. If $\al \in \pi$ is such that $H_\al \neq 0$, then
 $$
 k_1 \lu \, 1_\al = \delta_1 \, 1_\al = (\id_{H_\al} \otimes
 \delta_\ali) \cp{\al}{\ali}=k_\al (\id_{H_\al} \otimes
 \lai) \cp{\al}{\ali}= k_\al \lu \, 1_\al,
 $$
and thus $k_\al=k_1$ (since $\lu \neq0$ and $1_\al \neq 0$). If $\al \in \pi$ is such that $H_\al=0$, then
$\delta_\al=0=\la$ and so $\delta_\al=k_1 \, \la$. Hence we can conclude that $\delta$ is a scalar multiple of
$\lambda$.

To show the existence and the uniqueness of right \p-integrals for $H$, it suffices to consider the coopposite
Hopf \p-coalgebra $H^\cop$ to $H$ (see \S\ref{coop}). Indeed $\lambda=(\la)_{\al \in \pi} \in \Pi_{\al \in \pi}
H^*_\al$ is a right \p-integral for $H$ if and only if $(\lai )_{\al \in \pi}$ is a left \p-integral for $H^\cop$.
This completes the proof of the theorem.
\end{proof}

\begin{corollary}\label{corinteg}
Let $H=\{H_\al\}_{\al \in \pi}$ be a finite dimensional Hopf \p-coalgebra. Then
\begin{enumerate}
\renewcommand{\labelenumi}{{\rm (\alph{enumi})}}
 \item The antipode $S=\{S_\al\}_{\al \in \pi}$ of $H$ is bijective.
 \item Let $\al \in \pi$. Then $H_\al^*$ is a free left (resp.\@ right) $H_\al$-module
       for the action defined, for any $f \in H_\al^*$ and $a,x \in H_\al$, by
       $$
       (a \rightharpoonup f)(x)=f(xa)
        \quad \text{(resp.} \quad (f\leftharpoonup a)(x)= f(a x) \text{).}
       $$
       Its rank is 1 if $H_\al \neq 0$ and $0$ otherwise.
       Moreover, if $\lambda=(\lambda_\be)_{\be \in \pi}$ is a non-zero left
       (resp.\@ right) \p-integral for $H$, then $\la$ is a basis vector
       for $H^*_\al$.
\end{enumerate}
\end{corollary}
\begin{proof}
To show Part (a), let $\al \in \pi$. By Lemma~\ref{Sinj} and Theorem~\ref{piintegexists}, $S_\al:H_\al \to H_\ali$
and $S_\ali: H_\ali \to H_\al$ are injective. Thus $\dim H_\al= \dim H_\ali $ and so $S_\al$ is bijective. To show
Part~(b), let $\lambda=( \la )_{\al \in \pi}$ be a non-zero left \p-integral for $H$ and fix $\al \in \pi$. If
$H_\al=0$, then the result is obvious. Let us suppose that $H_\al \neq 0$. Recall that $ H^\pc_\ali=H_\al^*$ and
$f_\ali : (H^*)^{\co H}_\al \otimes H_\ali \to H_\al^*$ defined by $f \otimes h \mapsto S_\ali(h)\rightharpoonup
f$ is an isomorphism. Since $0 \neq \la \in (H^*)^{\co H}_\al$, $\dim (H^*)^{\co H}_\al=1$, and $S_\ali$ is
bijective, the map $H_\al \to H^*_\al$ defined by $h \mapsto h \rightharpoonup \la$ is an isomorphism. Thus
$(H_\al^*,\rightharpoonup)$ is a free left $H_\al$-module of rank 1 with vector basis $\la$. Using $H^{\opcop}$
(see \S\ref{opcoop}), one easily deduces the right case.
\end{proof}

\section{The distinguished $\pi$-grouplike element}\label{s:disting}
In this section, we extend the notion of a grouplike element of a Hopf algebra to the setting of a Hopf
\p-coalgebra. We show that a \p-grouplike element is distinguished in a finite dimensional Hopf \p-coalgebra and
we study its relations with the \p-integrals. As a corollary, for any $\al \in \pi$ of finite order, we give an
upper bound for the (finite) order of $S_\ali S_\al$.

\subsection{$\pi$-grouplike elements}
A \emph{\p-grouplike element} of a Hopf \p-coalgebra $H$ is a family $g=(g_\al)_{\al \in \pi} \in \Pi_{\al \in
\pi} H_\al$ such that $\cp{\al}{\be}(g_{\al \be})=g_\al \otimes g_\be$ for any $\al,\be \in \pi$ and
$\epsilon(g_1)=1_\kk$ (or equivalently $g_1 \neq 0$). Note that $g_1$ is then a (usual) grouplike element of the
Hopf algebra $H_1$.

One easily checks that the set $G(H)$ of \p-grouplike elements of $H$ is a group (with respect to the
multiplication and unit of the product monoid $\Pi_{\al \in \pi} H_\al$) and if $g=(g_\al)_{\al \in \pi} \in
G(H)$, then $g^{-1}=(S_\ali(g_\ali))_{\al \in \pi}$.

Remark that the group $\homo(\pi,\kk^*)$ acts on $G(H)$ by $\phi g=(\phi(\al) g_\al) _{\al \in \pi}$ for any
$g=(g_\al)_{\al \in \pi} \in G(H)$ and $\phi \in \homo(\pi,\kk^*)$.
\begin{lemma}\label{distinteg}
Let $H$ be a finite dimensional Hopf \p-coalgebra. Then there exists a unique \p-grouplike element $g=(g_\al)_{\al
\in \pi}$ of $H$ such that $
 (\id_{H_\al} \otimes \lb){\cp{\al}{\be}}= \lab \, g_\al
 $ for any right \p-integral $\lambda=(
\la)_{\al \in \pi}$ and all $\al,\be \in \pi$.
\end{lemma}
The \p-grouplike element $g=(g_\al)_{\al \in \pi}$ of the previous lemma is called the \emph{distinguished
\p-grouplike element of} $H$. Note that $g_1$ is the (usual) distinguished grouplike element of the Hopf algebra
$H_1$.
\begin{proof}
Let $\lambda=(\la)_{\al \in\pi}$ be a non-zero right \p-integral for $H$. Let $\ga \in \pi$. For any  $\varphi \in
H^*_\ga$, $( \varphi \lambda_{\gai \al} )_{\al \in \pi}$ is a right \p-integral for $H$ and thus, by
Theorem~\ref{piintegexists}, there exists a unique $k_\varphi \in \kk$ such that $\varphi \lambda_{\gai \al}=
k_\varphi \la$ for all $\al \in \pi$. Now $(\varphi \mapsto k_\varphi) \in H^{**}_\ga \cong H_\ga$ ($\dim H_\ga <
+\infty$). Thus there exists a unique $g_\ga \in H_\ga$ such that $\varphi \lambda_{\gai \al}=\varphi(g_\ga) \la$
for any $\al \in \pi$ and $\varphi \in H^*_\ga$. Then $\varphi \lb =\varphi(g_\al) \lab$ for any $\al,\be \in \pi$
and $\varphi \in H_\al^*$ and hence $(\id_{H_\al} \otimes \lb)\cp{\al}{\be}=\lab \, g_\al$ for all $\al, \be \in
\pi$. Let $\al,\be \in \pi$. If $H_{\al\be}=0$, then either $H_\al=0$ or $H_\be=0$ (by Corollary~\ref{subgroup})
and so $\cp{\al}{\be}(g_{\al\be})=0=g_\al \otimes g_\be$. If $H_{\al \be} \neq 0$, then, for any $\varphi \in
H^*_\al$ and $\psi \in H^*_\be$, $k_{\varphi \psi} \lab= (\varphi \psi) \lu=\varphi(\psi \lu)=k_\psi \varphi \lb=
k_\varphi k_\psi \lab$ and thus $k_{\varphi \psi}=k_\varphi k_\psi$ (since $\lab \neq 0$ by
Lemma~\ref{intnonzero}), that is $\cp{\al}{\be}(g_{\al \be})=g_\al \otimes g_\be$. Moreover $\epsilon(g_1)
\lu=\epsilon \lu= (\epsilon \otimes \lu) \cp{1}{1}= \lu$ and so $\epsilon(g_1)=1$ (since $\lu \neq 0$ by
Lemma~\ref{intnonzero}). Then $g=(g_\al)_{\al\in \pi}$ is a \p-grouplike element of $H$. Since all the right
\p-integrals for $H$ are scalar multiple of $\lambda$, the ``existence" part of the lemma is demonstrated. Let us
now show the uniqueness of $g$. Suppose that $h=(h_\al)_{\al \in \pi}$ is another such \p-grouplike element of
$H$. Let $\lambda=( \la )_{\al \in \pi}$ be a non-zero right \p-integral for $H$. Fix $\al \in \pi$. If $H_\al=0$,
then $h_\al=0=g_\al$. If $H_\al \neq 0$, then $\la \neq 0$ (by Lemma~\ref{intnonzero}) and so there exists $a \in
H_\al$ such that $\la(a)=1$. Therefore $g_\al=\la(a) g_\al=(\id_{H_1} \otimes \la)\cp{1}{\al}
(a)=\la(a)h_\al=h_\al$. This completes the proof of the lemma.
\end{proof}

\subsection{The distinguished $\pi$-grouplike element and $\pi$-integrals} Throughout this subsection,
$H=\{H_\al\}_{\al \in \pi}$ will denote a finite dimensional Hopf \p-coalgebra.

Since $H_1$ is a finite dimensional Hopf algebra, there exists (see e.g. \cite{Rad3}) a unique algebra morphism
$\nu: H_1 \to \kk$ such that if $\Lambda$ is a left integral for $H_1$, then $\Lambda x=\nu(x)\Lambda$ for all $x
\in H_1$. This morphism is a grouplike element of the Hopf algebra $H_1^*$, called the distinguished grouplike
element of $H_1^*$. In particular, it is invertible in $H_1^*$ and its inverse $\nu^{-1}$ is also an algebra
morphism and verifies that if $\Lambda$ is a right integral for $H_1$, then $x \Lambda =\nu^{-1}(x)\Lambda$ for
all $x \in H_1$.

For all $\al \in \pi$, we define a left and a right $H_1^*$-action on $H_\al$ by setting, for any $f \in H_1^*$
and $a \in H_\al$,
  $$
  f \rightharpoonup
  a =\aua f(\adu) \quad \text{and} \quad a \leftharpoonup
  f =f(\auu) \ada.
  $$
The next theorem generalizes \cite[Theorem 3]{Rad1}. It is used in Section~\ref{s:pitrace} to show the existence
of traces.
\begin{theorem}\label{intfortrace}
Let $\lambda=( \la )_{\al \in \pi}$ be a right \p-integral for $H$, $g=(g_\al)_{\al \in \pi}$ be the distinguished
\p-grouplike element of $H$, and $\nu$ be the distinguished grouplike element of $H_1^*$. Then, for any $\al \in
\pi$ and $x,y \in H_\al$,
\begin{enumerate}
\renewcommand{\labelenumi}{{\rm (\alph{enumi})}}
 \item $\la(xy)=\la(S_\ali S_\al(y \leftharpoonup \nu)\, x)$;
 \item $\la(xy)=\la(y \,S_\ali S_\al(\nu^{-1} \rightharpoonup g^{-1}_\al x g_\al ))$;
 \item $\lai(S_\al(x))=\la(g_\al x)$.
\end{enumerate}
\end{theorem}
Before proving Theorem~\ref{intfortrace}, we establish the following lemma.

\begin{lemma} \label{lemuni}
Let $\lambda=( \la )_{\al \in \pi}$ be a right \p-integral for $H$, $\al \in \pi$, and $a \in H_\al$.
\begin{enumerate}
\renewcommand{\labelenumi}{{\rm (\alph{enumi})}}
 \item If $\Lambda$ is a right integral for $H_1$ such that
       $\lu(\Lambda)=1$, then $$S_\al(a)=\la (\Lambda_{(1,\al)} a) \, \Lambda_{(2,\ali)};$$
 \item If $\Lambda$ is a left integral for $H_1$ such that
       $\lu(\Lambda)=1$, then $$S_\ali^{-1}(a)=\la (a \Lua) \, \Ldai.$$
\end{enumerate}
\end{lemma}
\begin{proof}
To show Part (a), let $\al \in \pi$. Define $f \in H_\al ^*$ by $f(x)=\la( \Lua x) \, \Ldai$ for any $x \in
H_\al$. If $*$ denotes the product in the convolution algebra $\conv(H,H_\ali)$ (see \S\ref{convo}), then, for any
$x \in H_1$,
\begin{eqnarray*}
(f * \id_{H_\ali})(x)
   & = & \la (\Lua \xua) \, \Ldai \xdai \\
   & = & \la ((\Lambda x)_{(1,\al)}) \, (\Lambda x)_{(2,\ali)}
         \text{ \quad by \eqref{hopfmor}}\\
   & = & \epsilon(x) \, \la ( \Lua) \, \Ldai
         \text{ \quad since $\Lambda$ is a right integral for
         $H_1$}\\
   & = & \epsilon(x) \, \lu(\Lambda) \, 1_\ali
         \text{ \quad by \eqref{defiint}}\\
   & = & \epsilon(x) \, 1_\ali \text{ \quad since
         $\lu(\Lambda)=1$.}
\end{eqnarray*}
Therefore, since $\id_{H_\ali}$ is invertible in $\conv(H,H_\ali)$ with inverse $S_\al$, we have that $f=S_\al$,
that is $S_\al(a)=\la (\Lambda_{(1,\al)} a) \, \Lambda_{(2,\ali)}$ for all $a \in H_\al$. Part (b) can be deduced
from Part (a) using the Hopf \p-coalgebra $H^\opp$ (see \S\ref{op}).
\end{proof}

\begin{proof}[Proof of Theorem~\ref{intfortrace}]
We use the same arguments as in the proof of \cite[Theorem 3]{Rad1}, even if we cannot use the duality (since the
notion a Hopf \p-coalgebra is not self dual). We can assume that $\lambda$ is a non-zero right \p-integral
(otherwise the result is obvious). To show Part (a), let $\al \in \pi$ and $x,y \in H_\al$. Since $\lu$ is a
non-zero right integral for the Hopf algebra $H_1^*$, there exists a left integral $\Lambda$ for $H_1$ such that
$\lu(\Lambda)=\lu(S_1(\Lambda))=1$ (cf \cite[Proposition 1] {Rad1}). By Lemma~\ref{lemuni}(b) for $a=S_\ali
S_\al(y \leftharpoonup \nu)$, we have that
\begin{equation}\label{equni1}
S_\al(y \leftharpoonup \nu) = \la (S_\ali S_\al (y \leftharpoonup \nu) \, \Lua) \, \Ldai.
\end{equation}
It is easy to verify that $( \nu^{-1} \lambda_\ga )_{\ga \in \pi}$ is a right \p-integral for $H$ and $\Lambda
\leftharpoonup\nu $ is a right integral for $H_1$ such that $(\nu^{-1} \lu)(\Lambda \leftharpoonup \nu)=1$. Thus
Lemma~\ref{lemuni}(a) for $a=y \leftharpoonup \nu$ gives  that
\begin{eqnarray*}
 S_\al( y \leftharpoonup \nu)
   & = & (\nu^{-1} \la) ((\Lambda \leftharpoonup
         \nu)_{(1,\al)} \, (y\leftharpoonup \nu) ) \, (\Lambda \leftharpoonup
         \nu)_{(2,\ali)} \\
   & = & (\nu^{-1} \la) ((\Lambda_{(1,\al)}y) \leftharpoonup
         \nu) \, \Lambda_{(2,\ali)} \text{ \quad by \eqref{hopfmor}}\\
   & = & \la (((\Lambda_{(1,\al)}y) \leftharpoonup
         \nu)\leftharpoonup \nu^{-1}) \, \Lambda_{(2,\ali)} \\
   & = & \la ((\Lambda_{(1,\al)}y) \leftharpoonup
         \epsilon) \, \Lambda_{(2,\ali)} \\
   & = & \la( \Lambda_{(1,\al)} y)
         \, \Lambda_{(2,\ali)} \text{ \quad by \eqref{counit}.}
\end{eqnarray*}
Hence, comparing with (\ref{equni1}), we obtain
\begin{equation}\label{equni2}
\la( \Lambda_{(1,\al)} \, y) \, \Lambda_{(2,\ali)}= \la (S_\ali S_\al (y \leftharpoonup \nu) \, \Lua) \, \Ldai.
\end{equation}
Now $(\lambda_\ga S_\gai )_{\ga \in \pi}$ is a right \p-integral for $H^\cop$ and $\Lambda$ is a left integral for
$H_1^\cop$ such that $(\lu S_1)(\Lambda)=1$. Thus $(S^\cop_\ali)^{-1}(S^{-1}_\ali(x))=\la S_\ali(S_\ali^{-1}(x)
\Ldai) \, \Lua$ (by applying Lemma~\ref{lemuni}(b) for $a=S_\ali^{-1}(x) \in H^\cop_\al$), that is
\begin{equation}\label{equni3}
x=\Lua \la (S_\ali(\Ldai) \,x).
\end{equation}
Finally evaluating (\ref{equni2}) with $\la(S_\ali(\cdot)x)$ and using (\ref{equni3}) gives that $\la (xy)=\la
(S_\ali S_\al (y \leftharpoonup \nu) \, x)$.

To show Part (b), let $\al \in \pi$ and $a,b \in H_\al$. For any $\ga \in \pi$, let us define $\phi_\ga \in
(H_\ga^\opcop)^*$ by $\phi_\ga(x)=\lambda_\gai(g_\gai x)$ for all $x \in H_\ga^\opcop$. Using
Lemma~\ref{distinteg}, one easily checks that $\phi=(\phi_\ga)_{\ga \in \pi}$ is a right \p-integral for
$H^\opcop$. Let us denote by $\times^{\opp}$ the multiplication of $H_\ali^\opcop$ and by $\leftharpoonup^\cop$
the right action of $(H^\opcop_1)^*$ on $H_\ali^\opcop$ defined by $h\leftharpoonup^\cop f=(f \otimes \id)
\cp{1}{\ali}^\cop$. Then, since $\nu^{-1}$ is the distinguished grouplike element of $(H^\opcop_1)^*$, Part (a)
with $x=g^{-1}_\al b$ and $y=g_\al^{-1} a g_\al$ gives that $\phi_\ali(y\times^{\opp} x)= \phi_\ali(S_\al^\opcop
S_\ali^\opcop(y \leftharpoonup^\cop \nu^{-1}) \times^\opp x)$, that is $\la(ab)=\la(b\, S_\ali S_\al (\nu^{-1}
\rightharpoonup g_\al^{-1} a g_\al))$.

Let us show Part (c). For any $\al \in \pi$, define $\phi_\al \in H_\al^*$ by $\phi_\al(x)=\lambda_\al(g_\al x)$
for all $x \in H_\al$. Since $(\phi_\al )_{\al \in \pi}$ and $( \lai S_\al )_{\al \in \pi}$ are left \p-integrals
for $H$ which are non-zero (because $\lambda$ is non-zero, $g_\al$ is invertible and $S_\al$ is bijective), there
exists $k \in \kk$ such that $\phi_\al=k \lai S_\al$ for all $\al \in \pi$ (by Theorem~\ref{piintegexists}). As
above, let $\Lambda$ be a left integral  for $H_1$ such that $\lu(\Lambda)=\lu(S_1(\Lambda))=1$. Recall that
$\epsilon(g_1)=1$. Thus $1=\lu(\Lambda)=\lu(\epsilon(g_1) \Lambda)=\lu(g_1 \Lambda) =k \lu(S_1(\Lambda))=k$. Hence
$\lai S_\al=\phi_\al$ for all $\al \in \pi$, that is $\lai(S_\al(x))=\la(g_\al x)$ for all $\al \in \pi$ and $x
\in H_\al$. This completes the proof of the theorem.
\end{proof}

The following corollary will be used later to relate the distinguished grouplike element of a finite dimensional
quasitriangular Hopf \p-coalgebra to the $R$-matrix.

\begin{corollary}\label{chatquipue}
Let $\Lambda$ be a left integral for $H_1$ and $g=(g_\al)_{\al \in \pi}$ be the distinguished \p-grouplike element
of $H$. Then $ \Lambda_{(1,\al)} \otimes \Lambda_{(2,\ali)} = S_\ali S_\al(\Lambda_{(2,\al)}) g_\al \otimes
\Lambda_{(1,\ali)} $ for all $\al \in \pi$.
\end{corollary}
\begin{proof}
We can suppose that $\Lambda\neq 0$. Let $\al \in \pi$. Remark that it suffices to show that, for all $f \in
H^*_\ali$,
\begin{equation}\label{groscul}
  f(\Lambda_{(2,\ali)}) \, \Lambda_{(1,\al)}
  =f( \Lambda_{(1,\ali)}) \, S_\ali
  S_\al(\Lambda_{(2,\al)})g_\al.
\end{equation}
Fix $f \in H^*_\ali$. Let $\lambda=(\lambda_\ga )_{\ga \in \pi}$ be a non-zero right \p-integral for $H$ (see
Theorem~\ref{piintegexists}). By multiplying $\lambda$ by some (non-zero) scalar, we can assume that
$\lu(\Lambda)=\lu(S_1(\Lambda))=1$. By Corollary~\ref{corinteg}(b), there exists $a \in H_\ali$ such that
$f(x)=\lai(a x)$ for all $x \in H_\ali$. By Lemma~\ref{lemuni}(b), $S_\ali(a)=\lai(a \Lambda_{(1,\ali)}) \, S_\ali
S_\al(\Lambda_{(2,\al)})$. Thus
\begin{equation}\label{eqo1}
 S_\ali(a) g_\al=f(\Lambda_{(1,\ali)}) \, S_\ali S_\al(\Lambda_{(2,\al)})
 g_\al.
\end{equation}
Since $(\lambda_\ga S_\gai )_{\ga \in \pi}$ is a right \p-integral for $H^{\opcop}$ and $\Lambda$ is a right
integral for $H_1^{\opcop}$ such that $(\lu S_1)(\Lambda)=1$, Lemma~\ref{lemuni}(a)  applied to $H^\opcop$ gives
that
 $$
 S_\ali(a)=\la S_\ali(a \Lambda_{(2,\ali)})
 \Lambda_{(1,\al)}.
 $$
Then
\begin{eqnarray}
 S_\ali(a) g_\al
   & = & \la S_\ali(a \Lambda_{(2,\ali)}) \, \Lambda_{(1,\al)} g_\al \nonumber\\
   & = & \lai(g_\ali a \Lambda_{(2,\ali)}) \, \Lambda_{(1,\al)} g_\al
         \text{ \quad by Theorem~\ref{intfortrace}(c).} \label{eqo2}
\end{eqnarray}
Now, since $\Lambda$ is left integral for $H_1$,
$$
\Lambda_{(1,\al)} g_\al \otimes \Lambda_{(2,\ali)} g_\ali = \cp{\al}{\ali}(\Lambda g_1) = \nu(g_1) \,
\Lambda_{(1,\al)} \otimes \Lambda_{(2,\ali)}.
$$
Therefore
$$
\Lambda_{(1,\al)} g_\al \otimes  g_\ali a \Lambda_{(2,\ali)}
  = \Lambda_{(1,\al)} \otimes \nu(g_1) \, g_\ali a \Lambda_{(2,\ali)}
  g_\ali^{-1},
$$
and so, using \eqref{eqo2} and then Theorem~\ref{intfortrace}(a),
\begin{eqnarray*}
 S_\ali(a) g_\al
   & = & \lai(\nu(g_1) \, g_\ali a \Lambda_{(2,\ali)}
         g_\ali^{-1}) \Lambda_{(1,\al)} \\
   & = & \lai(\nu(g_1) \, S_\al S_\ali(g_\ali^{-1}\leftharpoonup \nu)
         g_\ali a  \Lambda_{(2,\ali)}) \, \Lambda_{(1,\al)}
\end{eqnarray*}
Now, since $g_\be^{-1}=S_\bei(g_\bei)$ for any $\be \in \pi$, $g^{-1}=(g_\be^{-1})_{\be \in \pi} \in G(H)$, and
$\nu$ is an algebra morphism,
  $$
  S_\al S_\ali(g_\ali^{-1}\leftharpoonup \nu) =\nu(g_1^{-1}) \, S_\al
  S_\ali(g_\ali^{-1})= \nu(g_1)^{-1} g_\ali^{-1} .
  $$
Thus $S_\ali(a) g_\al= \lai(a  \Lambda_{(2,\ali)}) \, \Lambda_{(1,\al)} = f(\Lambda_{(2,\ali)}) \,
\Lambda_{(1,\al)}$. Finally, comparing with (\ref{eqo1}), we get \eqref{groscul}.
\end{proof}

\subsection{The order of the antipode}
It is known that the order of the antipode of a finite dimensional Hopf algebra $A$ is finite (by \cite[Theorem 1]
{Rad3}) and divides $4 \dim A$ (by \cite[Proposition 3.1]{NZ2}). Let us apply this result to the setting of a Hopf
\p-coalgebra.

Let $H=\{H_\al\}_{\al \in \pi}$ be a finite dimensional Hopf \p-coalgebra with antipode $S=\{S_\al\}_{\al \in
\pi}$. Let $\al \in \pi$ of finite order $d$ and denote by $\langle \al \rangle$ the subgroup of $\pi$ generated
by $\al$. By considering the (finite dimensional) Hopf algebra $\oplus_{\be \in \langle \al \rangle} H_\be$
(coming from the Hopf $\langle \al \rangle$-coalgebra $\{H_\be\}_{\be \in \langle \al \rangle}$, as in
\S\ref{pifinite}), we obtain that the order of $S_\ali S_\al \in {\rm Aut_{Alg}}(H_\al)$ is finite and divides $2
\sum_{\be \in \langle \al \rangle} \dim H_\be$. As a corollary of Theorem~\ref{intfortrace}, we give another upper
bound for the order of $S_\ali S_\al$.

\begin{corollary}\label{cororderantipode}
Let $H=\{H_\al\}_{\al \in \pi}$ be a finite dimensional Hopf \p-coalgebra with antipode $S=\{S_\al\}_{\al \in
\pi}$. Then
\begin{enumerate}
\renewcommand{\labelenumi}{{\rm (\alph{enumi})}}
 \item If $\al \in \pi$ has a finite order d, then $(S_\ali
       S_\al)^{2 d \dim H_1}=\id_{H_\al}$;
 \item If $\al \in \pi$ has order $2$, then $S_\al^{8 \dim H_1}=\id_{H_\al}$.
\end{enumerate}
\end{corollary}
Before proving Corollary~\ref{cororderantipode}, we establish the following lemma.
\begin{lemma}\label{antipode4calc}
Let $H$ be a finite dimensional Hopf \p-coalgebra, $g=(g_\al)_{\al \in \pi}$ be the distinguished \p-grouplike
element of $H$, and $\nu$ be the distinguished grouplike element of $H_1^*$. Then $(S_\ali S_\al)^2(x)= g_\al (\nu
\rightharpoonup x \leftharpoonup \nu^{-1}) g^{-1}_\al$ for all $\al \in \pi$ and $x \in H_\al$.
\end{lemma}
\begin{proof}
Let $\al \in \pi$ and $x,y \in H_\al$. If $H_\al=0$, then the result is obvious. Let us suppose that $H_\al \neq
0$. Let $\lambda=(\lambda_\ga)_{\ga \in \pi}$ be a non-zero right \p-integral for $H$. Then
\begin{eqnarray*}
\lefteqn{ \la(g_\al (\nu \rightharpoonup x \leftharpoonup \nu^{-1}) g^{-1}_\al y)}\\
  & = & \la (y S_\ali S_\al (\nu^{-1} \rightharpoonup
        g_\al^{-1} g_\al (\nu \rightharpoonup x \leftharpoonup
        \nu^{-1}) g^{-1}_\al g_\al))
        \text{ \quad by Theorem~\ref{intfortrace}(b)}\\
  & = & \la(y S_\ali S_\al (x \leftharpoonup \nu^{-1})) \\
  & = & \la(S_\ali S_\al(S_\ali S_\al (x \leftharpoonup \nu^{-1}) \leftharpoonup \nu
        )y) \text{ \quad by Theorem~\ref{intfortrace}(a)}\\
  & = & \la ((S_\al S_\ali)^2 (x\leftharpoonup \nu^{-1} \leftharpoonup
        \nu)y) \text{ \quad since $S_\ali S_\al$ is comultiplicative}\\
  & = & \la ((S_\al S_\ali)^2 (x) y).
\end{eqnarray*}
Now, by Corollary~\ref{corinteg}(b), $H_\al^*$ is a free right $H_\al$-module of rank 1 for the action defined by
$(f\vartriangleleft a)(x)=f(a x)$ for any $f \in H_\al^*$ and $a,x \in H_\al$. Moreover $\la$ is a basis vector of
$H_\al^*$. Thus, since the above computation says that
  $$
  \la \vartriangleleft g_\al (\nu \rightharpoonup x \leftharpoonup \nu^{-1}) g^{-1}_\al
  = \la \vartriangleleft (S_\al S_\ali)^2 (x),
  $$
we conclude that $(S_\ali S_\al)^2(x)= g_\al (\nu \rightharpoonup x \leftharpoonup \nu^{-1}) g^{-1}_\al$.
\end{proof}

\begin{proof}[Proof of Corollary~\ref{cororderantipode}]
To show Part (a), let $\al \in \pi$ of finite order $d$. Consider the distinguished \p-grouplike element
$g=(g_\al)_{\al \in \pi}$ of $H$ and the distinguished grouplike element $\nu$ of $H_1^*$. Using
Lemma~\ref{antipode4calc}, one easily shows by induction that
\begin{equation}\label{equarec}
  (S_\ali S_\al)^{2l}(x)=g_\al^l(\nu^l \rightharpoonup x \leftharpoonup \nu^{-l}) g_\al^{-l}
\end{equation}
for all $x \in H_\al$ and $l \in \mathbb{N}$. Recall that the order of a grouplike element of a finite dimensional
Hopf algebra $A$ is finite and divides $\dim A$ (see \cite[Theorem 2.2]{NZ2}). Therefore $g_1$ has a finite order
which divides $\dim H_1$ and $\nu$ has a finite order which divides $\dim H_1^*= \dim H_1$. Now $$g_\al ^{d \dim
H_1}=g^{\dim H_1}_{1 \, (1,\al)} \cdots g^{\dim H_1}_{1 \, (d,\al)} = 1_{1 \, (1,\al)} \cdots 1_{1 \, (d,\al)}
=1_\al^d =1_\al.$$ Then, for all $x \in H_\al$, by \eqref{equarec},
\begin{eqnarray*}
(S_\ali S_\al)^{2d \dim H_1}(x)
  & = & g_\al^{d \dim H_1}(\nu^{d \dim H_1} \rightharpoonup x \leftharpoonup
        \nu^{-d \dim H_1}) g_\al^{-d \dim H_1}\\
  & = & 1_\al(\epsilon \rightharpoonup x \leftharpoonup
        \epsilon)1_\al \; = \; x.
\end{eqnarray*}
Hence $(S_\ali S_\al)^{2d \dim H_1}=\id_{H_\al}$. Part (b) is Part (a) for $d=2$, since in this case $S_\al$ is an
endomorphism of $H_\al$.
\end{proof}

\section{Semisimplicity and cosemisimplicity}\label{semicosemi}
In this section, we define the semisimplicity and the cosemisimplicity for Hopf \p-coalgebras, and we give
criteria for a Hopf \p-coalgebra to be semisimple (resp.\@ cosemisimple).

\subsection{Semisimple Hopf $\pi$-coalgebras}
A Hopf \p-coalgebra $H=\{H_\al\}_{\al \in \pi}$ is said to be \emph{semisimple} if each algebra $H_\al$ is
semisimple.

Note that, since any infinite dimensional Hopf algebra (over a field) is never semisimple (see \cite[Corollary
2.7]{sweed2}), a necessary condition for $H$ to be semisimple is that $H_1$ is finite dimensional.

\begin{lemma}\label{semisimplecrit}
Let $H=\{H_\al\}_{\al \in \pi}$ be a finite dimensional Hopf \p-coalgebra. Then $H$ is semisimple if and only if
$H_1$ is semisimple.
\end{lemma}
\begin{proof}
We have to show that if $H_1$ is semisimple then $H$ is semisimple. Suppose that $H_1$ is semisimple and fix $\al
\in \pi$. Since $H_\al$ is a finite dimensional algebra, it suffices to show that all left $H_\al$-modules are
completely reducible. Thus let $M$ be a left $H_\al$-module and $N$ be a submodule of $M$. Since $H_1$ is a finite
dimensional semisimple Hopf algebra, there exists a left integral $\Lambda$ for $H_1$ such that
$\epsilon(\Lambda)=1$ (cf \cite[Theorem 5.1.8]{sweed}). Let $p:M \to N$ be any $\Bbbk$-linear projection which is
the identity on $N$. Let $P: M \to N$ be the $\Bbbk$-linear map defined by $$P(m)=\Lua \cdot p( S_\ali(\Ldai)
\cdot m)$$ for any $m\in M$, where $\cdot$ denotes the action of $H_\al$ on $M$. The map $P$ is the identity on
$N$ since, for any $n \in N$,
\begin{multline*}
  P(n)=\Lua \cdot p( S_\ali(\Ldai) \cdot n) \\ =\Lua \cdot ( S_\ali(\Ldai) \cdot n)=
  (\Lua S_\ali(\Ldai)) \cdot n =
  \epsilon(\Lambda) 1_\al \cdot n = n.
\end{multline*}
Let $h \in H_\al$. Using \eqref{counit} and the fact that $\Lambda$ is a left integral for $H_1$, we have
\begin{eqnarray*}
 \Lua \otimes \Ldai \otimes  h
   & = & \cp{\al}{\ali}(\epsilon(\huu) \,\Lambda) \otimes  \hda\\
   & = & \cp{\al}{\ali}(\huu \Lambda) \otimes  \hda \\
   & = & \hua \Lua \otimes \hdai \Ldai \otimes \hta,
\end{eqnarray*}
and so
\begin{eqnarray*}
\lefteqn{\Lua \otimes S_\ali(\Ldai) h \qquad \qquad}\\
   & = & \hua \Lua \otimes S_\ali(\hdai \Ldai) \hta \\
   & = & \hua \Lua \otimes S_\ali(\Ldai) S_\ali(\hdai) \hta \text{ \quad by Lemma~\ref{antipodepptes}(c)} \\
   & = & \hua \epsilon(\hdu) \, \Lua \otimes S_\ali(\Ldai) 1_\al
         \text{ \quad by \eqref{antipode}}\\
   & = & h \Lua \otimes S_\ali(\Ldai)
        \text{ \quad by \eqref{counit}}.
\end{eqnarray*}
Therefore, for all $h \in H_\al$ and $m \in M$,
\begin{multline*}
  P(h \cdot m)= \Lua \cdot p( S_\ali(\Ldai)h \cdot m) \\ = h \Lua \cdot
  p(S_\ali(\Ldai) \cdot m)=h \cdot P(m).
\end{multline*}
Hence $P$ is $H_\al$-linear and $\ker P$ is a $H_\al$-supplement of $N$ in $M$.
\end{proof}

\subsection{Cosemisimple $\pi$-comodules and $\pi$-coalgebras}
Let $C$ be a \p-coalgebra and $M$ be a right \p-comodule over $C$. If $\{N^i=\{N^i_\al\}_{\al \in \pi}\}_{i \in
I}$ is a family of \p-subcomodules of $M$, we define their \emph{sum } by $\{\sum_{i \in I} N^i_\al \}_{\al \in
\pi}$. It is easy to see that it is a \p-subcomodule of $M$. We denote it by $\sum_{i \in I} N^i$. This sum is
said to be \emph{direct} provided $\sum_{i \in I} N^i_\al $ is a direct sum for all $\al \in \pi$. In this case
$\sum_{i \in I} N^i$ will be denoted by $\oplus_{i \in I} N^i$.

A right \p-comodule $M=\{M_\al\}_{\al \in \pi}$ is said to be \emph{simple} if it is \emph{non-zero} (i.e., $M_\al
\neq 0$ for some $\al \in \pi$) and if it has no \p-subcomodules other than $0=\{0\}_{\al \in \pi}$ and itself.

\begin{lemma}\label{equivcodecomp}
Let $M$ be a right \p-comodule over a \p-coalgebra $C$. The following conditions are equivalent:
\begin{enumerate}
\renewcommand{\labelenumi}{{\rm (\alph{enumi})}}
\item $M$ is a sum of a family of simple \p-subcomodules;
\item $M$ is a direct sum of a family of simple \p-subcomodules;
\item Every \p-subcomodule $N$ of $M$ is a direct
      summand, i.e., there exists a \p-subcomodule $N'$ of $M$
      such that $M=N \oplus N'$.
\end{enumerate}
\end{lemma}
\begin{proof}
Let us show Condition (a) $\Rightarrow$ Condition (b). Suppose that $M=\sum_{i \in I} M^i$ is a sum of simple
\p-submodules. Let $J$ be a maximal subset of $I$ such that $\sum_{j \in J} M^j$ is direct. Let us show that this
sum is in fact equal to $M$. It suffices to prove that each $M^i$ ($i \in I$) is contained in this sum. The
intersection of our sum with $M^i$ is a \p-subcomodule of $M^i$, thus equal to 0 or $M^i$. If it is equal to $0$,
then $J$ is not maximal since we can adjoin $i$ to it. Hence $M^i$ is contained in the sum.

To show Condition (b) $\Rightarrow$ Condition (c), suppose that $M=\oplus_{i \in I} M^i$ and let $N$ be a
\p-subcomodule of $M$. Let $J$ be a maximal subset of $I$ such that the sum $N + \oplus_{j \in J} M^j$ is direct.
The same reasoning as before shows this sum is equal to $M$.

Let us show Condition (c) $\Rightarrow$ Condition (a). Let $N$ be the \p-subcomodule of $M$ defined as the sum of
all simple \p-subcomodules of $M$. Suppose that $M \neq N$. Then $M=N \oplus F$ where $F$ is a non-zero
\p-subcomodule of $M$. Let us show that there exists a simple \p-subcomodule of $F$, contradicting the definition
of $N$. By Theorem~\ref{theorational}(a), $\oF=\oplus_{\al \in \pi} \oF_\al$ (where $\oF_\al=F_\ali$) is a
rational \p-graded left $C^*$-module which is non-zero. Let $v \in \oF$, $v\neq 0$. The kernel of the morphism of
\p-graded left $C^*$-modules $C^* \to C^* v$ is a \p-graded left ideal $J \neq C^*$. Therefore $J$ is contained in
a maximal \p-graded left ideal $I\neq C^*$ (by Zorn's lemma). Then $I/J$ is a maximal \p-graded left
$C^*$-submodule of $C^*/J$ (not equal to $C^*/J$), and hence $Iv$ is a maximal \p-graded $C^*$-submodule of
$C^*v$, not equal to $C^*v$ (corresponding to $I/J$ under the \p-graded isomorphism $C^*/J \to C^*v$). Moreover it
is rational since it is a submodule of the rational module $\oF$ (see Theorem~\ref{theorational}(b)). So we can
consider the \p-subcomodule $\overline{Iv}$ of $M$ (see Lemma~\ref{rational2}). Write then $M=\overline{Iv} \oplus
L$ where $L$ is \p-subcomodule of $M$. Therefore $\oM=Iv \oplus \overline{L}$ and so $C^*v=Iv \oplus (\overline{L}
\cap C^*v)$. Now, since $Iv$ is a maximal \p-graded $C^*$-submodule of $C^*v$ (not equal to $C^*v$), we have that
$\overline{L} \cap C^*v$ is a non-zero \p-graded $C^*$-submodule of $\oF$ which does not contain any \p-graded
submodule other than $0$ and itself. Moreover $\overline{L} \cap C^*v$ is rational since it is a \p-graded
$C^*$-submodule of the rational \p-graded $C^*$-module $\oF$ (see Theorem~\ref{theorational}(b)). Finally
$\overline{\overline{L} \cap C^*v}$ is a simple \p-subcomodule of $F$.
\end{proof}

A right \p-comodule satisfying the equivalent conditions of Lemma~\ref{equivcodecomp} is said to be {\it
cosemisimple}. A \p-coalgebra is called \emph{cosemisimple} if it is cosemisimple as a right \p-comodule over
itself (with comultiplication as structure maps).

When $\pi=1$, one recovers the usual notions of cosemisimple comodules and coalgebras.

When $\pi$ is finite, a \p-coalgebra $C=\{C_\al\}_{\al \in \pi}$ is cosemisimple if and only if the \p-graded
coalgebra $\tilde{C}=\oplus_{\al\in \pi} C_\al$ (defined as in \S\ref{pifinite}) is \emph{graded-cosemisimple}
(i.e., is a direct sum of simple $\pi$-graded right comodules).

\begin{lemma}
Every \p-subcomodule or quotient of a cosemisimple right \p-comodule is cosemisimple.
\end{lemma}
\begin{proof}
Let $N$ be a \p-subcomodule of a cosemisimple right \p-comodule $M$. Let $F$ be the sum of all simple
\p-subcomodules of $N$ and write $M=F \oplus F'$. Therefore $N=F \oplus(F' \cap N)$. If $F'\cap N \neq 0$, it
contains a simple \p-subcomodule (see the demonstration of Lemma~\ref{equivcodecomp}). Thus $F' \cap N=0$ and
$N=F$, which is cosemisimple. Now write $M=N \oplus N'$. $N'$ is a sum of simple \p-subcomodules (it is a
\p-subcomodule of $M$ and thus cosemisimple) and the canonical projection $M \to M/N$ induces a \p-comodule
isomorphism between $N'$ onto $M/N$. Hence $M/N$ is cosemisimple.
\end{proof}

\subsection{Cosemisimple Hopf $\pi$-coalgebras} A Hopf \p-coalgebra $H=\{H_\al\}_{\al \in \pi}$ is said to be {\it
cosemisimple} if it is cosemisimple as a \p-coalgebra. A right \p-comodule $M=\{M_\al\}_{\al \in \pi}$ over $H$ is
said to be \emph{reduced} if, for all $\al \in \pi$,  $M_\al=0$ whenever $H_\al=0$.

The next theorem is the Hopf \p-coalgebra version of the dual Maschke theorem (see \cite[\S 14.0.3]{sweed}).

\begin{theorem}\label{cosemicrit}
Let $H$ be a Hopf \p-coalgebra. The following conditions are equivalent:
\begin{enumerate}
\renewcommand{\labelenumi}{{\rm (\alph{enumi})}}
 \item Every reduced right \p-comodule over $H$ is cosemisimple;
 \item $H$ is cosemisimple;
 \item There exists a right \p-integral $\lambda=( \la )_{\al\in \pi}$ for $H$
       such that $\la(1_\al)=1$ for some $\al \in \pi$;
 \item There exists a right \p-integral $\lambda=( \la )_{\al\in \pi}$ for $H$
       such that $\la(1_\al)=1$ for all $\al \in \pi$ with $H_\al \neq 0$.
\end{enumerate}
\end{theorem}
\begin{proof}
Condition (a) implies trivially Condition (b). Moreover Condition (c) is equivalent to Condition (d). Indeed
Condition (d) implies Condition (c) since $H_1 \neq 0$ (by Corollary~\ref{subgroup}). Conversely, suppose that
$\be \in \pi$ is such that $\lb(1_\be)=1$. Let $\al \in \pi$ such that $H_\al \neq 0$. Then $\la(1_\al) \, 1_{\bei
\al}=(\lb \otimes \id_{H_{\bei \al}}) \cp{\be}{\bei \al} (1_\al)=\lb(1_\be) \, 1_{\bei \al}=1_{\bei \al}$. Now
$1_{\bei \al} \neq 0$ by Corollary~\ref{subgroup}. Hence $\la(1_\al)=1$.

Let us show that Condition (b) implies Condition (d). Consider $H$ as a a right \p-comodule over itself (with
comultiplication as structure maps). For any $\al \in \pi$, set $N_\al=\kk 1_\al$. Since the comultiplication is
unitary, $N$ is a \p-subcomodule of $H$. Therefore $N$ is a direct summand of $H$ (since $H$ is cosemisimple),
that is there exists a \p-comodule morphism $p=\{p_\al \}_{\al \in \pi}: H \to N$ such that
$p_{\al|N_\al}=\id_{N_\al}$ for all $\al \in \pi$. For any $\al \in \pi$, since $N_\al = \kk 1_\al$, there exists
a (unique) $\kk$-form $\la \in H_\al^*$ such that $p_\al(h)=\la(h) \, 1_\al$ for all $h \in H_\al$. Let us verify
that $\lambda=(\la)_{\al \in \pi}$ is a right \p-integral for $H$. Let $\al,\be \in \pi$. Since $p$ is a
\p-comodule morphism, we have that
\begin{equation}\label{eqcosem1}
 \lab \, 1_\al \otimes 1_\be = (\la 1_\al \otimes \id_{H_\be}) \cp{\al}{\be}.
\end{equation}
If $H_\al=0$, then either $H_\be=0$ or $H_{\al \be}=0$ (by Corollary~\ref{subgroup}) and so $ \lab \, 1_\be = 0 =
(\la \otimes \id_{H_\be}) \cp{\al}{\be}$. If $H_\al \neq 0$, then there exists $f \in H_\al^*$ such that
$f(1_\al)=1$ and, by applying $(f \otimes \id_{H_\be})$ to both sides of \eqref{eqcosem1}, we get that  $ \lab \,
1_\be = (\la \otimes \id_{H_\be}) \cp{\al}{\be}$. Therefore $\lambda$ is a right \p-integral for $H$. Finally, let
$\al \in \pi$ such that $H_\al \neq 0$. Then $\la(1_\al) 1_\al=p_\al(1_\al)=1_\al$ (since $1_\al \in N_\al$) and
so $\la(1_\al)=1$ (since $1_\al \neq 0$).

To show that Condition (d) implies Condition (a), let $M=\{M_\al \}_{\al \in \pi}$ be a reduced right \p-comodule
over $H$ with structure maps by $\rho=\{\rh{\al}{\be} \}_{\al, \be \in \pi}$ and $N=\{N_\al \}_{\al \in \pi}$ be a
\p-subcomodule of $M$. We have to show that $N$ is a direct summand of $M$ (see Lemma~\ref{equivcodecomp}). Define
$\delta_\al : H_\ali \otimes H_\al \to \kk$ by $\delta_\al(x \otimes y)=\la(S_\ali(x)y)$ for all $\al \in \pi$. We
first prove that, for any $\al,\be,\ga \in \pi$,
\begin{equation}\label{cosem11}
 (\id_{H_\be} \otimes \delta_{\al \be})(\cp{\be}{(\al \be)^{-1}}
 \otimes \id_{H_{\al \be}})=(\delta_\al \otimes
 \id_{H_\be})(\id_{H_\ali} \otimes \cp{\al}{\be}).
\end{equation}
Indeed, for any $x \in H_\ali$ and $y \in H_{\al \be}$,
\begin{eqnarray*}
\lefteqn{(\id_{H_\be} \otimes \delta_{\al \be})(\cp{\be}{(\al
\be)^{-1}} \otimes \id_{H_{\al \be}})(x \otimes y)}\\
   & = & \xub \lab(S_{(\al \be)^{-1}}(x_{(2,(\al \be)^{-1})})y)\\
   & = & \xub (\la \otimes \id_{H_\be})\cp{\al}{\be}
         (S_{(\al \be)^{-1}}(x_{(2,(\al \be)^{-1})})y) \text{\quad by \eqref{defiint}}\\
   & = & \xub S_\bei(\xdbi)
         \ydb \, \la(S_\ali(\xtai) \yua) \text{\quad by Lemma~\ref{antipodepptes}(c)}\\
   & = & \ydb \, \la(S_\ali(\epsilon(\xuu) \xdai) \yua) \text{\quad by \eqref{antipode}}\\
   & = & \la(S_\ali(x) \yua) \, \ydb \text{\quad by \eqref{counit}}\\
   & = & (\delta_\al \otimes \id_{H_\be})(\id_{H_\ali} \otimes
         \cp{\al}{\be})( x \otimes y).
\end{eqnarray*}
Let $q:M_1 \to N_1$ be any $\Bbbk$-linear projection and define, for all $\al \in \pi$, $$p_\al=(\id_{N_\al}
\otimes \delta_\al)(\rh{\al}{\ali} \circ q \otimes \id_{H_\al})\rh{1}{\al}: M_\al \to N_\al.$$ For any $\al,\be
\in \pi$, using \eqref{coasscomod} and \eqref{cosem11}, we have
\begin{eqnarray*}
\lefteqn{\rh{\al}{\be} p_{\al \be}}\\
  & = & \rh{\al}{\be} (\id_{N_{\al \be}} \otimes \delta_{\al \be})
        (\rh{\al \be}{(\al \be)^{-1}} \circ q  \otimes
        \id_{H_{\al\be}})
        \rh{1}{\al\be}\\
  & = & (\id_{N_\al}\otimes \id_{H_\be} \otimes \delta_{\al \be})
        ((\rh{\al}{\be} \otimes \id_{H_{(\al \be)^{-1}}})\rh{\al \be}{(\al \be)^{-1}}  \circ q \otimes
        \id_{H_{\al\be}})
        \rh{1}{\al\be}\\
  & = & (\id_{N_\al}\otimes \id_{H_\be} \otimes \delta_{\al \be})
        ((\id_{N_\al} \otimes \cp{\be}{(\al \be)^{-1}})\rh{\al}{\ali} \circ q \otimes
        \id_{H_{\al\be}})
        \rh{1}{\al\be}\\
  & = & (\id_{N_\al} \otimes (\id_{H_\be} \otimes \delta_{\al \be})
        (\cp{\be}{(\al \be)^{-1}} \otimes \id_{H_{\al\be}}) )( \rh{\al}{\ali} \circ q \otimes
        \id_{H_{\al\be}})
        \rh{1}{\al\be}\\
  & = & (\id_{N_\al} \otimes (\delta_\al \otimes
        \id_{H_\be})(\id_{H_\ali} \otimes \cp{\al}{\be}) )( \rh{\al}{\ali}\circ q  \otimes
        \id_{H_{\al\be}})
        \rh{1}{\al\be}\\
  & = & (\id_{N_\al} \otimes \delta_\al \otimes \id_{H_\be})
        (\rh{\al}{\ali} \circ q \otimes \id_{H_\al} \otimes \id_{H_\be})
        (\id_{M_1}
        \otimes \cp{\al}{\be}) \rh{1}{\al\be}\\
  & = & (\id_{N_\al} \otimes \delta_\al \otimes \id_{H_\be})
        (\rh{\al}{\ali} \circ q \otimes \id_{H_\al} \otimes \id_{H_\be})
        (\rh{1}{\al} \otimes \id_{H_\be})
        \rh{\al}{\be}\\
  & = & (p_\al \otimes \id_{H_\be})\rh{\al}{\be}.
\end{eqnarray*}
Thus $p=\{p_\al\}_{\al \in \pi}$ is a \p-comodule morphism between $M$ and $N$. Let $\al \in \pi$ and $n \in
N_\al$. If $H_\al=0$, then $N_\al=0$ (since $M$ and thus $N$ is reduced) and so $p_\al(n)=0=n$. If $H_\al \neq 0$,
then
\begin{eqnarray*}
p_\al(n)
  & = & \nza \la(S_\ali(\nuai) \nda) \text{ \quad since $q_{|N_1}=\id_{N_1}$}\\
  & = & \nza \epsilon(\nuu) \la(1_\al) \text{ \quad by \eqref{antipode}}\\
  & = & n \text{ \quad by \eqref{counitcomod} and since $\la(1_\al)=1$.}
\end{eqnarray*}
Therefore $q$ is a \p-comodule projection of $M$ onto $N$ and consequently $N$ is a direct summand of $M$ (namely
$M=N \oplus \ker q$). This finishes the proof of the theorem.
\end{proof}

\begin{corollary}\label{equivcosemifd}
Let $H$ be a Hopf \p-coalgebra. Then
\begin{enumerate}
\renewcommand{\labelenumi}{{\rm (\alph{enumi})}}
 \item If $H$ is cosemisimple, then the Hopf algebra $H_1$ is cosemisimple;
 \item If $H$ is finite dimensional, then $H$ is cosemisimple if and only if
       $H_1$ is cosemisimple.
\end{enumerate}
\end{corollary}
\begin{proof}
To show Part (a), suppose that $H$ is cosemisimple. By Theorem~\ref{cosemicrit} and Corollary~\ref{subgroup},
there exists a right \p-integral $\lambda=( \la )_{\al\in \pi}$ for $H$ such that $\lu(1_1)=1$. Since $\lu$ is a
right integral for $H_1^*$ such that $\lu(1_1) \neq 0$, $H_1$ is cosemisimple (by \cite[Theorem 14.0.3]{sweed}).
Let us show Part (b). Suppose that $H$ is finite dimensional and $H_1$ is cosemisimple. By \cite[Theorem
14.0.3]{sweed}, there exists a right integral $\phi$ for $H_1^*$ such that $\phi(1_1) =1$. By
Theorem~\ref{piintegexists}, there exists a non-zero right \p-integral $\lambda=( \la )_{\al\in \pi}$ for $H$. In
particular, $\lu$ is a non-zero right integral for $H_1^*$. Therefore, since $H_1$ is finite dimensional, there
exists $k \in \Bbbk$ such that $\phi=k \lu$ (by \cite[Theorem 5.1.6]{sweed}). Thus $(k \la )_{\al\in \pi}$ is a
right \p-integral for $H$ such that $k \lu (1_1)=1$. Hence $H$ is cosemisimple by Theorem~\ref{cosemicrit}. This
completes the proof of the corollary.
\end{proof}

\begin{corollary}
Let $H$ be a finite dimensional Hopf \p-coalgebra over a field $\Bbbk$ of characteristic $0$. Then $H$ is
semisimple if and only if it is cosemisimple.
\end{corollary}
\begin{proof}
By Lemma~\ref{semisimplecrit}, $H$ is semisimple if and only if $H_1$ is semisimple, and by
Corollary~\ref{equivcosemifd}(b), $H$ is cosemisimple if and only if $H_1$ is cosemisimple. It is then easy to
conclude using the fact that, in characteristic 0, a finite dimensional Hopf algebra is semisimple if and only if
it is cosemisimple (see \cite[Theorem 3.3]{LR}).
\end{proof}

\begin{corollary}
Let $H$ be a finite dimensional  cosemisimple Hopf \p-coalgebra. If $g=(g_\al)_{\al \in \pi}$ is the distinguished
\p-grouplike element of $H$, then $g=1$ in $G(H)$, i.e., $g_\al=1_\al$ for all $\al \in \pi$. Consequently, the
spaces of left and right \p-integrals for $H$ coincide.
\end{corollary}
\begin{proof}
Let $\al \in \pi$. If $H_\al=0$, then $g_\al=0=1_\al$. Suppose that $H_\al \neq 0$. By Theorem~\ref{cosemicrit},
there exists a right \p-integral $\lambda=(\lambda_\ga)_{\ga \in \pi}$ for $H$ such that $\lambda_\al(1_\al)=1$
and $\lu(1_1)=1$. Then $g_\al = \la(1_\al) \, g_\al =(\id_{H_\al} \otimes \lambda_1) \cp{\al}{1}(1_\al)= \lu(1_1)
\, 1_\al = 1_\al$. By Theorem~\ref{piintegexists} and Lemma~\ref{distinteg}, the spaces of left and right
\p-integrals for $H$ coincide.
\end{proof}

\section{Quasitriangular Hopf $\pi$-coalgebras}\label{quasitriangularity}
In this section, we recall the definitions of crossed, quasitriangular, and ribbon Hopf \p-coalgebras given by
Turaev in \cite{Tur1}, and we generalize the main properties of quasitriangular Hopf algebras to the setting of
Hopf \p-coalgebras.

\subsection{Crossed Hopf $\pi$-coalgebras}\label{deficro}
Following \cite[\S 11.2]{Tur1}, a Hopf \p-coalgebra $H=(\{H_\al\},\Delta,\epsilon,S)$ is said to be \emph{crossed}
provided it is endowed with a family $\varphi=\{\varphi_\be : H_\al \to H_{\be \al \bei} \}_{\al,\be \in \pi}$ of
$\Bbbk$-linear maps (the \emph{crossing}) such that
\begin{defi}
   \item \label{phialgisom} each $\varphi_\be:H_\al \to H_{\be \al \bei}$ is an algebra
               isomorphism;
   \item \label{phicomult} each $\varphi_\be$ preserves the comultiplication, i.e., for all $\al,\be,\ga \in
               \pi$,
               $$(\varphi_\be \otimes \varphi_\be) \cp{\al}{\ga} =\cp{\be \al
               \bei}{\be \ga \bei} \varphi_\be;$$
   \item \label{phicounit} each $\varphi_\be$ preserves the counit, i.e.,
               $\epsilon \varphi_\be=\epsilon$;
   \item \label{phiaction} $\varphi$ is \emph{multiplicative} in the sense that
              $\varphi_{\be \be'}=\varphi_\be \varphi_{\be'}$ for all $\be, \be' \in \pi$.
\end{defi}

\begin{lemma}\label{resA}
Let $H$ be a crossed Hopf \p-coalgebra with crossing $\varphi$. Then
\begin{enumerate}
\renewcommand{\labelenumi}{{\rm (\alph{enumi})}}
 \item $\varphi_{1 \mid H_\al}= \id_{H_\al}$ for all $\al \in \pi$;
 \item $\varphi_{\be}^{-1} = \varphi_\bei$ for all $\be \in \pi$;
 \item $\varphi$ preserves the antipode, i.e.,
       $\varphi_\be S_\al = S_{\be \al \bei} \varphi_\be$ for all
       $\al,\be \in \pi$;
 \item If $\lambda=(\la)_{\al \in \pi}$ is a left
       (resp.\@ right) \p-integral for $H$ and
       $\be \in \pi$, then $(\lbabi\pb)_{\al \in \pi}$ is also a left (resp.\@ right)
       \p-integral for $H$;
 \item If $g=(g_\al)_{\al \in \pi}$ is a \p-grouplike element of $H$ and $\be \in \pi$, then $(\pb(g_{\bei \al
       \be}))_{\al \in \pi}$ is also a \p-grouplike element of $H$.
\end{enumerate}
\end{lemma}
\begin{proof}
Parts (a), (b), (d) and (e) follow directly from the axioms of a crossing. To show Part (c), let $\al, \be \in
\pi$. Using the axioms, it is easy to verify that $\varphi_\be^{-1} S_{\be \al \bei} \varphi_\be *
\id_{H_\al}=\epsilon \, 1_\ali$ in the convolution algebra $\conv(H,H_\ali)$ (see \S\ref{convo}). Thus, since
$S_\al$ is the inverse of $\id_{H_\ali}$ in $\conv(H,H_\ali)$, we have that $\varphi_\be^{-1} S_{\be \al \bei}
\varphi_\be=S_\al$ and so $S_{\be \al \bei} \varphi_\be = \varphi_\be S_\al$.
\end{proof}

\begin{corollary}\label{phitopi}
Let $H$ be a finite dimensional crossed Hopf \p-coalgebra with crossing $\varphi$. Then there exists a unique
group homomorphism $\widehat{\varphi} : \pi \to \kk^*$ such that if $\lambda=(\la)_{\al \in \pi}$ is a left or
right \p-integral for $H$, then $\lbabi \pb= \widehat{\varphi}(\be) \la$ for all $\al,\be \in \pi$.
\end{corollary}
\begin{proof}
Let $\lambda=(\la)_{\al \in \pi}$ be a non-zero left \p-integral for $H$. For any $\be \in \pi$, since
$(\lambda_{\be \al \bei} \pb)_{\al \in \pi}$ is a non-zero left \p-integral for $H$ (see Lemma~\ref{resA}(d)) and
by the uniqueness (within scalar multiple) of a left \p-integral in the finite dimensional case (see
Theorem~\ref{piintegexists}), there exists a unique $\widehat{\varphi}(\be) \in \kk^*$ such that $\lambda_{\be \al
\bei} \pb= \widehat{\varphi}(\be) \la$ for all $\al \in \pi$. Using \eqref{phiaction} and Lemma~\ref{resA}, one
verifies that $\widehat{\varphi}: \pi \to \kk^*$ is a group homomorphism. Since any left \p-integral for $H$ is a
scalar multiple of $\lambda$, the result holds for any left \p-integral. Finally, let $\lambda=(\la)_{\al \in
\pi}$ be a right \p-integral for $H$. Since the antipode is bijective ($H$ is finite dimensional), and using
Lemma~\ref{resA}(d) and the fact that $(\lambda_\ali S_\al)_{\al \in \pi}$ is a left \p-integral for $H$, we have
that, for all $\al,\be \in \pi$, $
  \lambda_{\be \al \bei} \varphi_\be =
  \lambda_{\be \al \bei} S_{\be \ali \bei} \varphi_\be
  S_\ali^{-1} = \widehat{\varphi}(\be) \la S_\ali S_\ali^{-1} = \widehat{\varphi}(\be) \la
$.
\end{proof}

\begin{lemma}\label{marre}
Let $H$ be a finite dimensional crossed Hopf \p-coalgebra with crossing $\varphi$. Let $\widehat{\varphi}$ be as
in Corollary~\ref{phitopi}. Then, for any $\be \in \pi$,
\begin{enumerate}
\renewcommand{\labelenumi}{{\rm (\alph{enumi})}}
 \item If $\Lambda$ is a left or right integral for $H_1$, then
       $\pb(\Lambda)=\widehat{\varphi}(\be) \Lambda$;
 \item If $\nu$ is the distinguished grouplike element of
       $H_1^*$, then $\nu \pb=\nu$;
 \item If $g=(g_\al)_{\al \in \pi}$ is
       the distinguished \p-grouplike element of $H$, then
       $\pb(g_\al)=g_{\be \al \bei}$ for all $\al\in \pi$.
\end{enumerate}
\end{lemma}
\begin{proof}
Let us show Part (a). Let $\Lambda$ be a left integral for $H_1$. We can assume that $\Lambda\neq 0$ (if
$\Lambda=0$, then the result is obvious). By Lemma~\ref{resA} and \eqref{phicounit}, $ x \, \pb(\Lambda) =
\pb(\pbi(x) \, \Lambda)= \pb( \epsilon \pbi(x) \, \Lambda) = \epsilon(x) \, \pb(\Lambda)$ for any $x \in H_1$.
Thus $\pb(\Lambda)$ is a left integral for $H_1$. Therefore, since $H_1$ is finite dimensional and $\Lambda \neq
0$, there exists $k \in \kk$ such that $\pb(\Lambda)=k \Lambda$. Let $\lambda=(\la)_{\al \in \pi}$ be a non-zero
right \p-integral for $H$. We have that $\widehat{\varphi}(\be) \lu(\Lambda)= \lu(\pb(\Lambda))= \lu(k\Lambda)= k
\lu(\Lambda)$. Now $\lu(\Lambda)\neq 0$ (because $\Lambda$ is a non-zero left integral for $H_1$ and $\lu$ is a
non-zero right integral for $H_1^*$). Hence $k=\widehat{\varphi}(\be)$ and so $\pb(\Lambda)=\widehat{\varphi}(\be)
\Lambda$. It can be shown similarly that the result holds if $\Lambda$ is a right integral for $H_1$.

Let us show Part (b). If $\Lambda$ is a left integral for $H_1$, then, for all $x \in H_1$, $\Lambda \, x =\pbi
(\pb(\Lambda) \, \pb(x))= \pbi (\nu(\pb(x)) \, \pb(\Lambda)) = \nu\pb(x) \, \Lambda$ (since $\pb(\Lambda)$ is a
left integral for $H_1$). Thus, by the uniqueness of the distinguished grouplike element of the Hopf algebra
$H_1^*$, we have that $\nu\pb=\nu$.

To show Part (c), let $\lambda=(\la)_{\al \in \pi}$ be a right \p-integral for $H$. By Lemma~\ref{resA}(d),
$(\lambda_{\bei \al \be}\pbi)_{\al \in \pi}$ is also a right \p-integral for $H$. Then, for any $\al, \ga \in
\pi$, using \eqref{phicomult} and Lemmas~\ref{distinteg} and \ref{resA},
\begin{eqnarray*}
(\id_{H_\al} \otimes \lambda_\ga) \cp{\al}{\ga}
  & = & \pbi (\id_{H_{\be \al \bei}} \otimes \lambda_\ga \pbi) \cp{\be \al
      \bei}{\be \ga \bei} \,\pb \\
  & = & \pbi (\lambda_{\al \ga} \pbi \pb \; g_{\be \al \bei})\\
  & = & \lambda_{\al \ga} \; \pbi(g_{\be\al\bei}).
\end{eqnarray*}
Hence, by the uniqueness of the distinguished \p-grouplike element (see Lemma~\ref{distinteg}), we have that
$\pbi(g_{\be \al\bei})=g_\al$ and so $\pb(g_\al)=g_{\be \al\bei}$ for all $\al \in \pi$.
\end{proof}

\subsubsection{The opposite (resp.\@ coopposite) Hopf $\pi$-coalgebra}\label{opcopphi}
Let $H$ be a crossed Hopf \p-coalgebra with crossing $\varphi$. If the antipode of $H$ is bijective, then the
opposite (resp.\@ coopposite) Hopf \p-coalgebra to $H$ (see \S\ref{op} and \S\ref{coop}) is crossed with crossing
given by $\varphi_{\be|H_\al^\opp}^\opp=\varphi_{\be|H_\al}$ (resp.
$\varphi_{\be|H_\al^\cop}^\cop=\varphi_{\be|H_\ali}$) for all $\al,\be \in \pi$.

\subsubsection{The mirror Hopf $\pi$-coalgebra}\label{mirror}
Let $H=(\{H_\al\},\Delta,\epsilon,S,\varphi)$ be a crossed Hopf \p-coalgebra. Following \cite[\S 11.6]{Tur1}, its
\emph{mirror} $\overline{H}$ is defined by the following procedure: set $\overline{H}_\al=H_\ali$ as an algebra,
$\overline{\Delta}_{\al,\be}= (\pb \otimes \id_{H_\bei})\cp{\bei \ali \be}{\bei}$, $\overline{\epsilon}=\epsilon$,
$\overline{S}_\al=\pa S_\ali$ and $\overline{\varphi}_{\be | \overline{H}_\al}=\varphi_{\be| H_\ali}$. It is also
a crossed Hopf \p-coalgebra.

\subsection{Quasitriangular Hopf $\pi$-coalgebras}\label{quasitrig}
Following \cite[\S 11.3]{Tur1}, a \emph{quasitriangular} Hopf \p-coalgebra is a crossed Hopf \p-coalgebra
$H=(\{H_\al\},\Delta,\epsilon,S,\varphi)$ endowed with  a family $R=\{R_{\al,\be} \in H_\al \otimes
H_\be\}_{\al,\be \in \pi}$ of invertible elements (the \emph{$R$-matrix}) such that
\begin{defi}
  \item \label{quasicop} for any $\al,\be \in \pi$ and $x \in H_{\al\be}$,
             $$ R_{\al,\be} \cdot \cp{\al}{\be}(x)= \sigma_{\be,\al}
             (\varphi_\ali\otimes \id_{H_\al}) \cp{\al \be \ali}{\al}(x) \cdot R_{\al,\be}$$
             where $\sigma_{\be,\al}$ denotes the flip
             $H_\be \otimes H_\al \to H_\al \otimes
              H_\be$;
  \item \label{quasicomult} for any $\al,\be \in \pi$,
  \begin{eqnarray*}
  && (\id_{H_\al} \otimes \cp{\be}{\ga})(R_{\al,\be \ga})=(R_{\al,\ga})_{1 \be 3}
     \cdot (R_{\al,\be})_{12 \ga} \\
  && (\cp{\al}{\be} \otimes \id_{H_\ga})(R_{\al \be, \ga})=[(\id_{H_\al} \otimes
     \varphi_{\bei})(R_{\al, \be \ga \bei})]_{1 \be 3} \cdot (R_{\be,\ga})_{\al 23}
  \end{eqnarray*}
  where, for $\Bbbk$-spaces $P,Q$ and $r=\sum_j p_j
  \otimes q_j \in P \otimes Q$, we
  set $r_{12\ga}=r \otimes 1_\ga \in P \otimes Q \otimes H_\ga$,
  $r_{\al 23}=1_\al \otimes r \in H_\al \otimes P \otimes Q$ and
  $r_{1 \be 3}=\sum_j p_j \otimes 1_\be \otimes q_j
       \in P \otimes H_\be \otimes Q$;
  \item \label{quasiphi} the family $R$ is invariant under the crossing,
  i.e., for any $\al,\be, \ga \in \pi$, $$(\varphi_\be \otimes
  \varphi_\be)(R_{\al,\ga})=R_{\be \al \bei, \be \ga \bei}.$$
\end{defi}

Note that $R_{1,1}$ is a (classical) $R$-matrix for the Hopf algebra $H_1$.

When $\pi$ is abelian and $\varphi$ is trivial, one recovers the definition of a quasitriangular \p-colored Hopf
algebra given by Ohtsuki in \cite{Oh}.

If $\pi$ is finite, then an $R$-matrix for $H$ does not necessarily give rive to a (usual) $R$-matrix for the Hopf
algebra $\tilde{H}=\oplus_{\al \in \pi} H_\al$ (see \S\ref{pifinite}). However, if the group $\pi$ is finite
abelian and if $\varphi$ is trivial, then $\tilde{R}=\sum_{\al,\be \in \pi} \rab $ is an $R$-matrix for
$\tilde{H}$.

\begin{notation}
In the proofs, when we write a component $\rab$ of an $R$-matrix as $\rab=a_\al \otimes b_\be$, it is to signify
that $\rab=\sum_j a_j \otimes b_j$ for some $a_j \in H_\al$ and $a_j \in H_\be$, where $j$ runs over a finite set
of indices.
\end{notation}

We now generalize the main properties of quasitriangular Hopf algebras (see \cite{Drin}) to the setting of
quasitriangular Hopf \p-coalgebras.

\begin{lemma}\label{B1}
Let $H=(\{H_\al\},\Delta,\epsilon,S,\varphi,R)$ be a quasitriangular Hopf \p-coalgebra. Then, for any $\al, \be,
\ga \in \pi$,
\begin{enumerate}
\renewcommand{\labelenumi}{{\rm (\alph{enumi})}}
  \item $(\epsilon \otimes \id_{H_\al} )(R_{1,\al}) = 1_\al =(\id_{H_\al} \otimes
  \epsilon) (R_{\al,1})$;
  \item $(S_\ali \varphi_\al \otimes \id_{H_\be} )(R_{\ali,\be})= R_{\al,\be}^{-1}$ \, and \,
  $(\id_{H_\al} \otimes S_\be)(R_{\al,\be}^{-1})= R_{\al,\bei}$;
  \item $(S_\al \otimes S_\be)(R_{\al,\be})= (\varphi_\al \otimes
  \id_{H_{\al^{-1}}})(R_{\al^{-1},\be^{-1}})$;
  \item $(R_{\be,\ga})_{\al 23} \cdot (R_{\al,\ga})_{1 \be 3} \cdot (R_{\al,\be})_{12 \ga}
  $\\
  \phantom{xxxxxxx} = $(R_{\al,\be})_{1 2 \ga} \cdot [(\id_{H_\al} \otimes \varphi_\bei)
  (R_{\al,\be \ga \bei})]_{1 \be 3}
  \cdot (R_{\be,\ga})_{\al 23} $.
\end{enumerate}
\end{lemma}
Part (d) of Lemma~\ref{B1}, which is the Yang-Baxter equality for $R=\{R_{\al,\be}\}_{\al,\be \in \pi}$, first
appeared in \cite[\S 11.3]{Tur1}. We prove it here for completeness sake.
\begin{proof}
Let us show Part (a). We have
\begin{eqnarray*}
R_{1,\al}  
  & = & (\epsilon \otimes \id_{H_1} \otimes \id_{H_\al})
        (\cp{1}{1} \otimes \id_{H_\al})(R_{1,\al}) \text{ \quad by \eqref{counit}}\\
  & = & (\epsilon \otimes \id_{H_1} \otimes \id_{H_\al})([(\id_{H_1} \otimes
        \varphi_1)(R_{1,\al})]_{1 {\scriptscriptstyle 1_\pi} 3}
        \cdot (R_{1,\al})_{{\scriptscriptstyle 1_\pi}23})
        \quad \text{ by \eqref{quasicomult}} \\
  & = & (\epsilon \otimes \id_{H_1} \otimes \id_{H_\al})((R_{1,\al})_{1{\scriptscriptstyle 1_\pi}3}
        \cdot (R_{1,\al})_{{\scriptscriptstyle 1_\pi}23})
        \quad \text{ by Lemma~\ref{resA}(a)} \\
  & = & (\epsilon \otimes \id_{H_1} \otimes \id_{H_\al})((R_{1,\al})_{1{\scriptscriptstyle 1_\pi}3})
        \cdot (\epsilon
        \otimes \id_{H_1} \otimes \id_{H_\al}) ((R_{1,\al})_{{\scriptscriptstyle 1_\pi}23})
        \text{ \; by \eqref{hopfmor}}\\
  & = & (1_1 \otimes (\epsilon \otimes \id_{H_\al})(R_{1,\al})) \cdot
        R_{1,\al}.
\end{eqnarray*}
Thus $1_1 \otimes(\epsilon \otimes \id_{H_\al})(R_{1,\al})=1_1 \otimes 1_\al$ (since $R_{1,\al}$ is invertible).
By applying $(\epsilon \otimes \id_{H_\al})$ on both sides, we get the first equality of Part (a). The second one
can be obtained similarly.

To show the first equality of Part (b), set $$\mathcal{E} = (m_\al \otimes \id_{H_\be})(S_\ali \otimes \id_{H_\al}
\otimes \id_{H_\be})(\cp{\ali}{\al} \otimes\id_{H_\be})(R_{1,\be}).$$ Let us compute $\mathcal{E}$ in two
different ways. On the first hand,
\begin{eqnarray*}
\mathcal{E}
  & = & (m_\al \otimes \id_{H_\be})(S_\ali \otimes \id_{H_\al} \otimes \id_{H_\be})\\
  &   & \phantom{XXX}
        ([(\id_{H_\ali} \otimes \varphi_\ali) (R_{\ali,\al \be \ali})]_{1 \al 3}
        \cdot (R_{\al,\be})_{\ali 23})  \quad \text{ by \eqref{quasicomult}} \\
  & = & (S_\ali \otimes \varphi_\ali)(R_{\ali, \al \be \ali}) \cdot R_{\al,\be}\\
  & = & (S_\ali \pa \otimes \id_{H_\be})(R_{\ali, \be}) \cdot R_{\al,\be}
        \quad \text{ by \eqref{quasiphi}}.
\end{eqnarray*}
On the second hand,
\begin{eqnarray*}
\mathcal{E} 
              & = & (\epsilon 1_\al \otimes \id_{H_\be})(R_{1,\be})
                    \text{ \quad by \eqref{antipode}}\\
              & = & 1_\al \otimes 1_\be \text{ \quad by
                    Part (a).}
\end{eqnarray*}
Comparing these two calculations and since $\rab$ is invertible, we get the first equality of Part (b). The second
one can be proved similarly by computing the expression $\mathcal{F}= (\id_{H_\al} \otimes m_\bei)(\id_{H_\al}
\otimes \id_{H_\bei} \otimes S_\be) (\id_{H_\al} \otimes \cp{\bei}{\be})(R_{\al,1}^{-1})$.

Part (c) is a direct consequence of Part (b) and Lemma~\ref{resA}(a) and (c).

Finally, Part (d) follows from axioms \eqref{quasicop} and \eqref{quasicomult}:
\begin{eqnarray*}
  \lefteqn{
  (R_{\be,\ga})_{\al 23} \cdot (R_{\al,\ga})_{1 \be 3} \cdot(R_{\al,\be})_{12\ga}
  }\\
  & = & (R_{\be,\ga})_{\al 23} \cdot (\id_{H_\al} \otimes
        \cp{\be}{\ga})(R_{\al,\be\ga}) \\
  & = & (\id_{H_\al} \otimes R_{\be,\ga} \cdot \cp{\be}{\ga})(R_{\al,\be\ga}) \\
  & = & (\id_{H_\al} \otimes \sigma_{\ga,\be}(\varphi_\bei \otimes \id_{H_\be}) \cp{\be \ga \bei}{\be}
        \cdot R_{\be,\ga})(R_{\al,\be\ga}) \\
  & = & (\id_{H_\al} \otimes \sigma_{\ga,\be}(\varphi_\bei \otimes \id_{H_\be}))
        ((R_{\al,\be})_{1 \be \ga \bei 3} \cdot (R_{\al, \be \ga
        \bei})_{12 \be}) \cdot (R_{\be,\ga})_{\al 23} \\
  & = & (R_{\al,\be})_{12 \ga} \cdot [(\id_{H_\al} \otimes \varphi_\bei)(R_{\al, \be \ga
        \bei})]_{1 \be 3} \cdot (R_{\be,\ga})_{\al 23}.
\end{eqnarray*}
This completes the proof of the lemma.
\end{proof}

\subsection{The Drinfeld elements}  Let $H=(\{H_\al,m_\al,1_\al\},\Delta,\epsilon,S,\varphi,R)$
be a quasitriangular Hopf \p-coalgebra. We define the {\it (generalized) Drinfeld elements} of $H$, for any $ \al
\in \pi$, by
  $$
  u_\al = m_\al (S_\ali \varphi_\al \otimes \id_{H_\al} ) \sigma_{\al, \ali} (R_{\al,\ali}) \in H_\al.
  $$
Note that $u_1$ is the Drinfeld element of the quasitriangular Hopf algebra $H_1$ (see \cite{Drin}).
\begin{lemma}\label{B2}
For any $\al,\be \in \pi$,
\begin{enumerate}
\renewcommand{\labelenumi}{{\rm (\alph{enumi})}}
  \item $u_\al$ is invertible and $u_\al^{-1}=m_\al (\id_{H_\al} \otimes S_\ali S_\al) \sigma_{\al,\al}
        (R_{\al,\al})$;
  \item $S_\ali S_\al (\varphi_\al (x))= u_\al x u_\al^{-1}$ for all
        $x \in H_\al$;
  \item The antipode of $H$ is bijective;
  \item $\varphi_\be (u_\al ) = u_{\be \al \bei}$;
  \item $S_\ali (u_\ali) u_\al = u_\al S_\ali (u_\ali)$ and
        this element, noted $c_\al$, 
        verifies\\ $c_\al \pai(x)=\pa(x) c_\al$ for all $x \in H_\al$;
  \item $\cp{\al}{\be}(u_{\al \be}) = [ \sigma_{\be,\al}(\id_{H_\be}
        \otimes \pa )(R_{\be,\al}) \cdot R_{\al,\be} ]^{-1} \cdot
        ( u_\al \otimes u_\be)$\\
        $\phantom{\cp{\al}{\be}(u_{\al \be})}= (u_\al \otimes u_\be) \cdot [\sigma_{\be,\al}
        (\varphi_\bei \otimes \id_{H_\al})(R_{\be,\al})
        \cdot (\varphi_\ali \otimes \varphi_\bei)(R_{\al,\be})]^{-1}$;
  \item $\epsilon (u_1)=1$.
\end{enumerate}
\end{lemma}

\begin{proof}
We adapt the methods used in \cite{Drin} to our setting. Let us prove Parts~(a) and (b). We first show that for
all $x \in H_\al$,
\begin{equation}\label{pipi}
 S_\ali S_\al (\varphi_\al (x)) u_\al= u_\al x.
\end{equation}
Write $R_{\al,\ali}=a_\al \otimes b_\ali$ so that $u_\al = S_\ali(\varphi_\al (b_\ali)) a_\al$. Let $x \in H_\al$.
Using \eqref{coass} and \eqref{quasicop}, we have that
\begin{multline*}
  (R_{\al,\ali})_{12\al} \cdot(\id_{H_\al} \otimes \cp{\ali}{\al}) \cp{\al}{1}(x) = \\
  (\sigma_{\ali,\al}(\varphi_{\ali} \otimes \id_{H_\al}) \cp{\ali}{\al} \otimes \id_{H_\al}) \cp{1}{\al}(x) \cdot
  (R_{\al, \ali})_{12\al}\, ,
\end{multline*}
that is $a_\al \xua \otimes b_\ali \xdai \otimes \xta = \xda a_\al \otimes \varphi_\ali (\xuai) b_\ali \otimes
\xta$. Evaluate both sides of this equality with $(\id_{H_\al} \otimes S_\ali \varphi_\al \otimes S_\ali S_\al
\varphi_\al)$, reverse the order of the tensorands and multiply them to obtain
\begin{multline*}
S_\ali S_\al \varphi_\al (\xta)S_\ali\varphi_\al (b_\ali \xdai) a_\al \xua \\
=S_\ali S_\al \varphi_\al (\xta) S_\ali \varphi_\al (\varphi_\ali(\xuai)b_\ali) \xda a_\al.
\end{multline*}
Now, by Lemmas \ref{antipodepptes}(a) and \ref{resA}(c), the left-hand side is equal to
\begin{eqnarray*}
  \lefteqn{S_\ali \varphi_\al S_\al (\xta) S_\ali\varphi_\al(\xdai)
            S_\ali (\varphi_\al (b_\ali)) a_\al \xua}\\
      & = & S_\ali \varphi_\al(\xdai S_\al(\xta)) u_\al \xua \\
      & = & S_\ali \varphi_\al(\epsilon (\xdu)1_\ali) u_\al \xua
            \text{ \quad by \eqref{antipode}} \\
      & = & u_\al \epsilon(\xdu) \xua
            \text{ \quad since $S_\ali \varphi_\al(1_\ali)=1_\al$}\\
      & = & u_\al x \text{ \quad by \eqref{counit},}
\end{eqnarray*}
and, by Lemma \ref{antipodepptes}(a), the right-hand side is equal to
\begin{eqnarray*}
  \lefteqn{S_\ali S_\al \varphi_\al (\xta)
             S_\ali(\varphi_\al(b_\ali)) S_\ali(\xuai) \xda a_\al}\\
       & = & S_\ali S_\al \varphi_\al (\epsilon(\xuu) \xda) S_\ali
             (\varphi_\al(b_\ali)) a_\al
             \text{ \quad by \eqref{antipode}}\\
       & = & S_\ali S_\al \varphi_\al(x) u_\al \text{ \quad by \eqref{counit}.}
\end{eqnarray*}
Thus \eqref{pipi} is proven. Let us show that $u_\al$ is invertible. Set
$$\widetilde{u}_\al=m_\al(\id_{H_\al} \otimes S_\ali S_\al) \sigma_{\al,\al} (R_{\al,\al}) \in H_\al.$$
By Lemma~\ref{B1}(b) and \eqref{quasiphi}, $R_{\al,\al}=(\id_{H_\al} \otimes S_\ali)(\pa \otimes
\pa)(R_{\al,\ali}^{-1})$. Write $R_{\al,\ali}^{-1}=c_\al \otimes d_\ali $. Then $\widetilde{u}_\al=S_\ali (\pa
(d_\ali)) S_\ali S_\al (\pa (c_\al))$ and $a_\al c_\al \otimes b_\ali d_\ali=1_\al \otimes 1_\ali$. Now
\begin{eqnarray*}
 \widetilde{u}_\al u_\al
   & = & S_\ali (\pa (d_\ali)) S_\ali S_\al (\pa (c_\al)) u_\al \\
   & = & S_\ali (\pa (d_\ali)) u_\al c_\al
         \text{ \quad by (\ref{pipi}) with $x=c_\al$} \\
   & = & S_\ali (\pa (d_\ali)) S_\ali(\varphi_\al(b_\ali))a_\al c_\al \\
   & = & S_\ali \pa (b_\ali d_\ali)) a_\al c_\al
         \text{ \quad by Lemma~\ref{antipodepptes}(a)}\\
   & = & S_\ali \pa (1_\ali)) 1_\al \; = \; 1_\al.
\end{eqnarray*}
It can be shown similarly that $u_\al  \widetilde{u}_\al =1_\al$. Thus $u_\al$ is invertible,
$u_\al^{-1}=\widetilde{u}_\al$, and so $S_\ali S_\al (\varphi_\al (x))= u_\al x u_\al^{-1}$ for any $x \in H_\al$.

Part (c) is a direct consequence of Part (b). Part (d) follows from \eqref{phialgisom}, \eqref{phiaction}, and
\eqref{quasiphi}.
Let us show Part (e). For any $x \in H_\al$,
\begin{eqnarray*}
 \lefteqn{S_\ali(u_\ali) u_\al \pai(x)}\\
    & = & S_\ali(u_\ali) S_\ali S_\al (x) u_\al
          \text{ \quad by Part (b)}\\
    & = & S_\ali(u_\ali) S_\ali S_\al S_\ali (\pai S_\ali^{-1}(\pa(x))) u_\al
          \text{ \quad by Lemma~\ref{resA}(c)} \\
    & = & S_\ali(u_\ali) S_\ali(u_\ali S_\ali^{-1}(\pa(x))
          u_\ali^{-1}) u_\al \text{ \quad by Part (b)}\\
    & = & \pa(x) S_\ali(u_\ali) u_\al \text{ \quad since $S_\ali$ is anti-multiplicative.}
\end{eqnarray*}
In particular, for $x=u_\al$, one gets that $S_\ali(u_\ali) u_\al = u_\al S_\ali(u_\ali)$.

For the proof of the first equality of Part (f), set $\rt_{\al,\be}=\sigma_{\be,\al}(\id_{H_\be} \otimes
\pa)(\rba)$. Using Lemma~\ref{resA} and \eqref{quasiphi}, we have also that $\rt_{\al,\be}=\sigma_{\be,\al}(\pai
\otimes \id_{H_\al})(R_{\al \be \ali,\al})$. We first show that for all $x \in H_{\al \be}$,
\begin{equation}\label{rt}
  \rtab \cdot \rab \cdot \cp{\al}{\be}(x)= (\pa \otimes \pb)
  \cp{\al}{\be}(\pabi(x)) \cdot \rtab \cdot \rab.
\end{equation}
By \eqref{quasicop}, $\rba \cdot \cp{\be}{\al}(\pai(x)) = \sigma_{\al,\be} (( \pbi \otimes \id_{H_\be})\cp{\be \al
\bei}{\be}(\pai(x))) \cdot \rba$. Evaluate both sides of this equality with the algebra homomorphism
$\sigma_{\be,\al}(\id_{H_\be} \otimes \pa)$ and multiply them on the right by $\rab$ to obtain
\begin{eqnarray*}
  \lefteqn{\sigma_{\be,\al} (\id_{H_\be} \otimes \pa)(\rba) \cdot
  \sigma_{\be,\al}(\id_{H_\be} \otimes \pa) \cp{\be}{\al}(\pai(x)) \cdot \rab}\\
  & = & (\pa \pbi \otimes \id_{H_\be}) \cp{\be \al \bei}{\be}(\pai(x)) \cdot \sigma_{\be,\al}
        (\id_{H_\be} \otimes \pa)(\rba) \cdot \rab.
\end{eqnarray*}
Then, using \eqref{quasicop} and \eqref{phicomult}, one gets equality \eqref{rt}. Set now
$$
  \mathcal{E}=
  \rtab \cdot \rab \cdot \cp{\al}{\be}(u_{\al \be}).
$$
We have to show that $\mathcal{E} = u_\al \otimes u_\be$. Write $R_{\al \be, (\al \be)^{-1}}=r \otimes s$, $\rab =
a_\al \otimes b_\be$, and $\rtab = c_\al \otimes d_\be$. Then $ u_{\al \be}=S_{(\al \be)^{-1}} (\pab(s)) r = \pab
S_{(\al \be)^{-1}}(s) r$. We have that
\begin{eqnarray*}
  \mathcal{E}
    & = & \rtab \cdot \rab \cdot \cp{\al}{\be}(\pab S_{(\al \be)^{-1}}(s) r)\\
    & = & \rtab \cdot \rab \cdot \cp{\al}{\be}(\pab S_{(\al \be)^{-1}}(s))
          \cdot \cp{\al}{\be}(r) \text{ \quad by \eqref{hopfmor}.}
\end{eqnarray*}
Therefore, using \eqref{rt} for $x=\pab S_{(\al \be)^{-1}}(s)$ and then Lemmas~\ref{antipodepptes}(c) and
\ref{resA}(c),
\begin{eqnarray*}
  \mathcal{E}
    & = & (\pa \otimes \pb) \cdot \cp{\al}{\be}(S_{(\al \be)^{-1}}(s))
          \cdot \rtab \cdot \rab \cdot \cp{\al}{\be}(r) \\
    & = & (\pa \otimes \pb)\tba (S_\bei \otimes S_\ali)
          \cp{\bei}{\ali}(s) \cdot \rtab \cdot \rab \cdot
          \cp{\al}{\be}(r) \\
    & = & \pa S_\ali(s_{(2,\ali)}) c_\al a_\al r_{(1,\al)} \otimes \pb S_\bei
          (s_{(1,\bei)}) d_\be b_\be r_{(2,\be)}\\
    & = &  S_\ali\pa(s_{(2,\ali)}) c_\al a_\al r_{(1,\al)} \otimes S_\bei\pb
          (s_{(1,\bei)}) d_\be b_\be r_{(2,\be)}.
\end{eqnarray*}
Now $H_\al \otimes H_\be$ is a right $H_\al \otimes H_\be \otimes H_\ali \otimes H_\bei$-module under the action
$$(x \otimes y) \twoheadleftarrow (h_1 \otimes h_2 \otimes h_3 \otimes h_4) = S_\ali(\pa(h_3)) x h_1 \otimes S_\bei
(\pb(h_4)) y h_2.$$  Then
\begin{equation*}
\begin{split}
  \mathcal{E} \;\,
    & = \;\;\, c_\al \otimes d_\be \twoheadleftarrow a_\al r_{(1,\al)} \otimes
          b_\be r_{(2,\ali)} \otimes s_{(2,\ali)} \otimes
          s_{(1,\bei)} \\
    & = \;\;\, \rtab \twoheadleftarrow (\rab)_{12 \ali \bei} \cdot
          (\cp{\al}{\be} \otimes \sigma_{\bei,\ali}
          \cp{\bei}{\ali})(R_{\al \be,(\al \be)^{-1}})\\
    & = \;\;\, \rtab \twoheadleftarrow (\rab)_{12 \ali \bei} \cdot
          (\cp{\al}{\be} \otimes \id_{H_\ali} \otimes \id_{H_\bei})\\
    & = \phantom{XXXXXXX}
           ((R_{\al \be,\ali})_{12 \bei}
          \cdot (R_{\al\be,\bei})_{1 \ali 3}) \text{ \quad by \eqref{quasicomult}}\\
    & = \;\;\, \rtab \twoheadleftarrow (\rab)_{12 \ali \bei} \cdot
          (\cp{\al}{\be} \otimes \id_{H_\ali} \otimes \id_{H_\bei})
          ((R_{\al \be,\ali})_{12 \bei}) \\
    &  \phantom{XXXXXXXx} \cdot (\cp{\al}{\be} \otimes \id_{H_\ali} \otimes \id_{H_\bei})
           ((R_{\al\be,\bei})_{1 \ali 3}) \text{ \quad by \eqref{hopfmor}.}
\end{split}
\end{equation*}
For any $\Bbbk$-spaces $P,Q$ and any $x=\sum_j p_j \otimes q_j \in P \otimes Q$, we set $x_{12\al \be}=x \otimes
1_\al \otimes 1_\be \in P \otimes Q \otimes H_\al \otimes H_\be$, $x_{\al 2 \be 4}=\sum_j 1_\al \otimes p_j
\otimes 1_\be \otimes q_j \in H_\al \otimes P \otimes H_\be \otimes Q$, etc. Therefore, by \eqref{quasicomult} and
Lemma~\ref{B1}(d),
\begin{equation*}
\begin{split}
  \mathcal{E} \;\,
    & = \;\;\, \rtab \twoheadleftarrow (\rab)_{12 \ali \bei} \cdot [(\id_{H_\al} \otimes
          \pbi)(R_{\al,\be\ali\bei})]_{1 \be 3 \bei} \\
    &  \phantom{XXXXXx} \cdot (R_{\be,\ali})_{\al 23\bei}
       \cdot [(\id_{H_\al} \otimes \pbi)(R_{\al,\bei})]_{1 \be \ali 4} \cdot (R_{\be,\bei})_{\al 2
          \ali 4} \\
    & = \;\;\, \rtab \twoheadleftarrow (R_{\be, \ali})_{\al 23 \bei} \cdot
          (R_{\al,\ali})_{1 \be 3 \bei} \cdot (R_{\al,\be})_{12 \ali \bei} \\
    & \qquad \qquad \quad \;\, \cdot
          [(\id_{H_\al} \otimes \pbi)(R_{\al,\bei})]_{1 \be \ali 4} \cdot (R_{\be,\bei})_{\al 2 \ali 4}.
\end{split}
\end{equation*}
Write $\rba=e_\be \otimes f_\al$ and $\rbai=g_\be \otimes h_\ali$. Then $ \rtab = \pa(f_\al) \otimes e_\be$ and so
\begin{eqnarray*}
 \lefteqn{ \rtab \twoheadleftarrow (\rbai)_{\al 23 \bei}} \\
     & = & S_\ali(\pa (h_\ali)) \pa(f_\al) \otimes e_\be g_\be \\
     & = & \tba (\id_{H_\be} \otimes \pa S_\ali) (
          (\id_{H_\be} \otimes S_\ali^{-1})(\rba) \cdot \rbai) \text{ \quad by Lemma~\ref{resA}(c)}\\
     & = & \tba (\id_{H_\be} \otimes \pa S_\ali) (
           \rbai^{-1} \cdot \rbai)
           \text{ \quad by Lemma~\ref{B1}(b)}\\
     & = & 1_\al \otimes 1_\be.
\end{eqnarray*}
If we write $R_{\al,\ali}=m_\al \otimes n_\ali$, then $$1_\al \otimes 1_\be \twoheadleftarrow (\raai)_{1 \ali 3
\bei} =S_\ali\pa(n_\ali)m_\al \otimes 1_\be = u_\al \otimes 1_\be.$$ Therefore $$\mathcal{E} = u_\al \otimes 1_\be
\twoheadleftarrow (R_{\al,\be})_{12 \ali \bei} \cdot [(\id_{H_\al} \otimes \pbi)(R_{\al,\bei})]_{1 \be \ali 4}
\cdot (R_{\be,\bei})_{\al 2 \ali 4}.$$ Write now $R_{\al,\bei}= p_\al \otimes q_\bei$. Then
\begin{eqnarray*}
 \lefteqn{u_\al \otimes 1_\be \twoheadleftarrow
(R_{\al,\be})_{12 \ali \bei} \cdot [(\id_{H_\al} \otimes
\pbi)(R_{\al,\bei})]_{1 \be \ali 4}} \\
   & = & u_\al a_\al p_\al \otimes S_\bei(q_\bei)b_\be\\
   & = & (u_\al \otimes 1_\be) \cdot (\id_{H_\al} \otimes S_\bei) ((\id_{H_\al} \otimes S_\bei^{-1})(\rab)
         \cdot R_{\al,\bei}) \\
   & = & (u_\al \otimes 1_\be) \cdot (\id_{H_\al} \otimes S_\bei)(R_{\al,\bei}^{-1} \cdot R_{\al,\bei})
         \text{ \quad by Lemma~\ref{B1}(b)}\\
   & = & u_\al \otimes 1_\be.
\end{eqnarray*}
Hence $\mathcal{E} = u_\al \otimes 1_\be \twoheadleftarrow (R_{\be,\bei})_{\al 2 \ali 4}$. Finally, write
$R_{\be,\bei} = x_\be \otimes y_\bei$. Then $\mathcal{E} = u_\al \otimes S_\bei(\pb(y_\bei))x_\be =u_\al \otimes
u_\be$. This completes the proof of the first equality of Part (f). Let us show the second one. Using the first
equality of Part~(f) and then Part~(b), we have
\begin{eqnarray*}
  \cp{\al}{\be}(u_{\al\be})
    & = & [\sigma_{\be,\al} (\id_{H_\be} \otimes \pa)(R_{\be,\al}) \cdot R_{\al,\be}]^{-1} \cdot (u_\al \otimes
           u_\be)\\
    & = & (u_\al \otimes u_\be) \cdot (\pai (S_\ali S_\al)^{-1} \otimes \pbi (S_\bei S_\be)^{-1})\\
    &   & \phantom{ffffff} ([\sigma_{\be,\al} (\id_{H_\be} \otimes \pa)(R_{\be,\al}) \cdot R_{\al,\be}]^{-1}),
\end{eqnarray*}
and so, by Lemmas~\ref{resA} and \ref{B1}(c),
$$
  \cp{\al}{\be}(u_{\al\be})
    =(u_\al \otimes u_\be) \cdot [\sigma_{\be,\al}(\varphi_\bei \otimes \id_{H_\al})(R_{\be,\al})
          \cdot (\varphi_\ali \otimes \varphi_\bei)(R_{\al,\be})]^{-1}.
$$

It remains to show Part (g). We have
\begin{eqnarray*}
  u_1 & = & (\epsilon \otimes \id_{H_1}) \cp{1}{1}(u_1) \text{ \quad by \eqref{counit}}\\
      & = & (\epsilon \otimes \id_{H_1}) ((\sigma_{1,1}(R_{1,1}) \cdot
            R_{1,1})^{-1} \cdot (u_1 \otimes u_1)) \text{ \quad by Part (f)}\\
      & = & (\epsilon \otimes \id_{H_1})(R_{1,1})^{-1}
            \cdot (\id_{H_1} \otimes \epsilon)(R_{1,1})^{-1}
            \cdot \epsilon(u_1) u_1 \text{ \quad by \eqref{hopfmor}}\\
      & = & \epsilon(u_1) u_1 \text{ \quad by Lemma~\ref{B1}(a).}
\end{eqnarray*}
Now $u_1 \neq 0$ since $u_1$ is invertible (by Part (a)) and $H_1 \neq0$ (by Corollary~\ref{subgroup}). Hence
$\epsilon(u_1)=1$. This finishes the proof of the lemma.
\end{proof}

\subsubsection{The coopposite Hopf $\pi$-coalgebra}\label{coopquasi}
Let $H$ be a quasitriangular Hopf \p-coalgebra with $R$-matrix $R=\{R_{\al,\be}\}_{\al,\be \in \pi}$. By
Lemma~\ref{B2}(c), the antipode of $H$ is bijective. Thus we can consider the coopposite crossed Hopf \p-coalgebra
$H^\cop$ to $H$ (see \S\ref{opcopphi}). It is quasitriangular by setting
$$R^\cop_{\al,\be}= (\pa \otimes \id_{H_\bei})(R^{-1}_{\ali,\bei})=(S_\al \otimes \id_{H_\bei})(R_{\al,\bei}).$$
The Drinfeld elements of $H$ and $H^\cop$ are related by $u_\al^\cop=u_\ali^{-1}$.

\subsubsection{The mirror Hopf $\pi$-coalgebra}\label{mirrorquasi}
Let $H$ be a quasitriangular Hopf \p-coalgebra with $R$-matrix $R=\{R_{\al,\be}\}_{\al,\be \in \pi}$. Following
\cite[\S 11.6]{Tur1}, the mirror crossed Hopf \p-coalgebra $\overline{H}$ to $H$ (see \S\ref{mirror}) is
quasitriangular with $R$-matrix given by
$$\overline{R}_{\al,\be}=\sigma_{\bei,\ali}(R_{\bei,\ali}^{-1}).$$ The Drinfeld elements associated to $H$ and
$\overline{H}$ verify $\overline{u}_\al=S_\al(u_\al)^{-1}$.

\begin{corollary}\label{ell}
Let $H$ be a quasitriangular Hopf \p-coalgebra. For all $\al \in \pi$, set $\ell_\al=S_\ali(u_\ali)^{-1}
u_\al=u_\al S_\ali(u_\ali)^{-1} \in H_\al$.  Then
\begin{enumerate}
\renewcommand{\labelenumi}{{\rm (\alph{enumi})}}
  \item $\ell=(\ell_\al)_{\al \in
        \pi}$ is a \p-grouplike element of $H$;
  \item $(S_\ali S_\al)^2(x)=\ell_\al x \ell_\al^{-1}$ for all
        $\al \in \pi$ and $x \in H_\al$.
\end{enumerate}
\end{corollary}
\begin{proof}
Let us show Part (a).  Denote by $\overline{u}_\al$ the Drinfeld elements of the mirror Hopf \p-coalgebra
$\overline{H}$ to $H$ (see \S\ref{mirrorquasi}). Since $\overline{u}_\al=S_\al(u_\al)^{-1}$, Lemma~\ref{B2}(f)
applied to $\overline{H}$ gives that, for any $\al,\be \in \pi$,
\begin{multline*}
  \cp{\al}{\be}(S_{(\al \be)^{-1}} (u_{(\al \be)^{-1}})^{-1}) =
  \sigma_{\be,\al}(\varphi_\bei \otimes \id_{H_\al})(R_{\be,\al}) \\ \cdot
  (\varphi_\ali \otimes \varphi_\bei)(R_{\al,\be}) \cdot(S_\ali(u_\ali)^{-1} \otimes S_\bei(u_\bei)^{-1}).
\end{multline*}
Now, by Lemma~\ref{B2}(f),
$$
  \cp{\al}{\be}(u_{\al\be})
    = (u_\al \otimes u_\be) \cdot [\sigma_{\be,\al}(\varphi_\bei \otimes \id_{H_\al})(R_{\be,\al})
          \cdot (\varphi_\ali \otimes \varphi_\bei)(R_{\al,\be})]^{-1}.
$$
Thus we obtain that $\cp{\al}{\be}(\ell_{\al\be})=\cp{\al}{\be}(u_{\al\be}) \cdot \cp{\al}{\be}(S_{(\al \be)^{-1}}
(u_{(\al \be)^{-1}})^{-1})=\ell_\al \otimes \ell_\be$. Moreover $\epsilon(\ell_1)=\epsilon(u_1
S_1(u_1)^{-1})=\epsilon(u_1) \epsilon (S_1(u_1))^{-1}=\epsilon(u_1) \epsilon(u_1)^{-1}=1$ by \eqref{hopfmor} and
Lemma~\ref{antipodepptes}(d). Hence $\ell=(\ell_\al)_{\al \in \pi} \in G(H)$.

To show Part (b), let $\al \in \pi$ and $x \in H_\al$. Applying Lemma~\ref{B2}(b) to $\overline{H}$ and then to
$H$ gives that
\begin{eqnarray*}
(S_\ali S_\al)^2(x)
  & = & S_\ali S_\al(S_\ali(u_\ali)^{-1}\pa(x) S_\ali(u_\ali))\\
  & = & S_\ali S_\al(S_\ali(u_\ali)^{-1}) \; S_\ali S_\al( \pa(x))
        \; S_\ali S_\al( S_\ali(u_\ali)) \\
  & = & u_\al S_\ali(u_\ali)^{-1} x  S_\ali(u_\ali) u_\al^{-1}\\
  & = & \ell_\al x \ell_\al^{-1}.
\end{eqnarray*}
This completes the proof of the corollary.
\end{proof}

\subsubsection{The double of a crossed Hopf $\pi$-coalgebra}
The Drinfeld double construction for Hopf algebras can be extended to the setting of crossed Hopf \p-coalgebras,
see \cite{Zu}. This yields examples of quasitriangular Hopf \p-coalgebras.

\subsection{Ribbon Hopf $\pi$-coalgebras}\label{ribpico}
Following \cite[\S 11.4]{Tur1}, a quasitriangular Hopf \p-coalgebra $H=(\{H_\al\},\Delta,\epsilon,S,\varphi,R)$ is
said to be \emph{ribbon} if it is endowed with a family $\theta=\{\ta \in H_\al \}_{\al \in \pi}$ of invertible
elements (the {\it twist}) such that
\begin{defi}
  \item \label{twist1} $\varphi_\al(x)=\ta^{-1} x \ta$ for all $\al \in
              \pi$ and $ x \in H_\al$;
  \item \label{twist2} $S_\al(\ta)=\tai$ for all $\al \in \pi$;
  \item \label{twist3} $\pb(\ta)=\theta_{\be \al \bei}$ for all $\al, \be \in \pi$;
  \item \label{twist4} for all $\al, \be \in \pi$, $$\cp{\al}{\be}(\theta_{\al \be})=(\ta \otimes \tb) \cdot
            \sigma_{\be,\al}( (\varphi_\ali \otimes \id_{H_\al} )(R_{\al \be \ali ,\al}))
            \cdot R_{\al,\be}.$$
\end{defi}

Note that $\theta_1$ is a (classical) twist of the quasitriangular Hopf algebra $H_1$.

\begin{lemma}\label{C1}
Let $H=(\{H_\al\},\Delta,\epsilon,S,\varphi,R,\theta)$ be a ribbon Hopf \p-coalgebra. Then
\begin{enumerate}
\renewcommand{\labelenumi}{{\rm (\alph{enumi})}}
  \item $\varphi_\ali(x)= \ta x
        \ta^{-1} $ for all $\al \in \pi$ and $x \in H_\al$;
  \item $\epsilon (\theta_1) = 1 $;
  \item If $\al \in \pi$ has a finite order $d$, then $\theta^d_\al$ is a central
        element of $H_\al$. In particular $\theta_1$ is central;
  \item $\ta u_\al = u_\al \ta $  for all $\al \in \pi$, where the $u_\al$
        are the Drinfeld elements of $H$.
\end{enumerate}
\end{lemma}
\begin{proof}
Part (a) is a direct consequence of \eqref{twist1}, \eqref{twist3}, and Lemma~\ref{resA}. Let us show Part (b). We
have
\begin{eqnarray*}
  \theta_1 & = & (\epsilon \otimes \id_{H_1}) \cp{1}{1}(\theta_1)
                 \text{ \quad by \eqref{counit}}\\
    & = & (\epsilon \otimes \id_{H_1}) ((\theta_1 \otimes \theta_1)
          \cdot \tuu (\ruu) \cdot \ruu)
          \text{ \quad by \eqref{twist4} and Lemma~\ref{resA}(a)}\\
    & = & (\epsilon \otimes \id_{H_1})(\theta_1 \otimes \theta_1)
          \cdot (\id_{H_1} \otimes \epsilon)(\ruu)
          \cdot (\epsilon \otimes \id_{H_1})(\ruu) \text{ \quad by \eqref{hopfmor}}\\
    & = & \epsilon(\theta_1) \theta_1 \text{ \quad by Lemma~\ref{B1}(a).}
\end{eqnarray*}
Now $\theta_1\neq 0$ since it is invertible and $H_1 \neq 0$ (by Corollary~\ref{subgroup}). Hence
$\epsilon(\theta_1)=1$. To show Part (c), let $\al \in \pi$ of finite order~$d$. For any $x \in H_\al$, using
\eqref{phiaction}, Lemma~\ref{resA} and \eqref{twist1}, we have that $x=\varphi_1(x)=\varphi_{\al^d}(x) =
\varphi_\al^d(x)=\theta_\al^{-d} x \theta_\al^d$ and so $\theta_\al^d x =x \theta_\al^d$. Hence $\theta_\al^d $ is
central in $H_\al$. Finally, let us show Part~(d). Using Lemma~\ref{B2}(d) and \eqref{twist1}, we have that $u_\al
= \pa (u_\al) = \ta^{-1} u_\al \ta$, and so $\ta u_\al = u_\al \ta$.
\end{proof}

For any $\al \in \pi$, we set $$G_\al=\ta u_\al = u_\al \ta \in H_\al.$$

\begin{lemma}\label{zig}
Let $H=(\{H_\al\},\Delta,\epsilon,S,\varphi,R,\theta)$ be a ribbon Hopf \p-coalgebra. Then
\begin{enumerate}
\renewcommand{\labelenumi}{{\rm (\alph{enumi})}}
  \item $G=(G_\al)_{\al \in \pi}$ is a \p-grouplike element of $H$;
  \item $\varphi_\be (G_\al ) = G_{\be \al \bei}$ for all $\al,\be \in\pi$;
  \item $S_\al(G_\al)=G_\ali^{-1}$  for all $\al \in \pi$;
  \item $\ta^{-2}=c_\al$ for all $\al \in \pi$, where $c_\al=S_\ali(u_\ali)u_\al
        =u_\al S_\ali(u_\ali)$ as in Lemma~\ref{B2}(e);
  \item $S_\al (u_\al)=G_\ali^{-1} u_\ali G_\ali^{-1} $  for all $\al \in \pi$;
  \item $S_\ali S_\al (x) = G_\al x
        G_\al^{-1}$  for all $\al \in \pi$ and $x \in H_\al$.
\end{enumerate}
\end{lemma}

\begin{proof}
Let us show Part (a). Firstly $\epsilon(G_1) = \epsilon(\theta_1 u_1) = \epsilon(\theta_1) \epsilon(u_1) = 1$ by
Lemmas~\ref{B2}(g) and \ref{C1}(b). Secondly, for any $\al,\be \in \pi$, using \eqref{twist4} and
Lemma~\ref{B2}(f),
\begin{eqnarray*}
   \cp{\al}{\be}(G_{\al\be}) & = & \cp{\al}{\be}(\theta_{\al\be} u_{\al\be}) \\
      & = & \cp{\al}{\be}(\theta_{\al\be}) \cdot \cp{\al}{\be}(u_{\al\be})\\
      & = & (\ta \otimes \tb) \cdot
            [\sigma_{\be,\al}( (\varphi_\ali \otimes \id_{H_\al} )(R_{\al \be \ali ,\al}))
             \cdot R_{\al,\be}] \\
      &   & \quad \quad \cdot [\sigma_{\be,\al}( (\varphi_\ali \otimes \id_{H_\al} )
            (R_{\al \be \ali ,\al}))
            \cdot R_{\al,\be}]^{-1} \cdot (u_\al \otimes u_\be)\\
      & = & G_\al \otimes G_\be.
\end{eqnarray*} Thus $G=(G_\al)_{\al \in \pi} \in G(H)$.
Part (b) follows directly from Lemma~\ref{B2}(d) and \eqref{twist3}, and Part (c) from the fact that $G$ is a
\p-grouplike element. By Part (c) and \eqref{twist2}, $\ta^{-2}=u_\al G_\al^{-1} \ta^{-1} = u_\al S_\ali(G_\ali)
\ta^{-1} = u_\al S_\ali(\tai u_\ali)  \ta^{-1} =c_\al$ and so Part (d) is established. Let us show Part (e). By
\eqref{twist2} and Part~(c), $ G_\ali^{-1} u_\ali=S_\al(\ta^{-1}) = S_\al(G_\al^{-1}u_\al) =S_\al(u_\al)
S_\al(G_\al)^{-1} =S_\al(u_\al) G_\ali $. Therefore $S_\al(u_\al)=G_\ali^{-1} u_\ali G_\ali^{-1}$. Finally, to
show Part (f), let $x \in H_\al$. Then, using Lemmas~\ref{B2}(b) and \ref{C1}(a), $S_\ali S_\al(x) = u_\al \pai(x)
u_\al^{-1} =u_\al \ta x \ta^{-1} u_\al^{-1}= G_\al x G_\al^{-1}$. This completes the proof of the lemma.
\end{proof}

\subsubsection{The coopposite Hopf $\pi$-coalgebra}\label{cooprib}
Let $H$ be a ribbon Hopf \p-coalgebra with twist $\theta=\{\theta_\al\}_{\al \in \pi}$. The coopposite
quasitriangular Hopf \p-coalgebra $H^\cop$ (see \S\ref{coopquasi}) is ribbon with twist
$\theta_\al^\cop=\theta_\ali^{-1}$.

\subsubsection{The mirror Hopf $\pi$-coalgebra}\label{mirrorrib}
Let $H$ be a ribbon Hopf \p-coalgebra with twist $\theta=\{\theta_\al\}_{\al \in \pi}$. Following \cite[\S
11.6]{Tur1}, the mirror quasitriangular Hopf \p-coalgebra $\overline{H}$ (see \S\ref{mirrorquasi}) is ribbon with
twist $\overline{\theta}_\al=\theta_\ali^{-1}$.

\subsection{The distinguished $\pi$-grouplike element from the $R$-matrix}
In this subsection, we show that the  distinguished \p-grouplike element of a finite dimensional quasitriangular
Hopf \p-coalgebra can be computed by using the $R$-matrix. This generalizes \cite[Theorem 2]{Rad2}.
\begin{theorem}\label{calcg}
Let $H$ be a finite dimensional quasitriangular Hopf \p-coalgebra. Let $g=(g_\al)_{\al \in \pi}$ be the
distinguished \p-grouplike element of $H$, $\nu$ be the distinguished grouplike element of $H_1^*$,
$\ell=(\ell_\al)_{\al \in \pi} \in G(H)$ be as in Corollary~\ref{ell}, and $\widehat{\varphi}$ be as in
Corollary~\ref{phitopi}. We define $h_\al=(\id_{H_\al} \otimes \nu)(R_{\al,1})$ for any $\al \in \pi$. Then
\begin{enumerate}
\renewcommand{\labelenumi}{{\rm (\alph{enumi})}}
  \item $h=(h_\al)_{\al \in \pi}$ is a \p-grouplike
        element of $H$;
  \item g = $\widehat{\varphi}^{-1}\ell h$ in $G(H)$, i.e.,
        $g_\al= \widehat{\varphi}(\al)^{-1}\ell_\al h_\al$ for all $\al \in \pi$.
\end{enumerate}
\end{theorem}
\begin{proof}
We adapt the technique used in the proof of \cite[Theorem 2]{Rad2}. Let us first show Part (a). For any $\al,\be
\in \pi$, using \eqref{quasicomult}, the multiplicativity of $\nu$, and Lemma~\ref{marre}(b), we have that
\begin{eqnarray*}
 \cp{\al}{\be}(h_{\al \be})
   & = & \cp{\al}{\be} (\id_{H_{\al \be}} \otimes \nu)(R_{\al
         \be,1})\\
   & = & (\id_{H_\al} \otimes \id_{H_\be} \otimes \nu)
         ([(\id_{H_\al} \otimes \pbi)(R_{\al,1})]_{1 \be 3} \cdot
         (R_{\be,1})_{\al 23}) \\
   & = & ((\id_{H_\al} \otimes \nu \pbi)(R_{\al,1}) \otimes
         1_\be) \cdot (1_\al \otimes (\id_{H_\be} \otimes
         \nu)(R_{\be,1})) \\
   & = & ((\id_{H_\al} \otimes \nu)(R_{\al,1}) \otimes
         1_\be) \cdot (1_\al \otimes h_\be)\\
   & = & h_\al \otimes h_\be.
\end{eqnarray*}
Moreover, using Lemma~\ref{B1}(a), $\epsilon(h_1)=(\epsilon \otimes \nu)(R_{1,1})=\nu(1_1)=1$. Thus $h \in G(H)$.

To show Part (b), let $\al \in \pi$ and $\Lambda$ be a non-zero left integral for $H_1$. We first show that, for
any $x \in H_\ali$,
\begin{equation}\label{fini1}
\Lambda_{(1,\al)} \otimes x \Lambda_{(2,\ali)} = S_\ali(x) \Lambda_{(1,\al)} \otimes \Lambda_{(2,\ali)}
\end{equation}
and
\begin{equation}\label{fini2}
\Lambda_{(1,\ali)}x \otimes \Lambda_{(2,\al)} = \Lambda_{(1,\ali)} \otimes \Lambda_{(2,\al)} S_\ali(x
\leftharpoonup\nu).
\end{equation}
Indeed
\begin{eqnarray*}
\lefteqn{\Lambda_{(1,\al)} \otimes x \Lambda_{(2,\ali)}}\\
  & = & \epsilon(\xuu) \, \Lambda_{(1,\al)} \otimes \xdai
        \Lambda_{(2,\ali)}
        \text{ \quad by \eqref{counit}}\\
  & = & S_\ali(\xuai) \xda \Lambda_{(1,\al)} \otimes  \xtai
        \Lambda_{(2,\ali)}
        \text{ \quad by \eqref{antipode}}\\
  & = & S_\ali(\xuai) (\xdu \Lambda)_{(1,\al)} \otimes  (\xdu
        \Lambda)_{(2,\ali)} \text{ \quad by \eqref{hopfmor},}
\end{eqnarray*}
and so, since $\Lambda$ is a left integral for $H_1$,
\begin{eqnarray*}
\Lambda_{(1,\al)} \otimes x \Lambda_{(2,\ali)}
  & = & S_\ali(\xuai \epsilon(\xdu)) \Lambda_{(1,\al)} \otimes
        \Lambda_{(2,\ali)}  \\
  & = & S_\ali(x) \Lambda_{(1,\al)} \otimes
        \Lambda_{(2,\ali)}\text{ \quad by \eqref{counit}.}
\end{eqnarray*}
Similarly,
\begin{eqnarray*}
\lefteqn{\Lambda_{(1,\ali)} x \otimes \Lambda_{(2,\ali)}}\\
  & = & \Lambda_{(1,\ali)} \xuai \otimes
        \Lambda_{(2,\al)} \epsilon(\xdu)
        \text{ \quad by \eqref{counit}}\\
  & = & \Lambda_{(1,\ali)} \xuai \otimes
        \Lambda_{(2,\al)} \xda S_\ali(\xtai)
        \text{ \quad by \eqref{antipode}}\\
  & = & (\Lambda \xuu)_{(1,\ali)} \otimes
        (\Lambda \xuu)_{(2,\al)} S_\ali(\xdai) \text{ \quad by \eqref{hopfmor},}
\end{eqnarray*}
and so, since $\Lambda$ is a left integral for $H_1$,
\begin{eqnarray*}
 \Lambda_{(1,\ali)} x \otimes \Lambda_{(2,\ali)}
  & = & \Lambda_{(1,\ali)} \otimes
        \Lambda_{(2,\al)} S_\ali(\nu(\xuu) \, \xdai) \\
  & = & \Lambda_{(1,\ali)} \otimes
        \Lambda_{(2,\al)} S_\ali(x \leftharpoonup \nu).
\end{eqnarray*}

Write $R_{\al,\ali}=a_\al \otimes b_\ali$. Recall that $u_\al=S_\ali\pa(b_\ali) a_\al$. By Lemma~\ref{B1}(c) and
\eqref{quasiphi}, $R_{\ali,\al}=S_\al(a_\al) \otimes \pa S_\ali(b_\ali)$. Thus $u_\ali=S_\al S_\ali(b_\ali)
S_\al(a_\al)$ and so, using Lemma~\ref{B2}(b) and (d), $S_\ali(u_\ali)=S_\al^{-1}(u_\ali)=a_\al S_\ali(b_\ali)$.
Then
\begin{eqnarray*}
 \lefteqn{
   \Lambda_{(2,\al)} S_\ali(\pa(b_\ali)\leftharpoonup \nu)
   a_\al \otimes \Lambda_{(1,\ali)}
         }\\
   & = & \Lambda_{(2,\al)} a_\al \otimes \Lambda_{(1,\ali)}
         \pa(b_\ali)
         \text{ \quad by \eqref{fini2} for $x=\pa(b_\ali)$}\\
   & = & (\id_{H_\al} \otimes \pa) (\Lambda_{(2,\al)} a_\al
         \otimes \pai(\Lambda_{(1,\ali)}) b_\ali)\\
   & = & (\id_{H_\al} \otimes \pa)( a_\al \Lambda_{(1,\al)}
         \otimes b_\ali \Lambda_{(2,\ali)})
         \text{ \quad by \eqref{quasicop}}\\
   & = & (\id_{H_\al} \otimes \pa) (a_\al S_\ali(b_\ali)
         \Lambda_{(1,\al)} \otimes \Lambda_{(2,\ali)})
         \text{ \quad by \eqref{fini1} for $x=b_\ali$}\\
   & = & S_\ali(u_\ali) \Lambda_{(1,\al)} \otimes
         \pa(\Lambda_{(2,\ali)}) \\
   & = & (\pai \otimes \id_{H_\ali})
         (\pa S_\ali(u_\ali) \pa(\Lambda_{(1,\al)}) \otimes
         \pa(\Lambda_{(2,\ali)})) \text{ \quad by \eqref{phiaction}}\\
   & = & (\pai \otimes \id_{H_\ali})
         (\pa S_\ali(u_\ali) \pa(\Lambda)_{(1,\al)} \otimes
         \pa(\Lambda)_{(2,\ali)})
         \text{ \quad by \eqref{phicomult}.}
\end{eqnarray*}
Now $\pa(\Lambda)=\widehat{\varphi}(\al)\,\Lambda$ by Lemma~\ref{marre}(a) and $$ \Lambda_{(1,\al)} \otimes
\Lambda_{(2,\ali)} = S_\ali S_\al(\Lambda_{(2,\al)}) g_\al \otimes \Lambda_{(1,\ali)} $$ by
Corollary~\ref{chatquipue}. Therefore
\begin{eqnarray*}
 \lefteqn{
   \Lambda_{(2,\al)} S_\ali(\pa(b_\ali)\leftharpoonup \nu)
   a_\al \otimes \Lambda_{(1,\ali)} }\\
   & = & \widehat{\varphi}(\al)\,(\pai \otimes \id_{H_\ali}) (\pa S_\ali(u_\ali)
         S_\ali S_\al(\Lambda_{(2,\al)}) g_\al \otimes
         \Lambda_{(1,\ali)}) \\
   & = & \widehat{\varphi}(\al)\,S_\ali(u_\ali) \pai S_\ali
         S_\al(\Lambda_{(2,\al)})  \pai(g_\al) \otimes
         \Lambda_{(1,\ali)} \\
   & = & \widehat{\varphi}(\al)\,S_\ali(u_\ali) \pai S_\ali
         S_\al(\Lambda_{(2,\al)})  g_\al \otimes
         \Lambda_{(1,\ali)}
         \text{ \quad by Lemma~\ref{marre}(c).}
\end{eqnarray*}
Let  $\lambda=(\lambda_\ga )_{\ga \in \pi}$ be left \p-integral for $H$ such that $\lu(\Lambda)=1$ (see the proof
of Corollary~\ref{chatquipue}). Applying $(\id_{H_\al} \otimes \lai)$ on both sides of the last equality, we get
\begin{multline*}
 \lai (\Lambda_{(1,\ali)}) \, \Lambda_{(2,\al)} S_\ali(\pa(b_\ali)\leftharpoonup \nu)
 a_\al\\
  = \widehat{\varphi}(\al)\, S_\ali(u_\ali) \pai
 S_\ali S_\al(\lai (\Lambda_{(1,\ali)}) \, \Lambda_{(2,\al)})
 g_\al,
\end{multline*}
and so, since $\lai(\Lambda_{(1,\ali)}) \, \Lambda_{(2,\al)}=\lu (\Lambda) 1_\al=1_\al$,
\begin{equation}\label{fini3}
S_\ali(\pa(b_\ali)\leftharpoonup \nu) a_\al = \widehat{\varphi}(\al)\, S_\ali(u_\ali) g_\al.
\end{equation}
Write $R_{\al,1}=c_\al \otimes d_1$ so that $h_\al=\nu(d_1) c_\al$. Since, by \eqref{phicomult} and
Lemma~\ref{marre}(b), $\pa(x) \leftharpoonup \nu = \pa(x \leftharpoonup \nu)$ for all $x \in H_\ali$, we have that
\begin{eqnarray*}
a_\al \otimes \pa(b_\ali) \leftharpoonup \nu
  & = & a_\al \otimes \pa(b_\ali \leftharpoonup \nu) \\
  & = & (\id_{H_\al} \otimes \nu \otimes \pa)
        (\id_{H_\al} \otimes \cp{1}{\ali})(R_{\al,\ali}) \\
  & = & (\id_{H_\al} \otimes \nu \otimes \pa)
        ((R_{\al,\ali})_{1 {\scriptscriptstyle 1_\pi} 3} \cdot (R_{\al,1})_{12 \ali})
        \text{ \quad by \eqref{quasicomult}}\\
  & = & a_\al \nu(d_1) \,c_\al \otimes \pa(b_\ali)\\
  & = & a_\al h_\al \otimes  \pa(b_\ali).
\end{eqnarray*}
Therefore $S_\ali(\pa(b_\ali) \leftharpoonup \nu) a_\al= S_\ali(\pa(b_\ali)) a_\al h_\al=u_\al h_\al$. Finally,
comparing with \eqref{fini3}, we get $\widehat{\varphi}(\al)\, S_\ali(u_\ali)g_\al = u_\al h_\al$. Hence $g_\al
=\widehat{\varphi}(\al)^{-1}\ell_\al h_\al$, since $\ell_\al=S_\ali(u_\ali)^{-1}u_\al$. This finishes the proof of
the theorem.
\end{proof}

\begin{corollary}\label{geqgcarre}
Let $H$ be a finite dimensional ribbon Hopf \p-coalgebra. Let $g=(g_\al)_{\al \in \pi}$ be the distinguished
\p-grouplike element of $H$, $h=(h_\al)_{\al \in \pi} \in G(H)$ as in Theorem~\ref{calcg}, $G=(G_\al)_{\al \in
\pi} \in G(H)$ as in Lemma~\ref{zig}, and $\widehat{\varphi}$ as in Corollary~\ref{phitopi}. Then
$\widehat{\varphi}\,g=G^2 h$ in $G(H)$, i.e., $\widehat{\varphi}(\al)g_\al=G^2_\al h_\al$ for all $\al \in \pi$.
\end{corollary}
\begin{proof}
For any $\al \in \pi$, $\widehat{\varphi}(\al) g_\al= S_\ali(u_\ali)^{-1} u_\al h_\al=\ta^2 u_\al^2 h_\al=G_\al^2
h_\al $ by Theorem~\ref{calcg}(b) and Lemma~\ref{zig}(d).
\end{proof}

\section{Existence of $\pi$-traces}\label{s:pitrace}
In this section, we introduce the notion of a \p-trace for a crossed Hopf \p-coalgebra and we show the existence
of \p-traces for a finite dimensional unimodular Hopf \p-coalgebra whose crossing $\varphi$ verifies that
$\widehat{\varphi}=1$. Moreover, we give sufficient conditions for the homomorphism $\widehat{\varphi}$ to be
trivial.

\subsection{Unimodular Hopf $\pi$-coalgebras}
A Hopf \p-coalgebra $H=\{H_\al\}_{\al \in \pi}$ is said to be \emph{unimodular} if the Hopf algebra $H_1$ is
unimodular (it means that the spaces of left and right integrals for $H_1$ coincide). If $H_1$ is finite
dimensional, then $H$ is unimodular if and only if $\nu=\epsilon$, where $\nu$ is the distinguished grouplike
element of $H_1^*$.

If $\pi$ is finite, then a left (resp.\@ right) integral for the Hopf algebra $\tilde{H}=\oplus_{\al \in \pi}
H_\al$ (see \S\ref{pifinite}) must belong to $H_1$, and so the spaces of left (resp.\@ right) integrals for
$\tilde{H}$ and $H_1$ coincide. Hence, when $\pi$ is finite, $H$ is unimodular if and only if $\tilde{H}$ is
unimodular.

One can remark that a semisimple finite dimensional Hopf $\pi$-coalgebra $H=\{H_\al\}_{\al \in \pi}$ is unimodular
(since the finite dimensional Hopf algebra $H_1$ is semisimple and so unimodular). Note that a cosemisimple Hopf
\p-coalgebra is not necessarily unimodular.

\subsection{$\pi$-traces}\label{pitra}

Let $H=(\{H_\al\},\Delta,\epsilon,S,\varphi)$ be a crossed Hopf \p-coalgebra. A \emph{\p-trace} for $H$ is a
family of $\Bbbk$-linear forms $\tr=(\tr_\al)_{\al \in \pi} \in \Pi_{\al \in \pi} H_\al^*$ such that, for any
$\al,\be \in \pi$ and $x,y \in H_\al$,
\begin{defi}
  \item \label{trace1} $\tr_\al(xy)=\tr_\al(yx)$;
  \item \label{trace2} $\tr_\ali(S_\al(x))=\tr_\al(x)$;
  \item \label{trace3} $\tr_{\be \al \bei}(\pb(x))=\tr_\al(x)$.
\end{defi}

This notion is motivated mainly by topological purposes: \p-traces are used in \cite{Vir} to construct
Hennings-like invariants (see \cite{He2,KR1}) of principal \p-bundles over link complements and over
$3$-manifolds.

Note that $\tr_1$ is a (usual) trace for the Hopf algebra $H_1$, invariant under the action $\varphi$ of $\pi$.

In the next lemma, generalizing \cite[Proposition~4.2]{He2}, we give a characterization of the \p-traces.
\begin{lemma}\label{int1}
Let $H=\{H_\al\}_{\al \in \pi}$ be a finite dimensional unimodular ribbon Hopf \p-coalgebra with crossing
$\varphi$. Let $\lambda=( \la)_{\al \in \pi}$ be a non-zero right \p-integral for $H$, $G=(G_\al)_{\al \in \pi}
\in G(H)$ be as in Lemma~\ref{zig}, and $\widehat{\varphi}$ be as in Corollary~\ref{phitopi}. Let
$\tr=(\tr_\al)_{\al \in \pi} \in \Pi_{\al \in \pi} H_\al^*$. Then $\tr$ is a \p-trace for $H$ if and only if there
exists a family $z=(z_\al)_{\al \in \pi} \in \Pi_{\al \in \pi} H_\al$ satisfying, for all $\al,\be \in \pi$,
\begin{enumerate}
\renewcommand{\labelenumi}{{\rm (\alph{enumi})}}
  \item $\tr_\al(x)=\la(G_\al z_\al x)$ for all $x \in H_\al$;
  \item $z_\al$ is central in $H_\al$;
  \item $S_\al(z_\al)=\widehat{\varphi}(\al)^{-1} z_\ali$;
  \item $\pb(z_\al)=\widehat{\varphi}(\be) z_{\be \al \bei}$.
\end{enumerate}
\end{lemma}
\begin{proof}
We first show that, for all $\al \in \pi$ and $x,y \in H_\al$,
\begin{equation}\label{symtra}
\la(G_\al x y)=\la(G_\al y x),
\end{equation}
and
\begin{equation}\label{dertra}
\widehat{\varphi}(\al) \lai(S_\al(x))=\la(G^2_\al x).
\end{equation}
Indeed, let $\nu$ be the distinguished grouplike element of $H_1^*$. Since $\nu=\epsilon$ ($H$ is unimodular),
Theorem \ref{intfortrace}(a) gives that $\la(G_\al xy)=\la(S_\ali S_\al(y)\, G_\al x)$. Now, by
Lemma~\ref{zig}(f), $S_\ali S_\al (y) = G_\al y G_\al^{-1}$. Thus $\la(G_\al x y)=\la(G_\al y x)$ and
\eqref{symtra} is proven. Moreover, Corollary~\ref{geqgcarre} gives that $\widehat{\varphi}(\al) g_\al=G^2_\al
h_\al$, where $g=(g_\al)_{\al \in \pi}$ is the distinguished \p-grouplike element of $H$ and $h_\al=(\id_{H_\al}
\otimes \nu)(R_{\al,1})$. Since $\nu=\epsilon$ and by Lemma~\ref{B1}(a), $h_\al=(\id_{H_\al} \otimes
\epsilon)(R_{\al,1})=1_\al$. Thus $\widehat{\varphi}(\al) g_\al=G^2_\al$. Now $\lai(S_\al(x))=\la(g_\al x)$ by
Theorem~\ref{intfortrace}(c). Hence $\widehat{\varphi}(\al)\lai(S_\al(x))=\la(G^2_\al x)$ and \eqref{dertra} is
proven.

Let us suppose that there exists $z=(z_\al)_{\al \in \pi} \in \Pi_{\al \in \pi} H_\al$ verifying
Conditions~(a)-(d). For any $\al, \be \in \pi$ and $x,y \in H_\al$,
\begin{eqnarray*}
\tr_\al(xy) & = & \la(G_\al z_\al x y) \text{ \quad by Condition (a)}\\
  & = & \la(G_\al y z_\al x) \text{ \quad by \eqref{symtra}}\\
  & = & \la( G_\al z_\al y x) \text{ \quad since $ z_\al$ is central}\\
  & = & \tr_\al(y x) \text{ \quad by Condition (a),}
\end{eqnarray*}
\begin{eqnarray*}
 \lefteqn{ \tr_\ali(S_\al(x))}\\
        & = & \lai (G_\ali z_\ali S_\al(x)) \\
        & = & \widehat{\varphi}(\al) \lai (S_\al(G_\al^{-1}) S_\al(z_\al) S_\al(x)) \text{ \quad by Condition (c) and
              Lemma~\ref{zig}(c)}\\
        & = &  \widehat{\varphi}(\al) \lai (S_\al(x z_\al G_\al^{-1})) \text{ \quad by Lemma~\ref{antipodepptes}(a)}\\
        & = & \la (G_\al^{\,2} x z_\al G_\al^{-1}) \text{ \quad by
              \eqref{dertra}}\\
        & = & \la (G_\al z_\al G_\al x G_\al^{-1}) \text{ \quad since $z_\al$ is central}\\
        & = & \tr_\al (G_\al x G_\al^{-1})\\
        & = & \tr_\al (x) \text{ \quad since $\tr_\al$
              is symmetric,}
\end{eqnarray*}
and
\begin{eqnarray*}
  \lefteqn{\tr_{\be \al \bei}(\pb (x))}\\
      & = & \lbabi (G_{\be \al \bei} z_{\be \al \bei} \pb (x))\\
      & = & \widehat{\varphi}(\be)^{-1} \lbabi (\pb (G_\al) \pb(z_\al) \pb (x)) \text{ \quad by Condition (d) and
            Lemma~\ref{zig}(b)}\\
      & = & \widehat{\varphi}(\be)^{-1} \lbabi (\pb (G_\al z_\al x))\\
      & = & \widehat{\varphi}(\be)^{-1} \widehat{\varphi}(\be)\la (G_\al z_\al x)
            \text{ \quad by Corollary~\ref{phitopi}}\\
      & = & \tr_\al (x).
\end{eqnarray*}
Hence $\tr$ is a \p-trace.

Conversely, suppose that $\tr$ is a \p-trace. Recall that $H_\al^*$ is a right $H_\al$-module for the action
defined, for all $f \in H_\al^*$ and $a,x \in H_\al$, by
       $$
       (f\leftharpoonup a)(x)= f(a x).
       $$
By Corollary~\ref{corinteg}(b), $(H_\al^*,\leftharpoonup)$ is free, its rank is $1$ (resp.\@ $0$) if $H_\al \neq
0$ (resp.\@ $H_\al=0$), and $\la$ is a basis vector for $H^*_\al$. Thus, for any $\al \in \pi$, there exists
$w_\al \in H_\al$ such that $\tr_\al = \la \leftharpoonup w_\al$. Set $z_\al =G_\al^{-1}w_\al$. Let us verify that
the family $z=(z_\al)_{\al \in \pi}$ verify Conditions (a)-(d). By the definition of $z_\al$, Condition (a) is
clearly verified. Let $\al \in \pi$ and $x \in H_\al$. For any $y \in H_\al$,
\begin{eqnarray*}
(\la \leftharpoonup G_\al z_\al x)(y)
   & = & \la(G_\al z_\al x y) \\
   & = & \tr_\al(xy) \\
   & = & \tr_\al(yx) \text{ \quad by \eqref{trace1}}\\
   & = & \la(G_\al z_\al y x) \\
   & = & \la(G_\al x z_\al y) \text{ \quad by \eqref{symtra}}\\
   & = & (\la \leftharpoonup G_\al x z_\al)(y).
\end{eqnarray*}
Therefore $\la \leftharpoonup G_\al z_\al x=\la \leftharpoonup G_\al x z_\al$. Hence $G_\al z_\al x= x z_\al$
(since $\la$ is a basis vector for $(H^*_\al,\leftharpoonup)$) and so $z_\al x = x z_\al$. Condition (b) is then
verified. Let $\al \in \pi$. For any $x \in H_\al$,
\begin{eqnarray*}
\lefteqn{(\lai \leftharpoonup G_\ali S_\al(z_\al))(x)}\\
  & = & \lai( G_\ali S_\al(z_\al)x) \\
  & = & \lai( S_\al (S_\al^{-1}(x)z_\al G_\al^{-1})) \text{ \quad by Lemmas~\ref{antipodepptes}(a) and \ref{zig}(c)}\\
  & = & \widehat{\varphi}(\al)^{-1} \la (G_\al^2 S_\al^{-1}(x)z_\al G_\al^{-1})  \text{ \quad by \eqref{dertra}}\\
  & = & \widehat{\varphi}(\al)^{-1} \la (G_\al z_\al S_\al^{-1}(x)) \text{ \quad by \eqref{symtra} and
        since $z_\al$ is central}\\
  & = & \widehat{\varphi}(\al)^{-1} \tr_\al (S_\al^{-1}(x)) \\
  & = & \widehat{\varphi}(\al)^{-1} \tr_\ali (x) \text{ \quad by \eqref{trace2}}\\
  & = & (\lai \leftharpoonup G_\ali \widehat{\varphi}(\al)^{-1} z_\ali)(x).
\end{eqnarray*}
We conclude as above that $S_\al(z_\al)= \widehat{\varphi}(\al)^{-1} z_\ali$, and so Condition (c) is satisfied.
Finally, let $\al,\be \in \pi$. For any $x \in H_\al$,
\begin{eqnarray*}
 \lefteqn{(\la \leftharpoonup G_\al \widehat{\varphi}(\be) \pbi(z_{\be \al \bei}))(x)}\\
   & = & \widehat{\varphi}(\be) \la ( G_\al \pbi(z_{\be \al \bei}) x) \\
   & = & \lambda_{\be \al \bei} \pb(G_\al \pbi(z_{\be \al \bei}) x)
         \text{ \quad by Corollary~\ref{phitopi}}\\
   & = & \lambda_{\be \al \bei}( G_{\be \al \bei} z_{\be \al \bei} \pb(x)) \text{ \quad by Lemma~\ref{zig}(b)} \\
   & = & \tr_{\be \al \bei}(\pb(x)) \\
   & = & \tr_\al(x) \text{ \quad by \eqref{trace3}}\\
   & = & (\la \leftharpoonup G_\al z_\al)(x).
\end{eqnarray*}
Thus $\widehat{\varphi}(\be) \pbi(z_{\be \al \bei})=z_\al$ and so $\pb(z_\al)= \widehat{\varphi}(\be) z_{\be \al
\bei}$. Hence Condition (d) is verified and the lemma is proven.
\end{proof}

In the setting of Lemma~\ref{int1}, constructing a \p-trace from a right \p-integral $\lambda=(\la)_{\al \in \pi}$ 
reduces to finding a family $z=(z_\al)_{\al \in \pi}$ which satisfies Conditions (b)-(d) of Lemma~\ref{int1}. Let
us give two possible choices of the family $z$.

Let $\Lambda$ be a left integral for $H_1$ such that $\lu(\Lambda)=1$. Set $z_1=\Lambda$ and $z_\al=0$ if $\al
\neq 1$. This family $z=(z_\al)_{\al \in \pi}$ verifies Conditions (b)-(d) (since $H$ is unimodular and so
$\Lambda$ is central and $S_1(\Lambda)=\Lambda$, and by Lemma~\ref{marre}(a)). The \p-trace obtained is given by
$\tr_1=\epsilon$ and $\tr_\al=0$ if $\al \neq 1$.

If the homomorphism $\widehat{\varphi}$ of Corollary~\ref{phitopi} is trivial (that is $\widehat{\varphi}(\al)=1$
for all $\al \in \pi$), then another possible choice is $z_\al=1_\al$. In the two next lemmas, we give sufficient
conditions for the homomorphism $\widehat{\varphi}$ to be trivial.

\begin{lemma}\label{phitrivial}
Let $H$ be a finite dimensional crossed Hopf \p-coalgebra with crossing $\varphi$. If $H$ is semisimple or
cosemisimple or if $\varphi_{\be|H_1}=\id_{H_1}$ for all $\be \in \pi$, then $\widehat{\varphi}=1$.
\end{lemma}
\begin{proof}
Let $\be \in \pi$. If $H$ is semisimple, then $H_1$ is semisimple and thus there exists a left integral $\Lambda$
for $H_1$ such that $\epsilon(\Lambda)=1$ (by \cite[Theorem~5.1.8]{sweed}). Now
$\pb(\Lambda)=\widehat{\varphi}(\be) \Lambda$ by Lemma~\ref{marre}(a). Therefore, using \eqref{phicounit},
$\widehat{\varphi}(\be) = \widehat{\varphi}(\be) \epsilon(\Lambda)=\epsilon(\widehat{\varphi}(\be) \Lambda)=
\epsilon\pb (\Lambda)= \epsilon(\Lambda)=1$. Suppose now that $H$ is cosemisimple. By Theorem~\ref{cosemicrit},
there exists a right \p-integral $\lambda=( \la )_{\al\in \pi}$ for $H$ such that $\lu(1_1)=1$. Then
$\widehat{\varphi}(\be) = \widehat{\varphi}(\be) \lu(1_1)=\lu(\pb(1_1))=\lu(1_1)=1$. Suppose finally that
$\varphi_{\be|H_1}=\id_{H_1}$. Let $\lambda=( \la)_{\al \in \pi}$ be a non-zero right \p-integral for $H$. Then
$\widehat{\varphi}(\be) \lu=\lu \varphi_{\be|H_1} =\lu$ and thus $\widehat{\varphi}(\be)=1$ (since $\lu \neq 0$ by
Lemma~\ref{intnonzero}).
\end{proof}

\begin{lemma}\label{thetaphitriv}
Let $H$ be a finite dimensional ribbon Hopf \p-coalgebra with crossing $\varphi$ and twist
$\theta=\{\theta_\al\}_{\al \in \pi}$. Let $\lambda=(\la)_{\al \in \pi}$ be a right \p-integral for $H$. If
$\lu(\theta_1) \neq 0$, then $\widehat{\varphi}=1$.
\end{lemma}
\begin{proof}
Let $\be \in \pi$. By (\ref{ribpico}.c) and Corollary~\ref{phitopi}, $\lu(\theta_1)=\lu(\pb(\theta_1))=
\widehat{\varphi}(\be) \lu(\theta_1)$. Thus, since $\lu(\theta_1) \neq 0$, $\widehat{\varphi}(\be)=1$.
\end{proof}

We conclude with the following theorem, which follows directly from  Lemma~\ref{int1} (by choosing $z_\al=1_\al$
for all $\al \in \pi$) and Lemmas~\ref{phitrivial} and \ref{thetaphitriv}.

\begin{theorem}\label{traceunimod}
Let $H$ be a finite dimensional unimodular ribbon Hopf \p-coalgebra with crossing $\varphi$. Let
$\lambda=(\la)_{\al \in \pi}$ be a right \p-integral for $H$ and $G=(G_\al)_{\al \in \pi} \in G(H)$ be as in
Lemma~\ref{zig}. Suppose that at least one of the following conditions is verified:
\begin{enumerate}
\renewcommand{\labelenumi}{{\rm (\alph{enumi})}}
  \item $H$ is semisimple;
  \item $H$ is cosemisimple;
  \item $\lu(\theta_1) \neq 0$;
  \item $\varphi_{\be|H_1}=\id_{H_1}$ for all $\be \in \pi$.
\end{enumerate}
Then $\tr=(\tr_\al)_{\al \in \pi}$, defined by $\tr_\al(x)=\la(G_\al x)$ for all $\al \in \pi$ and $x \in H_\al$,
is a \p-trace for $H$.
\end{theorem}

\bibliographystyle{amsalpha}
\bibliography{alg-bib}
\end{document}